\documentclass[12pt,a4paper,twoside,final,notitlepage, leqno]{article}
\usepackage[english]{babel}
\usepackage[T1]{fontenc}  
\usepackage{graphicx}
\usepackage{amsmath,amsthm,epsfig,amsfonts,bbm}  

\newcommand{\pequationdeb}{$$ \left\{ \begin{minipage}[c]{130mm}}
\newcommand{\pequationfin}{\end{minipage}
                           \right. $$}

\def \smb {{\scriptstyle \bullet }}

\newcommand{\beq}     {\begin{equation}}
\newcommand{\enq}     {\end{equation}}
\newcommand{\be}    {\begin{enumerate}}
\newcommand{\ee}    {\end{enumerate}}

\newcommand{\Bb}

\newcommand{\R}{\mathbb{R}}

\newcommand{\N}{\mathbb{N}}

\newcommand{\pet}   { & \scriptscriptstyle \mathbf \!\!\!\!\!\! }

\def\resume{\if@twocolumn
\section*{R\'esum\'e}
\else \small
\quotation{\bf \it R\'esum\'e \rule[1mm]{1.5mm}{0.2mm}\vspace{0pt}}
\fi}
\def\endresume{\if@twocolumn\else\endquotation\fi}
\def\abstract{\if@twocolumn
\noindent\section*{{\bf Abstract}}
\else \small
\quotation{\noindent \bf {Abstract.} \rule[1mm]{1.5mm}{0.2mm}\vspace{0pt}}
\fi}
\def\endabstract{\if@twocolumn\else\endquotation\fi}

\begin{document} 
\title{\bf \Large   Towards higher order   lattice Boltzmann schemes  ~\\~  }

\author { { \large  Fran\c{c}ois Dubois~$^{ab}$ and   Pierre Lallemand~$^{c}$}     \\ ~\\
{\it \small  $^a$  Department of Mathematics, Universit\'e Paris-Sud,} \\
{\it \small B\^at. 425, F-91405 Orsay Cedex, France} \\ 
{\it  \small $^b$   Conservatoire National des Arts et M\'etiers,}  \\
{\it \small Department of  Mathematics and EA3196,  Paris, France.} \\ 
{\it \small francois.dubois@math.u-psud.fr} \\ 
{\it \small  $^c$  Retired from Centre National de la Recherche Scientifique, Paris.} \\  
{\it \small  pierre.lal@free.fr}  \\  ~\\~\\ }  
\date{ 13 June  2009~\protect\footnote{~Published 12 June 2009, {\it Journal of Statistical
    Mechanics: theory and experiment},  P06006, 2009,  
doi: 10.1088/1742-5468/2009/06/P06006,
http://www.iop.org/EJ/abstract/1742-5468/2009/06/P06006/, arXiv:0811.0599.
Edition 14 December 2009.}}

\maketitle

\begin{abstract}
In this contribution we extend the Taylor expansion method proposed previously by one of
us and establish  equivalent partial  differential equations of the     
lattice Boltzmann scheme proposed by   d'Humi\`eres \cite{DDH92}    
at an arbitrary order of accuracy. 
We derive formally the  associated dynamical equations for classical thermal 
and linear fluid models in one to three space dimensions. 
We use this approach to adjust  ``quartic'' relaxation parameters in order to enforce fourth
order accuracy for thermal model and diffusive relaxation modes of the Stokes problem. We
apply the resulting scheme for numerical computation of associated eigenmodes,  compare
our results with analytical references and observe  fourth-order accuracy 
when using   ``quartic'' parameters. 
$ $ \\[2mm]
\end{abstract}

\newpage
\section{Introduction}

\smallskip \noindent $\bullet$ \quad 
The lattice Boltzmann scheme is a numerical method for simulation of a wide family of
partial differential equations associated with conservation laws of physics.
The principle is to mimic at a discrete level the dynamics of the Boltzmann equation. 
In this  paradigm, the number  $ \, f(x,\,t)  \,{\rm d}x \, \,{\rm d}v \,  $
of particles at position $x$, time $t$ and velocity $v$ with an uncertainty of 
$    \,{\rm d}x \, \,{\rm d}v \,  $ follows the Boltzmann partial differential equation 
in the phase space (see {\it e.g.}  Chapman and Cowling \cite{CC39}):   
\begin{equation}
\label{1.1}
{{\partial f}\over{\partial t}} \,+\, v \smb \nabla_x f \,=\, Q(f)   \,. \, 
\end{equation}

\bigskip \noindent $\bullet$ \quad 
Note that the left hand side is a simple advection equation whose solution is trivial
through the method of characteristics:
\begin{equation}
\label{1.2}
f (x,\, v ,\, t ) \,=\, f (x - v t ,\, v ,\, 0 )  \qquad {\rm if } \quad 
Q(f) \equiv 0 \, . 
\end{equation}
Remark also that the right hand side is a collision operator, local in space and integral
relative to velocities: 
\begin{equation} \label{1.3}
Q(f) (x,\, v ,\, t ) \,=\, \int  
{\cal C} \big( f(x,\, w ,\, t ), x ,\, v ,\, t \big) \,  {\rm d}w \, , 
\end{equation}
where $ \, {\cal C} (\smb) \,$ describes collisions at a microscopic level. 
Due to microscopic conservation of mass, momentum and energy, an equilibrium distribution
$ \, f^{\rm eq} (x,\, v ,\, t ) \,$ satisfy the nullity of first moments of the
distribution of collisions:
\begin{equation*} \label{1.4}  
\int  \, Q( f^{\rm eq})  (x,\, v ,\, t ) \,    
\begin{pmatrix}   1 \cr v \cr {1\over2} 
 \mid \!\!v \!\! \mid ^2  \end{pmatrix}   \, {\rm d}v \, \,=\,\, 0  \, .  
\end{equation*} 
Such an equilibrium distribution $ \, f^{\rm eq} \,$ satisfies 
classically the Maxwell-Boltzmann distribution. 

\bigskip \noindent $\bullet$ \quad 
The lattice Boltzmann method follows all these physical recommandations
with specific additional options. First, space $x$
is supposed to live in a lattice    $ \,  \cal{ L}  \, $ included in Euclidian  space 
 of dimension $d$. Second, velocity belongs to a finite set    $ \,  \cal{ V}  \, $
composed by given velocities $ \, v_j \,$ $(0 \leq j \leq J)$ chosen in such a way that 
\begin{equation*} \label{1.7}   
   x \in {\cal{ L}}  \,\,\,    {\mathrm  {and} } 
 \,\,     v_j \in  {\cal V}  \, 
 \, \, \Longrightarrow \,  \,    x \,+\,  \Delta t \, v_j \, \in {\cal L}    \,, \,  
\end{equation*} 
where $\, \Delta t \,$ is the time step of the numerical method. Then the distribution of
particles, $ \, f $,  is denoted by $ \, f_j(x,\,t)\,$ with $ \, 0 \leq j \leq J ,  $ 
$\, x \,$ in the lattice   $ \,  \cal{ L}  \, $ and $ \, t \,$ an integer multiple of time
step  $\, \Delta t . \,$

\bigskip \noindent $\bullet$ \quad 
In the pioneering work of cellular automata introduced by 
Hardy, Pomeau and   De Pazzis \cite{HPP73}, 
Frisch, Hasslacher and   Pomeau  \cite{FHP86} and developed by 
 d'Humi\`eres,  Lallemand and  Frisch  \cite{DLF86}, 
the distribution $ \, \, f_j(x,\,t) \,$ was chosen as Boolean. Since
the so-called lattice Boltzmann equation  of  
Mac Namara and Zanetti  \cite{MZ88},   
Higuera, Succi and Benzi \cite{HSB89},  Chen, Chen and  Matthaeus \cite{CCM92}, 
Higuera and Jimenez  \cite{HJ89} 
(see also Chen and  Doolen \cite{CD98}),   
the distribution $ \,  f_j(\smb ,\, \smb) \,$ 
takes real values in a continuum and the 
collision process follows a linearized approach of Bhatnagar,  Gross and  Krook  \cite{BGK}.
With  Qian, d'Humi\`eres and  Lallemand   \cite{QDL92}, the 
 equilibrium distribution  $ \, f^{\rm eq} \,$ is determined with 
a polynomial in velocity. 
In   the work of  Karlin {\it et al}  \cite{KGSB98},  
the equilibrium state is obtained with a general methodology of entropy  minimization.

\bigskip \noindent $\bullet$ \quad 
The numerical scheme is defined through  the evolution of a population  $ \, f_j(x,\,t) $, 
with $ \, x \in  \cal{ L}  \, $ and  $ \, 0 \leq j \leq J \, $ towards a distribution 
 $ \, f_j(x,\,t+\Delta t)\,$ at a new discrete time. The scheme is composed by two steps that
 take into account successively the left and right hand sides of the Boltzmann equation
 (\ref{1.1}). The first step describes the relaxation  $  \, f \longrightarrow f^* \,  $ 
of particle distribution $ \, f \,$ towards the equilibrium. It is local in space and
nonlinear in general. 
D.~d'Humi\`eres    first  introduced in  \cite{DDH92} 
the fundamental notion of moments in the context of lattice Boltzmann schemes.
He defines  an invertible matrix $ \, M \, $ with 
$\, (J+1)\,$ lines and $\, (J+1)\,$ columns and the  moments $ \, m \,$
through a simple linear relation 
\begin{equation} \label{1.8}   
m_k \,=\,   \sum_{j=0}^{J}  M_{k   j} \,  f_j  \,, \qquad 0 \leq k \leq J \,.
\end{equation} 
%

\bigskip \noindent $\bullet$ \quad 
The first $N$ moments are supposed to be at equilibrium: 
\begin{equation} \label{1.9}   
m^*_i = m_i\equiv  m^{\rm eq}_i \equiv W_i  \,, \qquad 0 \leq i \leq N-1 \, 
\end{equation} 
and we introduce the vector $\, W \in \R^N \, $ of conserved variables
composed of  the $ \, W_i \, $  for $ \,0 \leq i \leq N-1 $:
$  W_i  \equiv  m^{\rm eq}_i, \,$ $ 0 \leq i \leq N-1$. 
 The first moments at equilibrium are respectively the total density 
\begin{equation} \label{1.10}  
\rho \, \equiv \, \sum_{j=0}^J f_j \, \,, 
\end{equation} 
momentum
\begin{equation} \label{1.11}  
q_\alpha \, \equiv \, \sum_{j=0}^J v_j^{\alpha} \, f_j \,, 
\qquad 1 \leq \alpha \leq d \,   
\end{equation} 
and possibly the energy  \cite{LL03} 
 for Navier--Stokes fluid simulations. In consequence, we have  
\begin{eqnarray}   \label{1.12}
M_{0 j} \, &\equiv& \, 1 \,, \,\,\, \qquad 0 \leq j \leq J   \\
\label{1.13}
M_{\alpha  j} \, & \equiv & \, v_j^{\alpha} \,, \qquad 1 \leq \alpha  \leq d 
\,, \,\,\,  0 \leq j \leq J \, .  
\end{eqnarray}
For the other  moments, we suppose  given 
$ \, (J+1-N) \,$ (nonlinear) functions $ \, G_k(\smb)\,$ 
\begin{equation} \label{1.14}  
\R^N \ni W \, \longmapsto \, G_k(W) \,\in \R \,, \qquad N \leq k \leq J \,
\end{equation} 
that define equilibrium moments $ \, m^{\rm eq}_k \,$ according to the relation 
\begin{equation} \label{1.15}  
 m^{\rm eq}_k \,=\,  G_k(W)  \,, \qquad N \leq k \leq J \,. 
\end{equation} 
Note also that more complicated models have been developed 
in Yeomans's group (see {\it e.g.} Marenduzzo {\it at al}  \cite {MOCY07}) 
for modelling of  liquid crystals. 
 
\bigskip \noindent $\bullet$ \quad 
The relaxation process is related to the linearized collision operator introduced at
relation  (\ref{1.3}). In particular  intermolecular interactions
(Maxwell molecules with a  $\, 1/r^4$ potential), 
the collision operator is exactly solvable  in terms of so-called Sonine polynomials 
(see {\it e.g.} Chapman and Cowling \cite{CC39}) and the  eigenvectors are known. 
Moreover, the discrete model is highly constrained by symmetry and exchanges of
coordinates.  
%
In the work of   d'Humi\`eres \cite{DDH92},  relaxation parameters 
(also named as $s$-parameters in the following)
$ \, s_k \,$  $\, ( N \leq k \leq J) \,$ 
are introduced, satisfying for stability constraints 
(see {\it e.g.}  \cite{LL00})  the conditions 
\begin{equation*} \label{1.16}  
0 \,<\, s_k  \,<\, 2  \,, \qquad N \leq k \leq J \, . \,
\end{equation*} 
Then the nonconserved momenta  $ \, m^*_k \,$ after collision are supposed to satisfy 
\begin{equation} \label{1.17}  
m^*_k =  m_k + s_k \,(   m_k^{\rm eq} - m_k)  \,, \qquad k \geq N  \, 
\end{equation} 
and we will denote by $ \, S \,$ the diagonal matrix of order $ \, J +1 - N \,$ 
whose diagonal coefficients are equal to $ \, s_k $:
\begin{equation} \label{1.18}  
S_{k\ell} \,\equiv\, \delta_{k\ell} \, s_{\ell} \,,  \qquad k,\, \ell  \geq N \,  
\end{equation} 
with $ \, \delta_{k\ell} \,$ the Kroneker symbol equal to 1 if $ k=\ell$ and null in the
other cases.  
Remark that this framework is general:
when the matrix $S$ is proportional  to identity, the d'Humi\`eres 
 scheme degenerates to the popular
``BGK'' method characterized by a ``Single Relaxation Rate''.
In this particular case the relaxation operator is diagonal and there is no particular
diagonalization basis to work with.
The distribution $ \, f^* \,$ after collision is reconstructed by 
inversion of relation  (\ref{1.8}):
\begin{equation} \label{1.19}  
f_j^* \,=\,   \sum_{\ell=0}^{J}  M^{-1}_{j \ell} \, m_{\ell}^*   \,, 
\qquad 0 \leq j \leq J \,.
\end{equation} 

\bigskip \noindent $\bullet$ \quad 
We suppose also  that the set of velocities $ \, {\cal V } \,$ is invariant by space reflection: 
\begin{equation*} \label{1.20}  
v_j \in  {\cal V } \, \, \, \Longrightarrow \,  \,  \exists \,  \ell \,  \in \{0,\, \dots ,\, J\},
\,\, \, v_\ell \,=\, - v_j \,,\,\,  v_\ell   \in  {\cal V } \,. \, 
\end{equation*} 
The second step is the advection that mimic at the discrete level the free evolution
through characteristics (\ref{1.2}):
\begin{equation} \label{1.21}  
f_j(x,\, t+\Delta t) \,=\, f_j^*(x - v_j \, \Delta t , \, t) \,, \qquad  x \in {\cal L}
,\,\,  0 \leq j\leq J \,, \,\,  v_j \in {\cal V} \, . 
\end{equation} 
Note that all physical relaxation processes are described in space of
moments. Nevertheless, evolution equation (\ref{1.21}) is the 
key issue of forthcoming expansions.

\bigskip \noindent $\bullet$ \quad 
The asymptotic analysis of cellular automata (see {\it e.g.}  H\'enon \cite{He87}) 
provides  evidence supporting   asymptotic   partial differential equations
  and viscosity coefficients related to the induced
parameter defined by
\begin{equation} \label{1.22}  
\sigma_k \, \equiv \, {{1}\over{s_k}} \,-\, {1\over2} \, . 
\end{equation} 
The lattice Boltzmann  scheme (\ref{1.8}) to  (\ref{1.21}) has been analyzed by d'Humi\`eres
\cite{DDH92}  with a Chapman-Enskog method coming from statistical physics. 
Remark that the extension of the discrete Chapman-Enskog expansion to higher order 
already exists (Qian-Zhou \cite{QZ00}, d'Humi\`eres  \cite{DDH07}). But the calculation  in the
nonthermal case ($N > 1$) is quite delicate from an algebraic point of view and introduces
noncommutative formal operators. 
Recently,   Junk and  Rheinl\"ander \cite{JR08} developed a Hilbert type expansion
 for the analysis  of lattice Boltzmann schemes at high order of accuracy. 
We have proposed in previous works \cite{Du07, Du08} 
the Taylor expansion method which is  an extension to the 
lattice Boltzmann   scheme of
the  so-called equivalent partial differential equation method  proposed  independently
by   Lerat and  Peyret  \cite{LP74}  and  by  Warming and  Hyett  \cite{WH74}. 
In this framework, the parameter $ \, \Delta t \,$ is considered as the only infinitesimal
variable and we introduce a {\bf constant}  velocity ratio $ \, \lambda \,$ between 
space step  and time step: 
\begin{equation} \label{1.23}   
\lambda \,\equiv \, {{\Delta x }\over{\Delta t}}  \, . 
\end{equation} 
The lattice Boltzmann scheme is classically considered as second-order accurate 
(see {\it e.g.}  \cite{LL00}).  
In fact, the viscosity coefficients $ \, \mu \,$  relative to second-order terms are
recovered according to a relation of the type
%
\begin{equation*} \label{1.24} 
\mu \,=\,  \zeta \,   \lambda^2 \,  \Delta t \, \sigma_k \,  
\end{equation*} 
for a particular value of label $k$.
The coefficient $ \, \zeta \, $  is equal to $ \, {{1}\over{3}}\, $  for the simplest  
  models that are considered   hereafter. 

\bigskip \noindent $\bullet$ \quad 
A natural question is to extend this accuracy to third or higher orders. 
In the case of single   relaxation times (the BGK variant of  d'Humi\`eres  scheme),  
progresses in this direction have been proposed by 
Shan  {\it et al} \cite{SH98, SYC06}  
and Philippi   {\it et al}  \cite{PHSS06}
 using  Hermite polynomial methodology  for the approximation of the Boltzmann equation.   
The price to pay is an extension of the stencil of the numerical scheme and the practical
associated problems for the numerical treatment of boundary conditions. 
Note also the work of the Italian team (Sbragaglia {\it et al} \cite{SBBSST07},    
Falcucci {\it et al} \cite{FBCCSS07}) on application to  multiphase flows.  
In the context of scheme with multiple relaxation times,  
Ginzburg, Verhaeghe and d'Humi\`eres have analyzed with the
Chapman-Enskog method the ``Two Relaxation Times'' version of the  scheme  
 \cite{GVD08a, GVD08b}.  
A nonlinear extension of this  scheme, the so-called ``cascaded lattice Boltzmann method'' 
has been proposed by
Geier {\it et al} \cite{GGK06}. It gives  also high order accuracy and the analysis is 
under development (see {\it e.g.} Asinari \cite{As08}). 
%
%
  The general nonlinear extension of the  Taylor expansion method  
to third-order of accuracy  of  d'Humi\`eres scheme 
is presented in \cite{Du09}. It  provides 
   evidence of the importance of the so-called tensor of momentum-velocity defined by 
\begin{equation} \label{1.25} 
\Lambda_{k p}^{\ell} \equiv  
\sum_{j=0}^{J}  \, M_{k   j} \, M_{p   j} \, M^{-1}_{j \ell} 
\,, \qquad  0 \leq k , \, p ,\, \ell \leq J \,. 
\end{equation} 
Moreover, it shows also that for  athermal Navier Stokes equations, 
the mass conservation equation 
contains a remaining term of third-order accuracy that cannot be set 
 to zero by fitting relaxation parameters \cite{Du09}.  

\bigskip \noindent $\bullet$ \quad 
Our motivation in this contribution 
is to show that it is possible to extend 
the order of accuracy of an {\it existing a priori} second-order accurate 
 lattice Boltzmann scheme to higher orders. 
We use the Taylor expansion method  \cite{Du08}  to determine the 
equivalent partial differential equation of 
 the  numerical scheme   to higher orders of accuracy. 
Nevertheless, it is quite  impossible to determine 
explicity the entire expansion 
in all generality in the nonlinear case. In consequence, we restrict here to a first step. 
We propose in the following a general methodology for deriving the equivalent equation of
the  d'Humi\`eres  scheme at an arbitrary order when the collision process defined by the functions 
$\, G_k \,$ of relation  (\ref{1.14}) are {\bf linear}.  
This calculation leads to explicit developments that can be expanded   with the help of
formal calculation. This work is detailed in Section~2. 
In Section~3 we apply the general methodology to classical linear models of thermics and
linearized athermal Navier Stokes equations. We treat fundamental examples from one to
three space dimensions. When it is possible, the equivalent partial equivalent equations
are explicited. 
In Section~4, we use the fourth-order equivalent equation of two and three-dimensional
 models to enforce accuracy through  a proper  choice of ``quartic'' parameters. For a scalar
heat equation, the effect of the precision of the numerical computation of eigenmodes is
presented. For linearized  athermal Navier Stokes equations,  
we propose a method for enforcing the precision of the eigenmodes of the associated partial
differential equation. First numerical results show that for appropriate tuning values of
the parameters,  fourth-order precision is achieved.



\section{Formal development of linearized d'Humi\`eres  scheme}

\smallskip \noindent $\bullet$ \quad 
In what follows, we suppose that the collision process is {\bf linear} {\it i.e.} that the
$ \, G_k\,$ functions introduced in (\ref{1.14})  (\ref{1.15}) are  linearized
around some reference  state. With this hypothesis, we can write: 
\begin{equation*} \label{1.26} 
G_k(W) \,\equiv \, \sum_{j=0}^N  G_{k j } \, W_{j}  \,\,=\,\, 
 \sum_{j=0}^N  G_{k j} \, m_{j}  \,, \qquad k \geq N \, . 
\end{equation*} 
To be precise, putting together relations (\ref{1.15}) and  (\ref{1.17}), there
exists a $ \, (J+1) \times (J+1) \,$ matrix $ \, \Psi \,$ such that the 
collisioned momentum $ \, m^* \,$ defined in (\ref{1.17}) is a linear combination of the 
moments before collision: 
\begin{equation} \label{1.27} 
m^* \,=\, \Psi \, \smb \, m \,, \qquad m^*_k \,=\, \sum_{j=0}^J \Psi_{k\ell} 
\, \,  m_{\ell}  \,. 
\end{equation} 
Of course, the conservation (\ref{1.13}) implies that $ \, \Psi \,$ has a structure of the type 
\begin{equation} \label{1.28} 
\Psi \,=\, \begin{pmatrix}    {\rm I} & 0  \cr  \Phi  &  {\rm I} -  S   \end{pmatrix}  \, . 
\end{equation} 
The top left  block of the right hand side of  (\ref{1.28})
is the identity matrix of dimension $N$ and the bottom  left block is
described through  the $ \, G_k \,$ functions introduced in  (\ref{1.14})   (\ref{1.15}): 
\begin{equation} \label{1.29}  
 \Phi_{k j}  \, \equiv \,  \Psi_{k j}  \, = \, s_k \,\, G_{k j }   
\,, \qquad j < N \,, \quad k \geq N \, . 
\end{equation} 
The  bottom right  block  of the right hand side of  (\ref{1.28}) 
contains the coefficients 
$ \, 1 - s_k \, $ ($k \geq N$) related to relaxation (\ref{1.18}). 
 
\smallskip \noindent $\bullet$ \quad 
In order to make  our result explicit, we need some notations. 
We introduce multi-indices  $\gamma$, $\delta$, $\varepsilon$ in $\, \{1,\dots \, d\}^q
\,$ in order to represent  multiple differentiation with respect to space. If 
\begin{equation*} \label{2.3}   
\gamma \,=\, \big( \, \underbrace{1 ,\, \dots \, 1}_{\alpha_1 \, {\rm times}}  \,, 
\dots \,, \,  \underbrace{d ,\, \dots \, d}_{\alpha_d \, {\rm times}}  \, \big) \,, 
\end{equation*} 
then 
\begin{equation*} \label{2.4}   
\partial_{\gamma}  \,\equiv \, {{\partial^{\alpha_1}}\over{\partial x_1^{\alpha_1}}} 
\, \cdots \, 
 {{\partial^{\alpha_d}}\over{\partial x_d^{\alpha_d}}} \,  
\end{equation*} 
and we denote by $ \mid \! \gamma \!\mid \,$ the length of multi-index $ \, \gamma$: 
\begin{equation*} \label{2.5}   
 \mid \! \gamma \!\mid  \, \equiv \, \alpha_1 \, + \, \cdots \, + \, \alpha_d \, . 
\end{equation*} 
Then thanks to the binomial formula for iterated  differentiation, we can introduce coefficients 
$\, P_{\ell \gamma}\,$ in order to satisfy the identity 
\begin{equation} \label{2.6}   
 \Big(- \sum_{\alpha=1}^d M_{\alpha \ell} \, \partial_{\alpha}
 \Big)^q \,\equiv \, \sum_{ \mid   \gamma  \mid = q } P_{\ell \gamma} \, 
\partial_{\gamma}\, .  
\end{equation} 
for any integer $ \, q . \, $

\smallskip \noindent $\bullet$ \quad 
We first establish that at first-order of accuracy, we have a representation of nonconserved
moments in terms of conservative variables: 
\begin{equation} \label{2.11}   
m_k  \,=\, \sum_{j=0}^J B^0_{k j} \,  W_j  
 \,   + {\rm O}(\Delta t)   \,, \qquad k \geq N  \,.  
\end{equation}  
with 
\begin{equation} \label{2.10}   
B^0_{k j} \,\equiv \,  {{1}\over{s_k}} \,   \Psi_{k j}  \,, \qquad k \geq N \,, \,\,\,
0 \leq j \leq N-1   \, . 
\end{equation}  
We have also the first-order conservation law 
\begin{equation} \label{2.14}   
{{\partial W_i}\over{\partial t}} \, + \, \sum_{\mid \gamma \mid = 1}
A^{\gamma}_{ij} \,  \partial_{\gamma} W_j \,=\,  \, {\rm O}(\Delta t)   \,, 
\quad 0 \leq i \leq N-1  \, . \,   
\end{equation}    
with coefficients $ \, A^{\gamma}_{ij} \, $ given according to 
\begin{equation} \label{2.13}   
A^{\gamma}_{ij} \, \equiv \,  \sum_{p=0}^J \, 
 \Lambda_{\gamma i}^p \,  \Big(  \Psi_{p j} \,  + \, 
 \sum_{\ell \geq N}  \Psi_{p \ell} \,    {{1}\over{s_{\ell}}} \,   \Psi_{\ell j} 
 \Big) \,, \quad \mid \! \gamma \!\mid = 1 \,, \,\,\, 0 \leq i ,\, j \leq N-1  \, . 
\end{equation}  
The proof of this result and those that follow   of this Section are detailed in
Appendix~A.

\smallskip \noindent $\bullet$ \quad 
The expansion of moments (\ref{2.11}) can be extended to second-order accuracy:
\begin{equation} \label{2.16}    
m_k  \,=\,   \sum_{0 \leq \mid \gamma \mid \leq 1}  \Delta t ^{\mid \gamma \mid} \, \, 
 B^{\gamma}_{kj} \, \,  \partial_{\gamma} W_j   \,   + \, {\rm O}(\Delta t^2)  \,.  
\end{equation}    
with
\begin{equation} \label{2.15}  \left\{ \begin{array}{c} \displaystyle  
B^{\gamma}_{kj} \, = \,   {{1}\over{s_{k}^2}} \,  \sum_{i=0}^{N-1} 
   \Psi_{k i} \, A^{\gamma}_{ij} \, - \, 
 {{1}\over{ s_{k} }} \,  \sum_{r=0}^{J} \, {{1}\over{  s_{r}}}   \sum_{p=0}^{J} 
 \Lambda_{\gamma k}^p \,  \Psi_{p r} \,  \Psi_{r j} \, \,, 
 \qquad  \qquad \qquad    \qquad    \\   [1mm] 
 \displaystyle \qquad \qquad  \qquad \qquad  \qquad \qquad 
\qquad  \mid \! \gamma \!\mid = 1 \,, \,\,\, k \geq N \,, \,\,\, 0 \leq  j \leq N-1 \, . 
\end{array}  \right.  \end{equation}  
\smallskip \noindent $\bullet$ \quad 
Then we extend the previous expansions  (\ref{2.14}) and  (\ref{2.16}) to  any order
$\sigma.$
By induction, we establish that we have an equivalent partial differential equation of the
form 
\begin{equation} \label{2.17}    
{{\partial W_i}\over{\partial t}} \, + \, \sum_{1 \leq \mid \gamma \mid \leq \sigma} 
\,  \Delta t ^{\mid \gamma \mid - 1} \,  \sum_{j=0}^{N-1} \,
A^{\gamma}_{ij} \,  \, \partial_{\gamma} W_j \,=\,   {\rm O}(\Delta t^{\sigma}) \,, 
\quad 0 \leq i \leq N-1  \,,\,   
\end{equation}    
and an expansion of nonconserved moments as
\begin{equation} \label{2.18}    
m_k  \,=\,   \sum_{0 \leq \mid \gamma \mid \leq  \, \sigma}  
\Delta t ^{\mid \gamma \mid} \,  \sum_{j=0}^{N-1}\, 
 B^{\gamma}_{kj} \, \,  \partial_{\gamma} W_j   \, + \, {\rm O}(\Delta t^{\sigma + 1})
 \,, \quad k \geq N \, ,    
\end{equation}    
with the following recurrence relations for defining the coefficients $ \, A^{\gamma}_{ij} \,$
and $\,  B^{\gamma}_{kj}  $: 
%
\begin{equation} \label{2.56} 
C^{1, \gamma}_{ij} \, = \,  A^{\gamma}_{i j} \,,   \quad 
0 \leq i,\, j \leq N-1 \, ,   
\end{equation} 
\begin{equation} \label{2.21}  \left\{ \begin{array}{c} \displaystyle 
C^{q+1, \gamma}_{ij} \, = \, -  \!\!\! 
\sum_{\delta \geq q, \,   \varepsilon  \geq 1,  \, 
\delta +  \varepsilon =  \gamma} \, 
 \sum_{\ell = 0 }^J \, 
C^{q, \delta}_{i \ell} \, A^{\varepsilon}_{\ell j} \,,  
 \qquad  \qquad \qquad    \qquad    \\   [1mm] 
 \displaystyle \qquad \qquad  \qquad \qquad  \qquad \qquad 
 \quad    2 \leq q+1  \leq  \, \mid \! \gamma \! \mid   \, , 
\,\,\,  0 \leq i,\, j \leq N-1 \, ,    
\end{array}  \right.  \end{equation} 
\begin{equation}  \label{2.24}  \left\{ \begin{array}{c} \displaystyle
A^{\gamma}_{ij} \,= \displaystyle 
\,-   \sum_{q=2}^{\mid \gamma \mid } 
\,  {{1}\over{ q ! }} \, C^{q, \gamma}_{ij} \, 
  \qquad \qquad   \qquad \qquad  \qquad \qquad \qquad  \qquad \qquad    \qquad    \\   [1mm] 
 \displaystyle \qquad \, -  \!\!\!   
\sum_{1 \leq \mid \delta \mid \leq \mid \gamma \mid, \,  
0 \leq \mid \varepsilon \mid \leq  \mid \gamma \mid - 1,  \,  \delta +  \varepsilon = \gamma} \, 
\,  \sum_{\ell = 0 }^J  \sum_{p = 0 }^J  \sum_{r = 0 }^J  \,
    {{ 1 } \over{ \mid   \delta   \mid \, !   }} \,
 M_{i \ell} \, M^{-1}_{\ell p} \, \Psi_{p r} \, P_{\ell \delta} \, 
 B^{ \varepsilon}_{rj} \, \,,     
\end{array}  \right. \end{equation} 
\begin{equation} \label{2.57} 
D^{0,\gamma}_{kj} \,= \, B^{\gamma}_{kj}  \,, \qquad  k \geq N \,, \quad 0 \leq j \leq N-1 \,,  
\end{equation} 
\begin{equation}   \label{2.55}  \left\{ \begin{array}{c} \displaystyle 
D^{q+1,\gamma}_{kj} \,= \,  \displaystyle - 
 \sum_{ \mid \delta \mid  \geq q , \, 
 \mid \varepsilon \mid \geq 1,\, \delta +  \varepsilon =  \gamma} \,  \sum_{\ell = 0 }^J \, 
 D^{q,\delta}_{k \ell} \,  \, A^{\varepsilon}_{\ell j} \, \,,  
 \qquad  \qquad \qquad    \qquad    \\   [1mm]   
  \qquad \qquad \qquad  \qquad 
1 \leq  q+1 \leq    \mid \! \gamma \! \mid  \,, \quad  k \geq N \,, \,\, 
0 \leq j \leq N-1 \, ,   
\end{array}   \right. \end{equation}
\begin{equation} \label{2.25}  \left\{ \begin{array}{c} \displaystyle
B^{\gamma}_{kj} \,= \displaystyle 
\, {{1}\over{s_k}} \,  \bigg( -  \sum_{1 \leq q \leq  \mid \gamma \mid} \, 
 {{1}\over{ q \, ! }} \,  D^{q,\gamma}_{kj}    
\qquad   \qquad \qquad \qquad   \qquad \qquad \qquad  \qquad    \\   [1mm]   \displaystyle
\qquad  \,+  \!\!\!\!\!\!  
\sum_{1 \leq \mid \delta \mid \leq  \mid \gamma \mid, \,
0 \leq \mid \varepsilon \mid \leq   \mid \gamma \mid - 1 ,\,  \delta +  \varepsilon = \gamma}
 \,  \,  \sum_{\ell = 0 }^J  \sum_{p = 0 }^J  \sum_{r = 0 }^J  \, 
 \,    {{1} \over{\mid \! \delta \!  \mid ! }} \,  
 M_{k \ell} \, M^{-1}_{\ell p}  \, \Psi_{p r} \, P_{\ell \delta} \,  B^{\varepsilon}_{rj}
  \bigg) \,,    \qquad    \\   [1mm]    \displaystyle
  \qquad \qquad \qquad  \qquad    \qquad  \qquad \qquad  \qquad  \qquad \qquad  \qquad   
k \geq N  \,, \,\,   0 \leq  j \leq N-1 \, . 
\end{array}  \right. \end{equation} 

\smallskip \noindent $\bullet$ \quad  
Note that the  results  (\ref{2.24}) and  (\ref{2.25}) 
are coupled through the relations  (\ref{2.56})  (\ref{2.21})  (\ref{2.57}) 
and  (\ref{2.55}).  For example, the evaluation of coefficient 
$\, D^{q+1,\gamma}_{kj} \, $  uses explicitly $ \,  A^{\varepsilon}_{\ell j} \, ,\, $
the evaluation of  $ \,  A^{\gamma}_{i  j} \, $ uses  
 $ \,  B^{\varepsilon}_{r  j} \, $ 
and the computation of  $ \,  B^{\gamma}_{r  j} \, $ is impossible 
if $\, D^{q,\gamma}_{kj} \, $ is not known. 
The proof is detailed in Appendix~A.  
It is an elementary and relatively lengthy  algebraic calculation. 
In particular,  our mathematical framework is classical: all differential  operators commute
and the technical difficulties of noncommutative time derivative operators associated with
the use of formal Chapman-Enskog method \cite{DDH07}  vanish. 
As a result, the general expansion of a linearized  d'Humi\`eres    scheme 
at an arbitrary order can be obtained
by making explicit   the coefficients  $ \, A^{\gamma}_{ij} \,$
and $\,  B^{\gamma}_{kj} .\, $  Remark that the hypothesis of
linearity allows making the above formulae   explicit   and 
we have done this work with the help of formal calculation. 
Nevertheless, it is always
possible to suppose that the $\, G_k \,$ functions are linearized expansions of a
nonlinear equilibrium. In this case, the previous  equivalent high order  partial
differential equations  (\ref{2.17}) give  a very good information concerning the behavior 
 of the  scheme.

\section{Equivalent Thermics and Fluid equations }
 
\smallskip \noindent $\bullet$ \quad 
We make explicit in this section the fourth-order equivalent equation 
(\ref{2.17})  of some   lattice Boltzmann schemes for
two fundamental problems of mathematical physics:  thermics and linearized
athermal Navier-Stokes equations. We treat first advective thermics in one space
dimension with the so-called D1Q3 lattice Boltzmann scheme. 
In order to obtain  results presentable on a sheet of paper, we simplify the
model and omit the advective term for two (D2Q5) and three (D3Q7) space dimensions. 
Secondly  we study  linearized  athermal  Navier-Stokes equations in 
one (D1Q3),  two (D2Q9)  and three (D3Q19) space dimensions. 
Note that we have to define precisely our results. 
 First the numbering of degrees of freedom
{\it via}   corresponding graphics  is  specified; see   Appendix~B.   
The  choice of moments, {\it id est} the $M$
matrix, is also made precise    in  Appendix~B.
Secondly    the $ \, \Psi \, $ matrix 
of relation (\ref{1.27}) is specified,  later in this section.

\bigskip  \noindent $\bullet$ \quad 
{\bf   D1Q3 for advective thermics at fourth-order }

\smallskip \noindent     
For a thermics problem, we have only one conserved quantity. Then $\, N=1 \,$ 
in relation  (\ref{1.9}). The two   nonconserved moments (momentum 
$ \, q^{\rm eq} \, $   and energy $ \, \epsilon^{\rm eq}$; see (\ref{7.1}))  
at equilibrium are supposed to be {\bf linear} functions 
of the conserved momentum $\, \rho $: 
\begin{equation} \label{3.2}  
q^{\rm eq}  \,=\, u \, \lambda \, \rho \,, \quad \epsilon^{\rm eq}  \,=\, 
\alpha  \, {{\lambda^2}\over{2}} \, \rho    \, . 
\end{equation} 
Due to  (\ref{1.29}) and (\ref{3.2}),  the matrix     $ \, \Psi \,$  for
dynamics relation  (\ref{1.27}) is given according to 
\begin{equation*} \label{3.3}  
 \Psi  \,=\,  \begin{pmatrix}   1 & 0 & 0 \cr s_1 \, u \, \lambda  & 1-s_1 & 0 \cr 
\alpha \, s_2 \, \lambda^2 / 2 & 0 & 1 - s_2  \end{pmatrix}   \, . 
\end{equation*} 
We  determine without difficulty the equivalent partial differential equation for
this lattice Boltzmann scheme at order four, to fix the ideas. 
For $ \, i = 1 ,\, 2 ,\,$ we introduce $ \, \sigma_i \,$ from 
relaxation time  $ \, s_i \,$ 
according to relation (\ref{1.22}). 
When a drift in velocity $\, u  \,$ is present, note that 
the diffusion coefficient is a function of mean value velocity. We have 
\begin{equation} \label{3.4}   
{{\partial \rho}\over{\partial t}} \,+\, 
  u \, \lambda \,   {{\partial \rho}\over{\partial x}}  
\,-\, \sigma_1 \, \Delta t \, \lambda^2 \, ( \alpha  -  u^2) \,  
  {{\partial^2 \rho}\over{\partial x^2}} \,
 + \,  \kappa_{3}   \, {{\Delta t^2 \, \lambda^3}\over{12}} \, 
    {{\partial^3 \rho}\over{\partial x^3}} \,+\,  \kappa_{4} \,
 {{\Delta t^3 \, \lambda^4}\over{12}}
\,   {{\partial^4 \rho}\over{\partial x^4}}  \,=\,   {\rm O}(\Delta t^4)
\end{equation}  
with parameters $\, \kappa_{3} \,$ and  $\, \kappa_{4} \,$ given according to 
%
\begin{equation*}   \begin{array} {rcl}    \label{3.5} 
\kappa_{3}   & \,  = \, &  \,  -u \,  \Big(  2\, \big( 1 - 12 \,  \sigma_1^2 \big) \, u^2  \,+\,  
  1 - 3 \, \alpha   \,-\, 12 \, \sigma_1 \, \sigma_2 \, (1 - \alpha) 
\,+\,  24 \,  \sigma_1^2 \, \alpha    \,  \Big)  \,  \\  
\kappa_{4} &  \,   =  \, &   \big( -9 \,+\, 60\,\sigma_1^2 \big) \, \sigma_1 \,u^4 
\,+\,  \big( - 5  \, (1 - 3 \, \alpha)  \,  \sigma_1
  \,-\, 3  \, (1 -  \alpha)  \,  \sigma_2 \,+\, 
  \\   &    &    \displaystyle     
\, +\,12 \,   (1 -  \alpha)  \,  \sigma_1 \,  \sigma_2^2  
\,+\, 36 \,   (1 -  \alpha)  \,  \sigma_1^2   \,  \sigma_2  
\,-\, 72 \,  \sigma_1^3 \, \alpha \big) \,   \,u^2 \, 
  \\    &  &    \displaystyle  
 +\,  \alpha \, \sigma_1 \, \big ( 2 - 3 \, \alpha  
\,-\, 12 \,    (1 -  \alpha)  \,  \sigma_1 \,  \sigma_2 
\,+\,  12 \,  \alpha   \,  \sigma_1^2 \big) \,     .    
 \end{array} \end{equation*}  
If $ \, u=0 ,\,$ then $ \, \kappa_{3} = 0 \,$ and the   scheme is
equivalent to an advection-diffusion equation up to  third-order accuracy. 
In this particular case, the scheme is fourth-order accurate 
in the previous sense if we set 
\begin{equation*} 
 \sigma_2  \,=\, \frac{ 2 - 3 \alpha 
+ 12  \,\alpha \,  \sigma_1^2}{12 \, \sigma_1  \, (1 - \alpha)}  \, .  
\end{equation*}

 \bigskip    \noindent $\bullet$ \quad 
{\bf   D2Q5 for pure  thermics at fourth-order }
 
\smallskip \noindent
%
%
We have $\, J=4$ and $N=1$.  The equilibrium 
energy (momentum $m_3$  in (\ref{7.14}) with the labelling conventions of Section~1)
is the only one to be non equal to zero. The  matrix$ \,\Psi \,$   of relation  (\ref{1.27})
is now given by the relation
\begin{equation}  \label{3.15}  
 \Psi   \,=\,  \begin{pmatrix}   1 & 0  & 0  & 0 & 0  \cr 
                   0   & 1 - s_1  & 0  & 0 & 0   \cr 
                 0  &0  &  1 - s_1  & 0 & 0   \cr 
                  \alpha \, s_3   &0  &  0  &  1 - s_3  & 0 \cr 
                  0   &0  &  0  & 0 &    1 - s_4   \end{pmatrix} \, . 
\end{equation}  
We have developed the conservation law up to fourth-order:  
\begin{equation}   \left\{ \label{3.16}  \begin{array} {c} \displaystyle 
  {{\partial \rho}\over{\partial t}}    \,-\, 
  {{ \lambda^2 \, \Delta t }\over{10}} \, \sigma_1 \, (4 + \alpha) \,  
 \bigg( {{\partial^2 \rho}\over{\partial x^2 }} \,+\,
  {{\partial^2 \rho}\over{\partial y^2 }} \bigg)  
  \qquad   \qquad   \qquad   \qquad    \qquad   \qquad       \\   [2.7mm] 
   \displaystyle    \qquad   \qquad     
 \,+\,   {{\Delta t^3 \,  \lambda^4}\over{1200}} \, \sigma_1 \,  
 (4 + \alpha) \, \bigg(  \kappa_{40} \,
 \Big( {{\partial^4 \rho}\over{\partial x^4 }} +
{{\partial^4 \rho}\over{\partial y^4 }}  \Big)     
 \,+\, \kappa_{22} \,   {{\partial^4 \rho}\over{\partial x^2 \partial y^2}}  \bigg) 
  \,=\,   {\rm O}(\Delta t^4) 
\end{array} \right.  \end{equation} 
and the $\, \kappa \, $ coefficients are explicited as follows: 
\begin{eqnarray}    \label{3.17}    
\kappa_{40}   \,=\, &  8 \,-\,3 \, \alpha \, 
+  \, 12 \,  (\alpha + 4) \, \sigma_1^2 
\, - \, 12 \, (1 - \alpha)\, \sigma_1 \,  \sigma_3 \,-\,  60 \,  \sigma_1 \, \sigma_4  
  \\   \label{3.18} 
\kappa_{22}   \,= \,&  -6 \,   (\alpha + 4) \, +\,  24\,   (\alpha + 4)  \, \sigma_1^2 
 \, - \, 24 \, (1 - \alpha)\, \sigma_1 \,  \sigma_3 \,+\, 120 \,\sigma_1 \, \sigma_4 \,
 . \,   
\end{eqnarray}

\bigskip \noindent $\bullet$ \quad 
{\bf  D2Q9 for  advective thermics at fourth-order }

\smallskip \noindent  
The lattice Boltzmann model  D2Q9 
for a passive scalar (see   
Chen, Ohashi and  Akiyama  \cite {COA94},  
Shan \cite {Sh97}, Ginzburg \cite{Gu05}) 
is obtained from the 
  D2Q5 model by adding four velocities 
along the diagonals (Figure \ref{fd2q5q9}, right). 
The evaluation of matrix $M$ is absolutely nontrivial and is precised at 
(\ref{7.19}).
The dynamics is given by  
\begin{equation}  \label{3.20} 
   \Psi =    \begin{pmatrix}  1 &  \!\! \!\!0  & \!\!0  & \!\!0 & \!\!0  & \!\!0  
& \!\!0  & \!\!0  & \!\!0 \cr 
u \, \lambda \, s_1  &  \!\! \!1 \!\!-\!\! s_1  & \!\!0  &\!\!0  &  \!\!0  & \!\!0 
&  \!\!0    &  \!\!0  &  \!\!0  \cr 
v \, \lambda \, s_1  &  \!\! \!\!0  &  \!1 \!\!-\!\! s_1  & \!\!0 & \!\!0  & \!\!0  
& \!\!0  & \!\!0  & \!\!0 \cr 
a_3 \, s_3   &  \!\!\!\!0  &  \!\!0  &  \!1 \!\!-\!\! s_3  & \!\!0  & \!\!0  & \!\!0  
& \!\!0  & \!\!0\cr 
a_4 \, s_4   &  \!\!\!\!0  &  \!\!0  & \!\!0 &    \!1 \!\!-\!\! s_4  & \!\!0  
& \!\!0  & \!\!0  & \!\!0\cr 
a_5   \, u \, s_5   &  \!\!\!\!0  &  \!\!0  & \!\!0 &  \!\!0    &  \!1 \!\!-\!\! s_5  
& \!\!0  & \!\!0  & \!\!0 \cr 
a_6  \, v \, s_5   &  \!\!\!\!0  &  \!\!0  & \!\!0 &  \!\!0    &  \!\!0  
&  \!1 \!\!-\!\! s_5  & \!\!0  & \!\!0 \cr 
a_7 \, s_7   &  \!\!\!\!0  &  \!\!0  & \!\!0 &  \!\!0    &  \!\!0  &  \!\!0  
&  \!1 \!\!-\!\! s_7 & \!\!0 \cr 
a_8 \, s_8   &  \!\!\!\!0  &  \!\!0  & \!\!0 &  \!\!0    &  \!\!0  &  \!\!0  
&  \!\!0  &  \!1 \!\!-\!\! s_8    \end{pmatrix}    \, . 
\end{equation} 
The coefficients $\,a_3 \,$ to $ \, a_8 \,$ in relation  (\ref{3.20})   are chosen in order to obtain
the advection diffusion equation  at  order 2:  
\begin{equation}  \label{3.21}  
 {{\partial \rho}\over{\partial t}} \,+\, 
  \lambda \, \bigg( u \,    {{\partial \rho}\over{\partial x}} \,+\, 
 v \,    {{\partial \rho}\over{\partial y}} \bigg) 
\,-\,  \lambda^2 \, \xi \, \sigma_1 \, \Delta t \,  
 \bigg( {{\partial^2 \rho}\over{\partial x^2 }} \,+\,
  {{\partial^2 \rho}\over{\partial y^2 }} \bigg) \,=\, {\rm O}(\Delta t)^2 \, .  
\end{equation} 
We have precisely: 
\begin{equation*}  \label{3.22}  
 a_3 \,=\,  3 \, (u^2 + v^2) \,-\, 4 \,+\,  6\, \xi \,,\quad 
 a_7 \,=\,   u^2 - v^2   \,,\quad a_8\,=\,   u\, v \,    
\end{equation*} 
as explained in  \cite{Du09}. 
When $ \, u = v = 0 ,\,$ the equation  (\ref{3.21}) takes the form 
\begin{equation*}  \label{3.24}   
 {{\partial \rho}\over{\partial t}}  
\,-\,  \lambda^2 \, \xi \, \sigma_1 \, \Delta t \,  
 \bigg( {{\partial^2 \rho}\over{\partial x^2 }} \,+\,
  {{\partial^2 \rho}\over{\partial y^2 }} \bigg)  \, 
\,+\, {{\lambda^4  \,  \Delta t^3 \, \xi}\over{36}} \,  \bigg( 
\kappa_{40} \,  \bigg( 
  {{\partial^4 \rho}\over{\partial x^4 }} \,+\,   {{\partial^4 \rho}\over{\partial y^4 }}
\bigg)  \,+\, \kappa_{22} \,   {{\partial^4 \rho}\over{\partial x^2 \,\partial y^2  }}    \bigg)  
 \,=\,   {\rm O}(\Delta t^4)
 \end{equation*}  
with coefficients $ \, \kappa_{40} \,$ and  $ \, \kappa_{22} \,$ evaluated according to 
\begin{equation*}   \begin{array} {rcl}    \label{3.25bis} 
\kappa_{40} & \,=\, &  \sigma_1   \,  \Big(  
2 \, \sigma_5 \, (\sigma_7 -\sigma_3) \, (a_4 - 4) \,+\, 6 \, \xi \, 
\big( 1 - \sigma_1 \, \sigma_7 - 5 \,  \sigma_1 \, \sigma_3 
 + 2 \,  \sigma_5 \,  (\sigma_7 -\sigma_3)     \big) \Big)  \\  
\kappa_{22} & \,=\, &  2 \, \big(  \sigma_1 +  \sigma_5  
-2 \,  \sigma_1 \,  \sigma_5 \, ( \sigma_3  + \sigma_7  + 4 \,\sigma_8 )   \big)  \, (a_4 - 4) \\
&& + 12 \, \xi \,  \big(  \sigma_5 + 3 \,  \sigma_1  
-2 \, \sigma_1  \, \sigma_5 \, ( \sigma_3  + \sigma_7 ) 
-2 \, \sigma_1  \, \sigma_3 \,  \sigma_5  
-8 \, \sigma_1  \, \sigma_8 \,  ( \sigma_1  + \sigma_5  )   
+  \sigma_1^2  \, \sigma_7   \big)   \, . \, %
  \end{array}  \end{equation*}    
Remark that the equivalent partial differential equation of this general 
lattice Boltzmann scheme  has been exactly derived in a
complex case where all the time relaxations are {\it a priori} distinct. 
The coefficients $ \,  \kappa_{40} \,$  and  $ \, \kappa_{22} \,$ of the fourth-order
terms are polynomials of degree 3 in the $ \, \sigma$   coefficients. 
When we make the ``BGK hypothesis'' {\it id est} that all the  $\sigma$ coefficients   are
equal,  
%
%
a first possibility for killing  the coefficients $ \, \kappa_{40} \,$  and 
$ \, \kappa_{22} \,$  is given by: 
\begin{equation*} 
\sigma_1 \,=\,  \sigma_1 \,=\,\sigma_3 \,=\,\sigma_4 \,=\,\sigma_5 
\,=\,\sigma_7 \,=\,\sigma_8 \,=\, \frac{1}{6} \,, \qquad \xi = 0  \, . 
\end{equation*}   
We observe that this choice of parameters 
is without any practical interest because the diffusion term in   (\ref{3.21}) is null. 
We observe that a second possibility 
\begin{equation*} 
\xi \,=\, \frac{2}{3} \,  \frac{1 - 6\, \sigma_1^2}{1 - 8\, \sigma_1^2} \,, \qquad 
a_4 \,=\, -2 \,  \frac{1 - 2\, \sigma_1^2}{1 - 8\, \sigma_1^2} 
\end{equation*}   
induces also a fourth-order accurate lattice Boltzmann scheme. 
If we replace the strong  ``BGK hypothesis''  
by the weaker one associated to ``Two Relaxation Times'' as suggested by Ginzburg,
Verhaeghe and d'Humi\`eres in \cite{GVD08a, GVD08b},  {\it id est} 
\begin{equation*}  \label{d2q9-trt} 
\sigma_1 \,=\,\sigma_5  \,, \qquad  \sigma_3 \,=\, \sigma_4 \,=\,\sigma_7 \,=\,\sigma_8  \,,   
\end{equation*}   
we can achieve formal  fourth-order accuracy for 
\begin{equation*} 
  \sigma_1  \,=\, \frac{1}{\sqrt{12}}   \quad {\rm and } \quad 
  \sigma_3  \,=\, \frac{1}{\sqrt{3}} \, .  
\end{equation*}   
%

\bigskip \noindent $\bullet$ \quad 
{\bf   D3Q7 for pure thermics }

\smallskip \noindent 
For three-dimensional thermics, one only needs  a seven point scheme and use the so-called
D3Q7 lattice Boltzmann scheme whose stencil is described in  Figure \ref{fd3q7q19}. 
The matrix 
$M$ is given at relation (\ref {7.36}). 
The dynamics of this Boltzmann scheme uses the following matrix for computation 
of out of  equilibrium moments, according to relation (\ref{1.27}): 
\begin{equation*}  \label{3.37}   
  \Psi  \,=\,   \begin {pmatrix} 1 & 0  & 0  & 0 & 0  & 0  & 0   \cr  
0  &  1-s_1  & 0       & 0     & 0      & 0      & 0   \cr
0  &  0      & 1-s_1   & 0     & 0      & 0      & 0   \cr
0  &  0      & 0       & 1-s_1 & 0      & 0      & 0   \cr
0  &  0      & 0       & 0     & 1-s_4  & 0      & 0   \cr
0  &  0      & 0       & 0     & 0      & 1-s_4  & 0   \cr
\alpha \, s_6   &  0      & 0       & 0     & 0      & 0   & 1-s_6   \end {pmatrix}   \, . 
\end{equation*} 
The equivalent 
thermal scalar conservation law now takes the following form at fourth-order of accuracy: 
\begin{equation*}  \label{3.38}   \begin{array}{rcl} \displaystyle
{{\partial \rho}\over{\partial t}} \,-\, 
  {{ \lambda^2 \, \Delta t }\over{21}} \, \sigma_1 \, (\alpha + 6) \,  \Delta  \rho  \, 
+ \,  {{\Delta t^3 \,  \lambda^4}\over{1764}} \, \sigma_1 \,   (\alpha+6) \, \bigg( 
\kappa_{400} \,  \Big( {{\partial^4 \rho}\over{\partial x^4 }} \,+\,  
   {{\partial^4 \rho}\over{\partial y^4 }}  \,+\,  
   {{\partial^4 \rho}\over{\partial z^4 }} \Big)   &&
\\  \displaystyle
\,+\, \kappa_{220} \,  
 \Big( {{\partial^4 \rho}\over{\partial x^2 \partial y^2}} \,+\, 
 {{\partial^4 \rho}\over{\partial y^2 \partial z^2}} \,+\,    
 {{\partial^4 \rho}\over{\partial z^2 \partial x^2}}   \Big) \bigg)  
  & =  &   {\rm O}(\Delta t^4) 
\end{array}  \end{equation*} 
where the $\, \kappa \,$ coefficients are given by
\begin{align}     \displaystyle
 \label{3.39}  \kappa_{400} \,=\, & \, 8 \,- \, \alpha \,+\, 4 \, \sigma_1^2 \, ( \alpha + 6)  
\,-\, 56 \, \sigma_1 \,  \sigma_4
\,-\, 4 \, (1 - \alpha) \,  \sigma_1 \,  \sigma_6 \, \\  
 \label{3.40}   \displaystyle
\kappa_{220}   \,=\, & -2  \, ( \alpha + 6)   
  \,+\, 8 \, \sigma_1^2 \, ( \alpha + 6)   
\,+\, 56 \, \sigma_1 \,  \sigma_4
\,-\, 8 \, (1 -  \alpha ) \,  \sigma_1 \,  \sigma_6 \, .  
  \end{align}  
%

\bigskip \noindent $\bullet$ \quad 
After these examples where only one partial differential equation is present, we 
consider the case of two (D1Q3), three (D2Q9) or four (D3Q19)
 partial differential equations ``emerging''
from the lattice Boltzmann algorithm. These equations model macroscopic 
conservation of mass and momentum of a linearized fluid in our approach in this
contribution.

\bigskip \noindent $\bullet$ \quad 
{\bf    D1Q3 for athermal linearized Navier--Stokes at   fifth-order }

\smallskip \noindent   
We have in this case two conservation laws ($N=2$ in  relation  (\ref{1.9}))
 and the equilibrium energy is supposed to be given simply by  
\begin{equation} \label{3.7}    
 \epsilon ^{\rm eq}  \,=\,  
\alpha  \, {{\lambda^2}\over{2}} \, \rho   \,. 
\end{equation}  
Due to  (\ref{1.29})  and (\ref{3.7}), the  matrix     $ \, \Psi \,$  for
dynamics relation  (\ref{1.27}) is now given according to 
\begin{equation*} \label{3.8}    
  \Psi   \,=\,   \begin{pmatrix} 1 & 0 & 0 \cr  0  & 1  & 0 \cr 
\alpha \, s \, \lambda^2 / 2 & 0 & 1 - s    \end{pmatrix}   \,, 
\end{equation*}  
and $\, \sigma \,$ is related to parameter $ \, s \,$ according to  (\ref{1.22}):  
$ \,\,  \sigma \equiv {{1}\over{s}} \,-\, {1\over2} \, . $ 
Then equivalent mass conservation at the order 5 looks like equation (\ref{3.4}). We have
precisely: 
\begin{equation}   \left\{ \label{3.9}  \begin{array} {c} \displaystyle 
{{\partial \rho }\over{\partial t}} \,+\,  \displaystyle 
{{\partial q }\over{\partial x}}    \,-\, 
 {{\lambda^2 \, \Delta t^2}\over{12}} \, (1-\alpha)  \,
{{\partial^3 q }\over{\partial x^3}}  \,-\, 
 {{\lambda^4 \, \Delta t^3}\over{12}} \, \alpha \,  (1-\alpha)  \,\sigma\,  
{{\partial^4 \rho }\over{\partial x^4}}  \, 
  \qquad   \qquad   \qquad   \qquad   \\   [2mm] 
   \displaystyle    \qquad   \qquad   \qquad   \qquad 
+ \,  {{\lambda^4 \, \Delta t^4}\over{120}} \,  (1-\alpha)  \, 
\big( 1   \,+\, \alpha  \,+\,10 \, (1 - 2 \,\alpha) \,  \sigma^2 \big  ) \, 
{{\partial^5 q  }\over{\partial x^5}}  \,=\,   {\rm O}(\Delta t^5) \, .      
\end{array}  \right. \end{equation} 
Conservation of momentum takes the form: 
\begin{equation}  \left\{ \label{3.10}  \begin{array} {c} \displaystyle  
{{\partial q }\over{\partial t}} \,+\,  \displaystyle 
\alpha\, \lambda^2\, {{\partial \rho }\over{\partial x}}  \,-\, 
  \lambda^2 \,  \Delta t \,  (1-\alpha)  \, \sigma \,
 {{\partial^2 q }\over{\partial x^2}}  \,  + \,  
 \zeta_3 \, {{\lambda^4 \, \Delta t^2}\over{6}} \, 
    {{\partial^3 \rho }\over{\partial x^3}}  
  \qquad   \qquad   \qquad   \qquad   \\   [2.7mm] 
   \displaystyle    \qquad   \qquad   \qquad   \qquad  
+ \,  \zeta_4 \,  {{\lambda^4 \, \Delta t^3}\over{12}} \,   
 {{\partial^4 q }\over{\partial x^4}}   \,  
\,+ \,    \zeta_5 \, {{\lambda^6 \, \Delta t^4}\over{120}} \,   
 {{\partial^5 \rho  }\over{\partial x^5}}  \,=\,   {\rm O}(\Delta t^5)      
\end{array}  \right.  \end{equation}   
with parameters $\, \zeta_{3} \,$ to   $\, \zeta_{5} \,$  given by  
\begin{eqnarray*}    \label{3.11}   
\zeta_3 \, & = & \,   \alpha  \, (1-\alpha)  \,   (1 \,-\, 6 \, \sigma^2) 
  \\   \label{3.12} 
\zeta_4 \, & = & \,   -    \, (1-\alpha)   \, \sigma \, 
\big( 1 \,-\, 4 \alpha  \,-\, 12 \, (1 - 2\, \alpha) \, \sigma^2   \big) \,
  \\   \label{3.13} 
\zeta_5  \, & = & \,  \alpha  \, (1-\alpha)   \,  
\big( 1 \,-\,  4 \, \alpha 
\,-\, 10 \,  ( 5 - 9 \,  \alpha  ) \,  \sigma^2 
\,+\, 120 \,   ( 2 - 3 \,  \alpha  ) \,  \sigma^4    \big) \, .
\end{eqnarray*}   
When 
 $ \, \sigma = {{1}\over{\sqrt{6}}} , \,  $ 
the coefficient $ \, \zeta_3 \,$ of  relation (\ref{3.10}) is null. 
In this case, the lattice Boltzmann scheme is formally third-order accurate for the
momentum equation. But, as  remarked in \cite{Du09}, 
  the mass conservation  (\ref{3.9}) remains formally second-order
accurate, except for the 
(without any practical interest as it leads to a null viscosity)
case $ \, \alpha = 1 .\,$

\bigskip  \noindent $\bullet$ \quad 
{\bf  D2Q9 for linearized athermal Navier--Stokes at order four }

\smallskip \noindent
The D2Q9 lattice Boltzmann scheme can be used also for simulation of fluid dynamics.
For the particular case of conservation of mass and momentum, 
we just replace matrix $ \, \Psi \,$ of (\ref{3.20})  by the following one, assuming the aim
is to simulate an athermal fluid with speed of sound  $\sqrt{1/3}$: 
\begin{equation*}  \label{3.27}  
   \Psi =    \begin {pmatrix}  1 &  \!\! \!\!0  & \!\!0  & \!\!0 & \!\!0  & \!\!0  
& \!\!0  & \!\!0  & \!\!0 \cr  
0  &  \!\! \!\! 1  & \!\!0  &\!\!0  &  \!\!0  & \!\!0 
&  \!\!0    &  \!\!0  &  \!\!0  \cr 
0   &  \!\! \!\!0  & \!\!1  & \!\!0 & \!\!0  & \!\!0  
& \!\!0  & \!\!0  & \!\!0 \cr 
-2 \, s_3   &  \!\!\!\!0  &  \!\!0  &  \!1 \!\!-\!\! s_3  & \!\!0  & \!\!0  & \!\!0  
& \!\!0  & \!\!0\cr 
  s_4   &  \!\!\!\!0  &  \!\!0  & \!\!0 &    \!1 \!\!-\!\! s_4  & \!\!0  
& \!\!0  & \!\!0  & \!\!0\cr 
0   &   \!\! \!\! \!\! -s_5 / \lambda &  \!\!0     & \!\!0 &  \!\!0    &  \!1 \!\!-\!\! s_5  
& \!\!0  & \!\!0  & \!\!0 \cr 
0   &  \!\!\!\!0  &   \!\! \!\! \!\! -s_5 / \lambda   & \!\!0 &  \!\!0    &  \!\!0  
&  \!1 \!\!-\!\! s_5  & \!\!0  & \!\!0 \cr 
0  &  \!\!\!\!0  &  \!\!0  & \!\!0 &  \!\!0    &  \!\!0  &  \!\!0  
&  \!1 \!\!-\!\! s_7 & \!\!0 \cr 
0  &  \!\!\!\!0  &  \!\!0  & \!\!0 &  \!\!0    &  \!\!0  &  \!\!0  
&  \!\!0  &  \!1 \!\!-\!\! s_7    \end  {pmatrix}  \, . 
\end{equation*} 
We have conservation of mass at fourth-order of accuracy:
\begin{equation}  \label{3.28}   {{\partial \rho}\over{\partial t}} \,+\, 
   {{\partial q_x  }\over{\partial x}}  \,+\,  {{\partial q_y  }\over{\partial y}}  \,-\,  
{{1}\over{18}} \, \lambda^2 \, \Delta t^2 \,  \Delta \, 
\Big(   {{\partial  q_x  }\over{\partial x}}  \,+\, 
 {{\partial q_y }\over{\partial y}}  \Big)+ \, {{\lambda^4 \, \Delta t^3}\over{108}} \, (  \sigma_3 + 
 \sigma_7) \, \Delta^2 \rho   \,=\, {\rm O}(\Delta t^4)  \, 
\end{equation} 
and conservation of two components of momentum: 
\begin{equation}   \label{3.29}  \left\{  \begin{array}{rcl}    \displaystyle
 {{\partial q_x }\over{\partial t}} \,+\, {{\lambda^2}\over{3}} \, 
 {{\partial \rho}\over{\partial x}}  \,-\,  {{\lambda^2}\over{3}} \, \Delta t \, \Big[ \sigma_3 \,  
 {{\partial }\over{\partial x}} \Big(   {{\partial  q_x  }\over{\partial x}}  \,+\, 
 {{\partial q_y }\over{\partial y}}  \Big)
 \,+\,  \sigma_7 \, \Delta q_x  \Big] \,   && 
  \\   [1mm]   \displaystyle   \!\!\!\!
- {{ \lambda^4 \, \Delta t^2}\over{27}} \, 
\Big( 3\, ( \sigma_3^2 + \sigma_7^2)  -1 \Big) \, 
 {{\partial }\over{\partial x}}  \Delta \rho     
   \,-\,  {{ \lambda^4 \, \Delta t^3}\over{108}} \, \bigg( 
\zeta_{40}  \,  {{\partial^4 q_x}\over{\partial x^4}} 
+ \zeta_{31}  \,  {{\partial^4 q_y}\over{\partial x^3 \, \partial y }}  &&
  \\   [1mm]   \displaystyle 
+ \zeta_{22}  \,  {{\partial^4 q_x}\over{\partial x^2 \, \partial y^2 }}
+ \zeta_{13}  \,  {{\partial^4 q_y}\over{\partial x \, \partial y^3 }}  
+ \zeta_{04}  \,  {{\partial^4 q_x}\over{\partial y^4 }} \bigg)  
  &= & {\rm O}(\Delta t^4) 
 \end{array}  \right.  \end{equation} 
\begin{equation}  \label{3.30} \left\{   \begin{array}{rcl}    \displaystyle
 {{\partial q_y }\over{\partial t}} \,+\, {{\lambda^2}\over{3}} \, 
 {{\partial \rho}\over{\partial y}}  \,-\,  {{\lambda^2}\over{3}}\, \Delta t \, \Big[ \sigma_3 \,  
 {{\partial }\over{\partial y}} \Big(   {{\partial  q_x  }\over{\partial x}}  \,+\, 
 {{\partial q_y }\over{\partial y}}  \Big)  
 \,+\,  \sigma_7 \, \Delta q_y  \Big] \,    && 
  \\   [1mm]   \displaystyle   \!\!\!\! 
- {{ \lambda^4 \, \Delta t^2}\over{27}} \, 
\Big( 3\, ( \sigma_3^2 + \sigma_7^2)  -1 \Big) \, 
 {{\partial }\over{\partial y}}  \Delta \rho              
   \,-\,  {{ \lambda^4 \, \Delta t^3}\over{108}} \, \bigg( 
\eta_{40}  \,  {{\partial^4 q_y}\over{\partial x^4}}                 
+ \eta_{31}  \,  {{\partial^4 q_x}\over{\partial x^3 \, \partial y }}   &&
  \\   [1mm]   \displaystyle 
+ \eta_{22}  \,  {{\partial^4 q_y}\over{\partial x^2 \, \partial y^2 }}
+ \eta_{13}  \,  {{\partial^4 q_x}\over{\partial x \, \partial y^3 }}  
+ \eta_{04}  \,  {{\partial^4 q_y}\over{\partial y^4 }} \bigg)  
  & = & {\rm O}(\Delta t^4) 
\end{array}   \right.  \end{equation} 
where the coefficients $ \, \zeta \,$ are given by  
\begin{equation} \label{3.31}  \left\{ \begin{array}{rcl}
 \zeta_{40} & \,=\, &  \eta_{04} \,=\,   -  \, \sigma_3 \, -  \, \sigma_7 
\,  - \, 12  \, \sigma_3^2 \, \sigma_7 
\,  - \, 12  \, \sigma_3 \, \sigma_7^2
\,  +  \, 18 \, \sigma_3^2  \, \sigma_5 \\  
&& \,  +  \,  6 \, \sigma_5 \sigma_7^2   \,
\,  - \,   12  \, \sigma_3 \, \sigma_4 \, \sigma_5
\,  - \,   24  \, \sigma_3 \, \sigma_5 \, \sigma_7
\, +  \,  12  \, \sigma_4 \, \sigma_5 \, \sigma_7  \,    \\  [1mm] 
 \zeta_{31} & \,=\, &  \eta_{13} \,=\, 
 - \, 4 \, \sigma_3 \, - \, 7 \, \sigma_7 
  \,  +  \, 18 \, \sigma_3^2  \, \sigma_5 
  \,  +  \, 18 \, \sigma_5  \, \sigma_7^2 
\,  - \, 12  \, \sigma_3^2 \, \sigma_7  \\
&& \,  - \, 12  \, \sigma_3 \, \sigma_7^2 
\,  - \,   12  \, \sigma_3 \, \sigma_4 \, \sigma_5
\, +  \,  12  \, \sigma_3 \, \sigma_5 \, \sigma_7  
\, +  \,  12  \, \sigma_4 \, \sigma_5 \, \sigma_7  
\, +  \,  12  \, \sigma_7^3    \\  [1mm]  
 \zeta_{22} & \,=\, &   \eta_{22} \,=\,  
 - \, 13 \, \sigma_3   \,  + \,  6 \, \sigma_4   \,  - \,  10 \, \sigma_7  
\, +  \,  18  \, \sigma_3^2   \, \sigma_5  
\,  - \,   12  \, \sigma_3^2 \, \sigma_7  
 \,  - \,   12  \, \sigma_3 \, \sigma_7^2    \\ 
&& \, +  \,  30  \, \sigma_5   \, \sigma_7^2  
\,  - \,   12  \, \sigma_3 \, \sigma_4 \, \sigma_5
  \, +  \,  120  \, \sigma_3 \, \sigma_5 \, \sigma_7  
\, -  \,  60  \, \sigma_4 \, \sigma_5 \, \sigma_7  
\, -  \,  12  \, \sigma_7^3     \\  [1mm]  
\zeta_{13} & \,=\, &  \eta_{31} \,=\, 
 - \, 10 \, \sigma_3   \,  + \,  6 \, \sigma_4   \,  - \,  7 \, \sigma_7  
\, +  \,  18  \, \sigma_3^2   \, \sigma_5  
\,  - \,   12  \, \sigma_3^2 \, \sigma_7 
  \,  - \,   12  \, \sigma_3 \, \sigma_7^2  \\ 
&&  \, +  \,  18  \, \sigma_5   \, \sigma_7^2  
\, +  \, 12  \, \sigma_3 \, \sigma_4 \, \sigma_5
\, +  \,  84  \, \sigma_3 \, \sigma_5 \, \sigma_7
  \, -  \,  60  \, \sigma_4 \, \sigma_5 \, \sigma_7  
\, +  \,  12  \, \sigma_7^3   \\  [1mm] 
  \zeta_{04} & \,=\, &  \eta_{40} \,=\, 
 - 3 \, \sigma_7   \, +  \, 24 \,  \sigma_5 \,  \sigma_7^2  
\, - \,   12 \, \sigma_7^3 \, . 
\end{array} \right.  \end{equation}
%

  \bigskip \noindent $\bullet$ \quad 
{\bf  D3Q19 for linearized Navier--Stokes } 

\noindent 
The D3Q19 Lattice Boltzmann scheme is described with details 
{\it e.g.} in J.~T\"olke {\it et al} 
 \cite{TKSR02}. The construction of matrix $M$ that parameterizes 
the transformation (\ref{1.8})  is presented in full  detail  
in  relations (\ref{3-kinetic}) to  (\ref{3-third}) 
in  Appendix~B.  
The associated matrix $\, \Psi \,$ is also of order 19 
and therefore quite difficult to write on a A4 paper sheet.
Due to constitutive relations (\ref{1.27}) and  (\ref{1.28}), it is
easily obtained from the expression of equilibrium moments. We have taken for this D3Q19
 scheme 
\begin{equation} \label{d3q19-equi}   \left\{ \begin{array}{rcl}  \displaystyle 
m_4^{\rm eq} &=&  \theta \, \lambda^2 \\
m_5^{\rm eq} &=& m_6^{\rm eq}  =  m_7^{\rm eq} =  m_8^{\rm eq}  =  m_9^{\rm eq} = 0  \\ 
m_{10}^{\rm eq}  &=&  m_{11}^{\rm eq} =  m_{12}^{\rm eq} =  0 \\ 
m_{13}^{\rm eq} &=& \beta \lambda^4  \\
m_{14}^{\rm eq}  &=&  m_{15}^{\rm eq}  =  0 \\ 
m_{16}^{\rm eq}  &=&   m_{17}^{\rm eq} =  m_{18}^{\rm eq} =  0  \, . 
\end{array} \right.    \end{equation}  
In order to obtain physical equations at first-order of accuracy with a sound velocity
$ \, c_0 \, $  given by
%
$ \,  c_0  =  \alpha \, \lambda  \, $ 
%
the relation 
$ \,\theta  \, = \,  57 \, \alpha^2  - 30  \, $ 
must be imposed to obtain correct fluid second-order partial differential equations 
and the parameter $ \, \beta \, $ remains free.

\bigskip \noindent $\bullet$ \quad  
When the number of velocities of the Boltzmann scheme is reduced (up to D2Q9 scheme  typically),
it is possible to expand the dispersion equation formally and to derive equivalent partial
differential equations up to an arbitrary order. 
We have done the comparison for one-dimensional and bi-dimensional schemes.  
%
%
The process has been extended to models with more velocities and various
conserved quantities; however the equations
become very complicated and thus will not be reproduced here.
Let us just
mention that the expressions found 
are quite similar to those obtained for the
previous test cases.


\section{The fourth-order accurate  lattice Boltzmann scheme }
 
\smallskip \noindent $\bullet$ \quad 
In this section, we precise how to choose particular ``quartic''
values of relaxation parameters in order to increase the  accuracy of the  scheme. 
We verify with the help of precise numerical experiments for analytical test cases that
the numerical precision follows our prediction.
We focus first on classical thermics at two and
three space dimensions. 
Then we propose two numerical experiments for athermal linearized Navier-Stokes equations
at two and three space dimensions for a nontrivial geometry.

 \bigskip 
\noindent $\bullet$ \quad 
{\bf  The D2Q5   lattice Boltzmann scheme for a thermal problem }

\noindent 
We obtain the order 4 by setting $ \, \kappa_{40} = 0 \,$  and $ \, \kappa_{22} = 0 \,$ 
in relations (\ref{3.17})  and  (\ref{3.18})  respectively. We obtain :
\begin{equation}  \label{4.1}   
 \sigma_3   \,=\,    \displaystyle  \sigma_1 \, {{ \alpha + 4}\over{1 - \alpha}}   \,-\,  
 {{1}\over{12 \, \sigma_1}} \,  {{2 + 3\, \alpha}\over{1 - \alpha}} \,,\qquad   
\sigma_4   \,=\,   \displaystyle {{1}\over{6 \, \sigma_1}} \, .     
 \end{equation} 
The BGK condition  $ \, \sigma_1 = \sigma_3 = \sigma_4 \,$ leads to 
 $\, \sigma_1= \frac{1}{\sqrt{12}} \,$ and   $ \, \alpha = -4 \, $ and
thus to a thermal diffusivity equal to 0. 
Note that the intermediate TRT presented in Ginzburg {\it et al}   \cite{GVD08a, GVD08b} 
supposes simply   $ \,  \sigma_3 = \sigma_4 . \,$ If we insert this constraint inside
 relations (\ref{4.1}), we get 
\begin{equation*} 
 \sigma_1  \,=\,  \frac{1}{\sqrt{12}} \,,\qquad  
 \sigma_3  \,=\,  \frac{1}{\sqrt{3}}  
 \end{equation*} 
to enforce fourth-order accuracy. 
Then the  d'Humi\`eres version of lattice Boltzmann scheme
is mandatory for this improvement of the method with a wide family of admissible
parameters. 
In order to study the fourth-order accuracy of the D2Q5  lattice Boltzmann scheme for
thermal problem, we use three different approaches.
The first two  consider the interior scheme and the third one incorporates  boundary
conditions.

\begin{figure}[!t] 
\centerline {\includegraphics[width=.76  \textwidth]{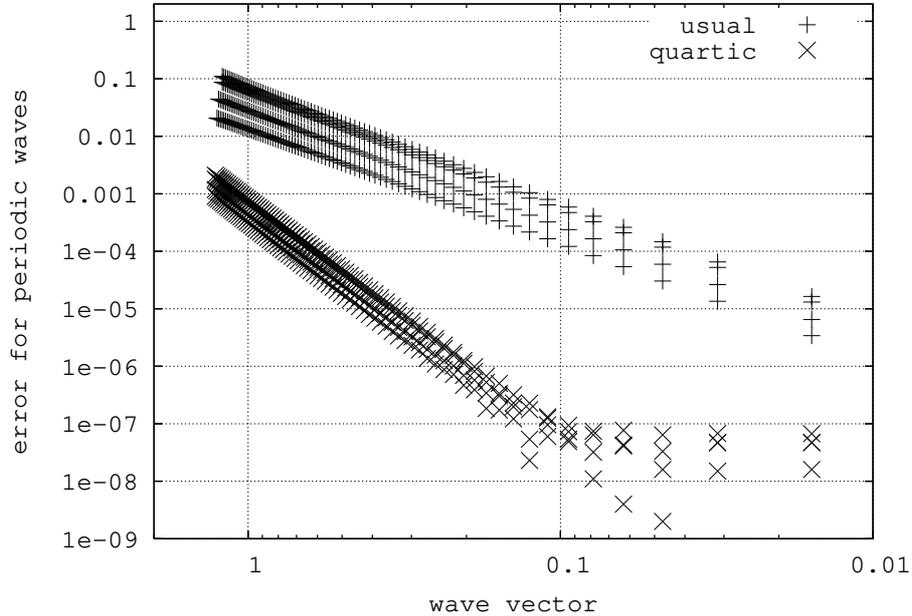}} 
 \caption {  Error $ \mid \!  \frac { \Gamma_{\rm num}} { \Gamma_{\rm th}} -1  \! \mid  $ 
 of D2Q5 scheme for thermic test case, ``one point''
  simulation. Different curves correspond to different orientation of the wave-vector
with respect to the axis, showing the angular dependence of the next order.    } 
\label{f-therm-unpt}
\end{figure}

\bigskip \noindent $\bullet$ \quad 
First of all, we study homogeneous plane waves with a ``one point computation''. 
In that case, we can derive numerically a dispersion equation for scheme 
(\ref{1.21})  associated
with (\ref{1.8}),  (\ref{1.27}),  (\ref{7.14})  and (\ref{3.15}), 
as proposed in \cite{LL00}. 
We introduce a wave in the Boltzmann scheme, {\it id est} 
%
$ \,\, f(x,\,t)~\equiv~\widehat{f} (k_x,\, k_y) \,\, {\rm exp} \big(   
i \, k_x \, x \,+ \, i \, k_y \, y \big)   . \,\,$ 
%
Then we have
%
$ \,\, f(x, \, t+\Delta t) \,=\, G \, f(x,\,t)  \,\,$ 
%
with the so-called amplification matrix $\, G \, $ 
(see {\it e.g.} Richtmyer and Morton \cite{RM67})
obtained without difficulty from matrices $M$, $\Psi$ and $B$ defined 
respectively in (\ref{7.14}) (\ref{3.15}) and 
\begin{equation*}  \label{4.5}     
B \,=\, { \rm diag} \, \Big( 1 ,\,  {\rm e}^{   \displaystyle i \, k_x \, \Delta x }  ,\, 
{\rm e}^{   \displaystyle i \, k_y \, \Delta x } ,\,
{\rm e}^{   \displaystyle -i \, k_x \, \Delta x } ,\,
{\rm e}^{   \displaystyle -i \, k_y \, \Delta x } \Big) \, 
\end{equation*}  
for the D2Q5 scheme displayed in Figure \ref{fd2q5q9} (left). 
Then $ \, G = B \, M^{-1} \, \Psi \, M \, . $  
Then if $ \, {{\partial}\over{\partial t}} \,$ is formally given by  relation (\ref{3.16})
and operators   $ \, {{\partial}\over{\partial x}} \,$ and  
 $ \, {{\partial}\over{\partial y}} \,$ replaced by $\, i \, k_x \,$ and 
 $\, i \, k_y \,$ respectively,  the number 
$ \, z =  {\rm exp } \big(  \Delta t \, 
 {{\partial}\over{\partial t}} \big)   \, $  
is an eigenvalue of matrix $G$ at  fourth-order of accuracy. 
The numerical experiment (see Figure \ref{f-therm-unpt})
 confirms the theoretical development of the dispersion equation. 
Note that for situations relaxing to uniform state, the eigenvalues that we determine below
are negative;  however we shall 
express results in terms of 
{\it positive} relaxation rates with adequate
sign changes.
 
\bigskip \noindent $\bullet$ \quad 
For inhomogeneous situations, with $\, N_{\cal L} \, $ lattice points (and 
$\,   N_{\cal L} \, (J+1) \, $ degrees of freedom), 
one can study the time evolution starting from
some initial state. Another approach for linear situations
considers that the state $\, X(t) \, $ that belongs to 
$\, \R^{ N_{\cal L} \, (J+1)} \, $ 
can be decomposed as a sum of eigenmodes of the operator $ \, A \, $ 
defined using  the discrete evolution scheme:
\begin{equation}  \label{4.8}  
X(t + \Delta t) \,\equiv \, A \, \smb \, X(t) \,.   
\end{equation}  
 The matrix $ \, A \, $ being of very large size, one can look for some 
of its eigenmodes using for instance the method proposed by Arnoldi \cite{Ar51}. 
To accelerate the Arnoldi computations, following a suggestion by L. Tuckerman
\cite{Tu02},  we replace
the determination of the eigenvalues of equation (\ref{4.8}) 
by the determination of the eigenvalues of
\begin{equation}  \label{4.9}   
X(t + (2 \ell +1)\,\Delta t) \,\equiv \, A^{2 \ell +1} \, \smb \, X(t)   \,,    
\end{equation}  
for some $ \ell \in  \N$,  
using the fact that the lattice Boltzmann scheme is very fast compared to the inner
``working''  of the Arnoldi procedure. Replacing problem  (\ref{4.8})  by problem
 (\ref{4.9}) not only
increases the splitting between various eigenmodes, but  also helps to discriminate
against the acoustic modes by multiplying the 
logarithm of the  imaginary part of the eigenvalues
by $2 \ell +1$. Note that choosing an even number of time steps would bring in
the ``checker-board'' type modes. We denote by $\, \Gamma_{\rm num} \,$ 
any eigenvalue computed with this methodology. 

\begin{figure}[!h] 
\centerline {\includegraphics[width=.76  \textwidth]{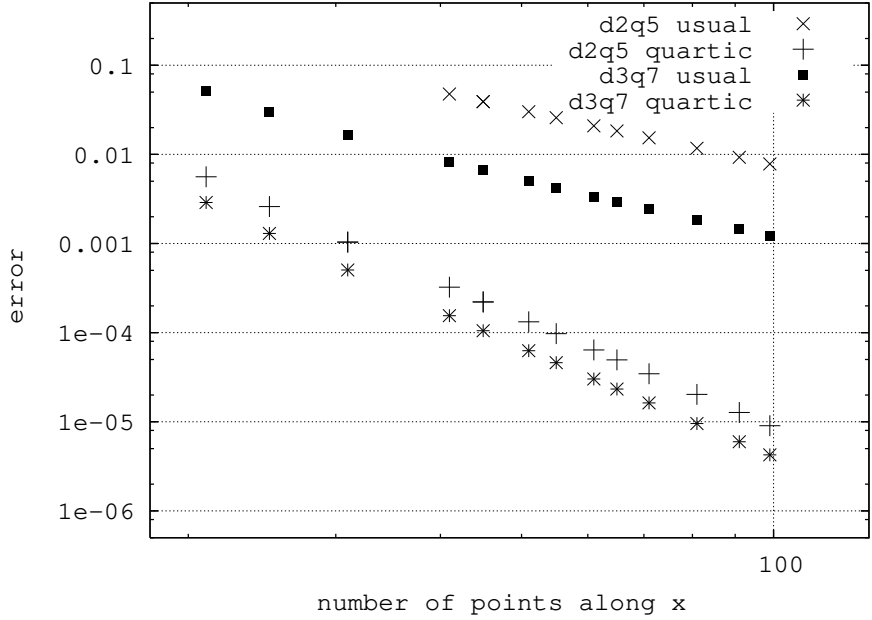}} 
 \caption {  Arnoldi test case for periodic thermics, $I_x=5$,  $I_y=0$.   
  Various parameters for lattice Boltzmann schemes D2Q5 and D3Q7.     } 
\label{f-arno-therm}
\end{figure}

\bigskip \noindent $\bullet$ \quad 
We first test this method for  ``internal''  lattice, 
{\it id est} with  a periodic $\,  N_{\cal L} \equiv N_x \times  N_y  \, $ 
situation and find the same results
as those derived from the ``one-point'' analysis 
(see Figure \ref{f-therm-unpt}) with very good accuracy. 
For this periodic situation, the eigenmodes are plane waves for the 
wave vector
%
$ \,  k_x = 2 \pi {{I_x}\over{N_x}} , \, 
 k_y =  2 \pi {{I_y}\over{N_y}} , $  
%
 where $I_x$ and $I_y$ are integers.
We compare the numerical relaxation rates $\Gamma_{\rm num}(I_x,I_y,N_x,N_y)$ 
to $\Gamma_{\rm th} \equiv \kappa (k_x^2+k_y^2)$
and show in Figure  \ref{f-arno-therm} the relative difference 
between those two quantities (called the ``error'')
for the particular values $I_x=5$ and $I_y=0$ and $N_x$ from 11 to 91. With 
arbitrarily chosen values of the ``non-hydrodynamic'' $s$-parameters, we observe
second-order convergence. However for the quartic $s$-parameters the convergence
is of order four with a large decrease in the absolute value of the error.
Analogous results are displayed in Figure \ref{f-arno-therm}  for D3Q7.
 
\bigskip \noindent $\bullet$ \quad 
We now consider a second case with boundary conditions: 
exact solution for the modes of the Laplace equation in a circle of radius $R$
with homogeneous Dirichlet boundary conditions.
Density is defined with (\ref{1.10})
applied with $J=4$ in this particular case.
Recall that density follows  heat equation
%
$ \,\, {{\partial \rho}\over{\partial t}} \,-\, \kappa \, \Delta \rho =  0 \,\, $
%
with
$\,  \kappa =   {{ \lambda^2 \, \Delta t }\over{10}} \, \sigma_1 \, (4 + \alpha) \,  $
and   homogeneous boundary conditions at $r=R$. 
The solution of this problem is standard (see {\it e.g.} Landau and Lifchitz \cite{LL59}, 
 Abramowitz and Stegun  \cite{AS65} 
or   Carslaw and Jaeger \cite{CJ67})  
and is parameterized by a pair $ \,  (\ell,\, n) \,$ of integers. We introduce 
 $ \,\zeta_\ell^n  $, 
the $ n^{\rm th} \,$ zero of the Bessel function $\, J_\ell .  \,$
Then a solution with time dependence as exp($-\Gamma t$) defines 
a   corresponding eigenvalue $\, \Gamma \,$ (also denoted as 
$\Gamma_{\rm th}$ in the following)  that satisfies  
\begin{equation}  \label{4.17}
\Gamma \,=\,  \kappa \,  \Big(  {{ \zeta_\ell^n } \over{ R }}  \Big)^2 \, .    \,
\end{equation}

\begin{figure}[!t]  
\centerline {\includegraphics[width=.76  \textwidth]{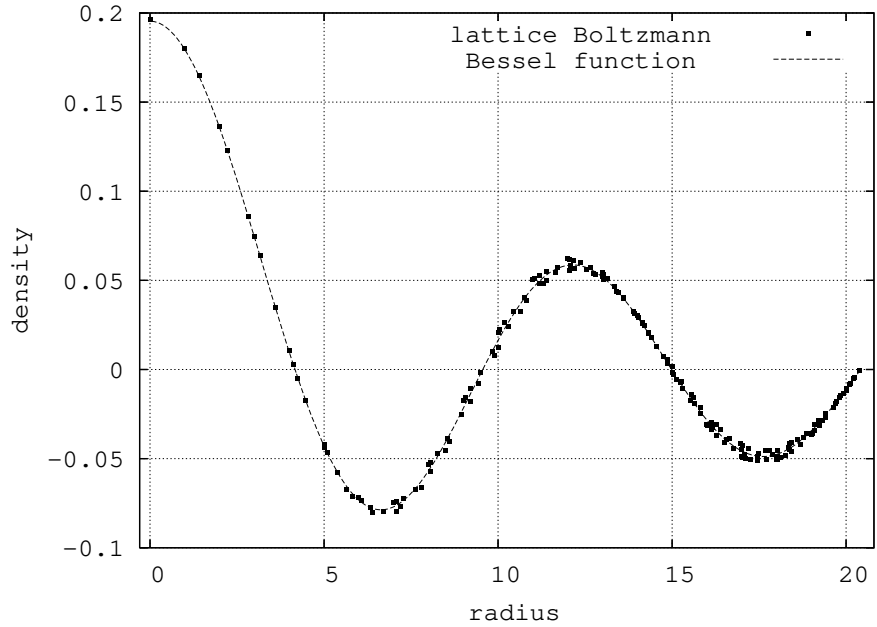}} 
 \caption { D2Q5 scheme for thermics inside a circle.  Eigenmode   $n=4$,  $\ell=0$ 
 for  heat equation with  Dirichlet boundary conditions. 
Second-order accuracy with  usual parameters for  lattice Boltzmann scheme.  } 
\label{f-d2q5-th-n4m0-u}
\end{figure}

\begin{figure}[!h]   
\centerline {\includegraphics[width=.76  \textwidth]{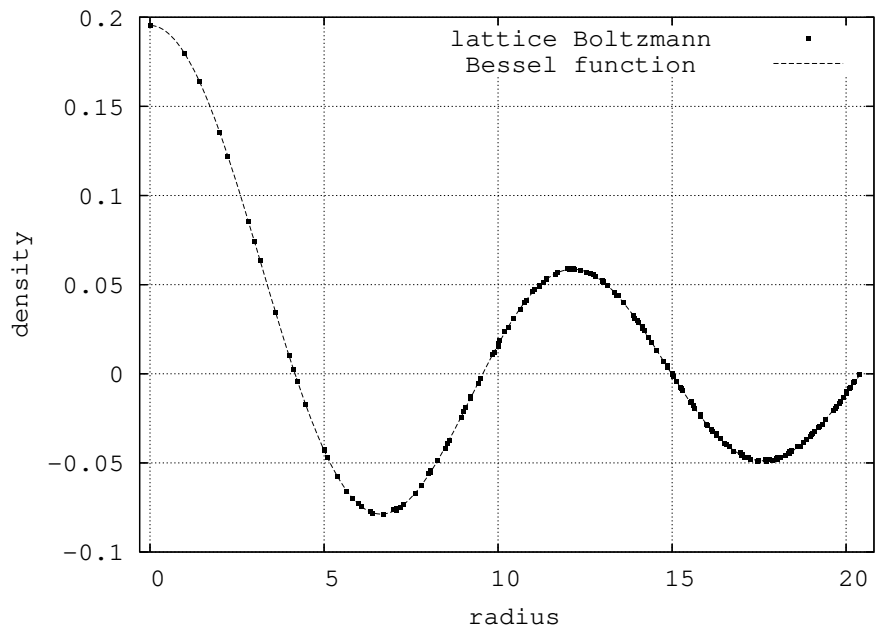}} 
 \caption {  D2Q5 scheme for thermics inside a circle.   Eigenmode   $n=4$, $\ell=0$.
  Quartic parameters for   lattice Boltzmann scheme.    } 
\label{f-d2q5-th-n4m0-q}
\end{figure}

\begin{figure}[!t]   
\centerline {\includegraphics[width=.72  \textwidth]{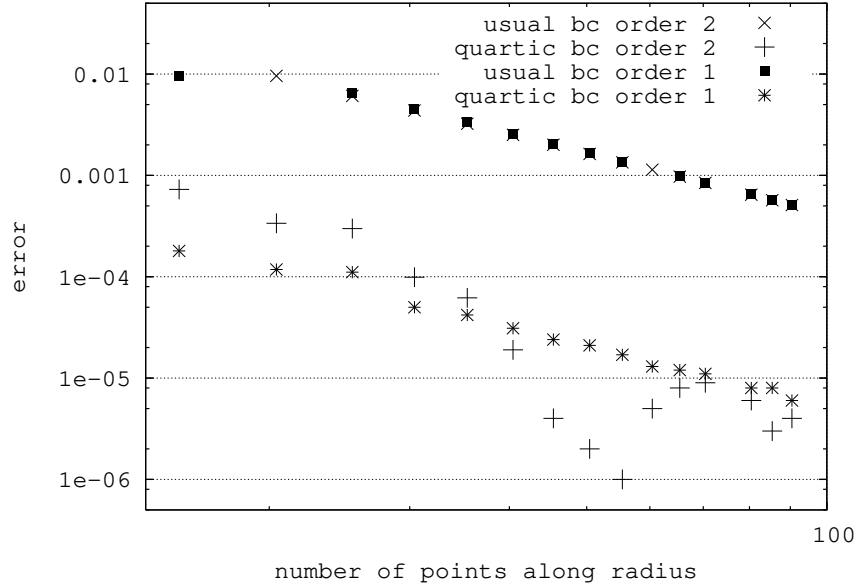}} 
 \caption { D2Q5 scheme for thermics in a circle.  
  Eigenmode   $n=1$, $\ell=5$.    Errors for  various parameters 
 for lattice Boltzmann  and boundary  schemes. }
\label{f-d2q5-th-n1m5}
\end{figure}

\bigskip \noindent $\bullet$ \quad 
The effect of fourth-order accuracy Boltzmann scheme in computing the eigenfunction 
is spectacular: just compare Figures \ref{f-d2q5-th-n4m0-u} and  \ref{f-d2q5-th-n4m0-q}. 
Nevertheless, the effect of boundary conditions 
(we use anti-bounce-back with interpolation {\it \`a la} Bouzidi  {\it et al}   \cite{bfl})
cannot be neglected. In Figure \ref{f-d2q5-th-n1m5}, we have compared 
the error defined by $ \, \mid \! {{\Gamma_{\rm num}}\over{\Gamma_{\rm th}}} -1 \! \mid \, $ 
for two internal schemes (with
usual and quartic parameters) and two versions (first and second-order) of 
simple numerical  boundary conditions introduced by Bouzidi {\it et al}   \cite{bfl}. 
We still observe a better numerical 
precision of the schemes (by two orders of magnitude typically) whereas the 
convergence still remains second-order accurate. 
We conclude that the effect of boundary conditions is crucial for the determination of the
order of convergence. Nevertheless, the choice of quartic parameters gives a higher
precision for the lattice Boltzmann scheme.

%
\bigskip   \noindent $\bullet$ \quad 
{\bf The  D3Q7 lattice Boltzmann scheme for a thermal problem } 

\noindent 
We obtain the order 4 by setting $ \, \kappa_{400} = 0 \,$  and $ \, \kappa_{220} = 0 \,$ 
in relations (\ref{3.39}) and  (\ref{3.40}). We obtain : 
\begin{equation*}    \displaystyle 
\sigma_4  \,=\,   \displaystyle  {{1}\over{6 \, \sigma_1}}  \,, \qquad 
 \sigma_6   \,=\,    \displaystyle  {{\alpha + 6}\over{1 - \alpha}} \,  \sigma_1 \,-\, 
{{4 \,+\, 3 \, \alpha}\over{12 \, (1 - \alpha) }} \, {{1}\over{\sigma_1}} \, .  
\end{equation*}
As for D2Q5, the ``BGK condition''  $\, \sigma_1=\sigma_4=\sigma_6 \,$ leads to 
$\, \sigma_6 = {{1}\over{\sqrt{6}}} \,$ and $ \alpha=-6$ and thus to thermal diffusivity equal to 0.
%
Theoretical modes of the Laplace equation in a sphere of radius $R$ with homogeneous
Dirichlet boundary conditions are parameterized through   the $n^{\rm th}$ zero
$ \,\eta_{\ell + 1/2} ^n \,  $ of semi-integer Bessel function $\,  J_{\ell + 1/2} \,$ 
and  the eigenvalue $ \Gamma  $ is given by: 
\begin{equation}      \label{4.26}   \displaystyle
\Gamma \,=\,  \kappa \,  \Big(  {{  \eta_{\ell + 1/2} ^n } \over{ R }}  \Big)^2
\,, \qquad  \ell \in \N \,, \quad n \geq  1    \,.
\end{equation}

\begin{figure}[!t]   
\centerline {\includegraphics[width=.76  \textwidth]{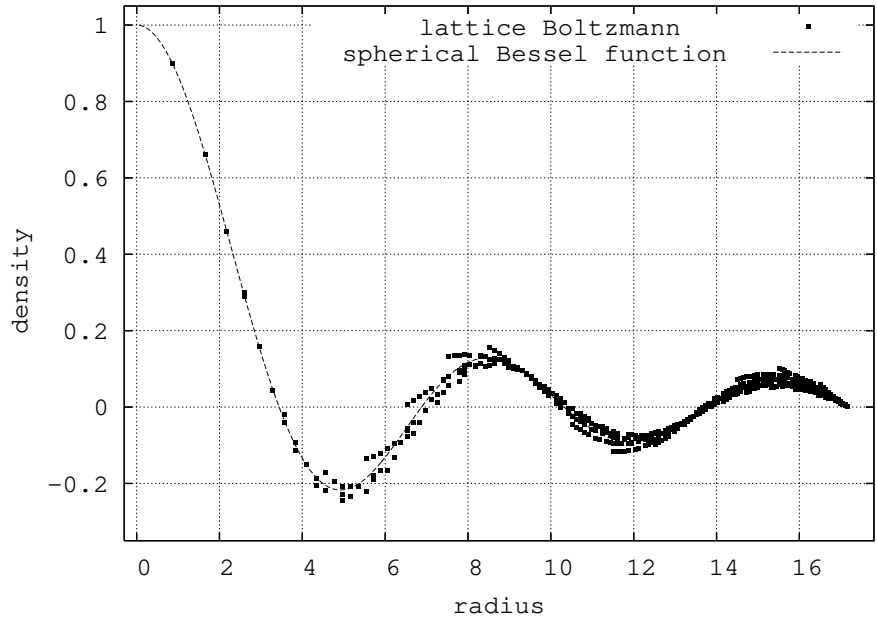}} 
 \caption {  D3Q7 lattice Boltzmann scheme for thermics in a sphere.
 Eigenmode $n=5$,   $\ell=1$,   $m=0$   with usual parameters.   }
\label{f-d3q7-th-u}
\end{figure}

\begin{figure}[!h]    
\centerline {\includegraphics[width=.76  \textwidth]{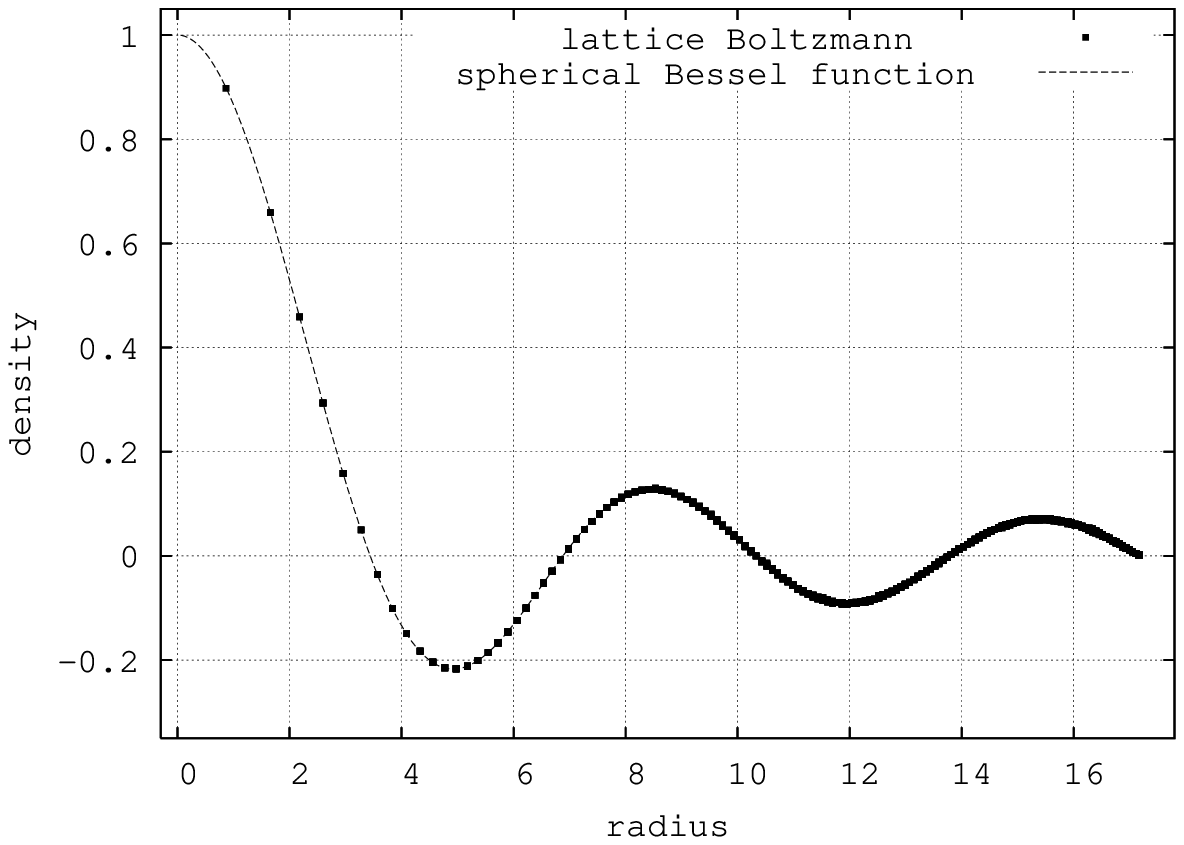}} 
 \caption {   D3Q7 lattice Boltzmann scheme for thermics in a sphere.
 Eigenmode $n=5$,   $\ell=1$,   $m=0$   with quartic  parameters.    }
\label{f-d3q7-th-q}
\end{figure}

\bigskip \noindent $\bullet$ \quad 
The results of Figures \ref{f-d3q7-th-u}  and \ref{f-d3q7-th-q} 
have been obtained with $ \, R = 17.2 \, $ 
and  $ \, n = 5.  \, $ The theoretical value of the eigenvalue is 
$ \, \Gamma = 5^2 \, \pi^2 \, \kappa / R^2 $ (as for $m=0$, the zeros of the
semi-integer Bessel function are simply $\pi \, n $).
We have used 
parameters 
$ \, s_1 = 1.26795 , \, $ $  s_4 = 1.2  , $ $ \,  s_6 = 1.3 \, $ 
for the usual computations. 
%
%
The quartic parameters have been chosen as 
\begin{equation*}      \label{4.28}   \displaystyle       
s_1 = 1.26795 \,, \quad s_4 = s_6 = 0.92820 \, .  
\end{equation*}  
%

\begin{figure}[!t]     
\centerline {\includegraphics[width=.76  \textwidth]{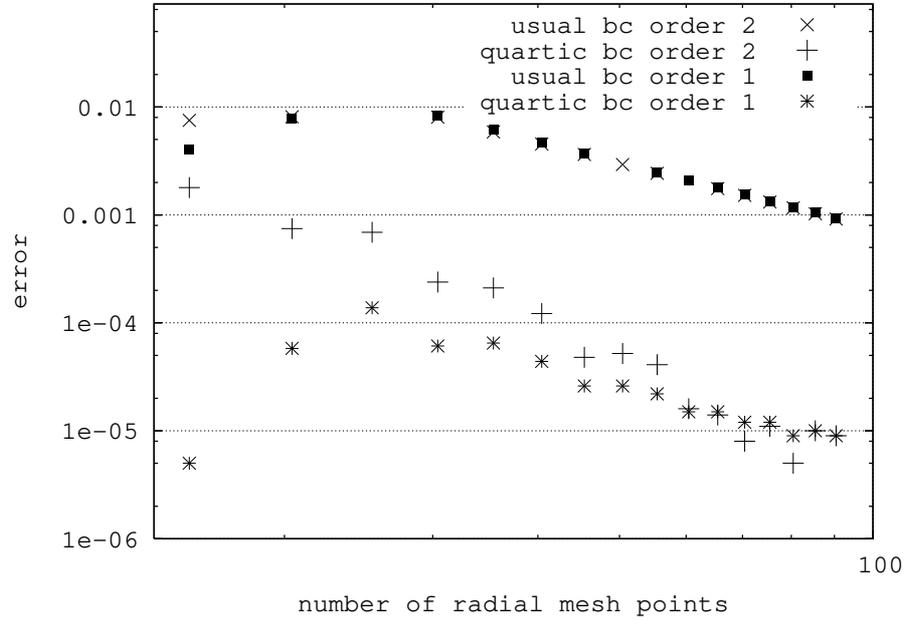}} 
 \caption {    D3Q7  scheme for thermics in a sphere  with   Dirichlet
 boundary conditions. Eigenmode   $n=1$,    $\ell=0$.  
 Errors for  various parameters for lattice Boltzmann   and boundary  schemes. } 
\label{f-d3q7-th-err-mode}
\end{figure}

\smallskip \noindent 
>From results presented in Figure \ref{f-d3q7-th-err-mode},  
the conclusion is essentially the same as that observed for two-dimen\-sional thermics: 
the results are improved by two orders of magnitude typically, but the 
rate of convergence
cannot be rigorously measured or  still remains of second-order. 
 
\begin{figure}[!h]      
\centerline {\includegraphics[width=.76  \textwidth]{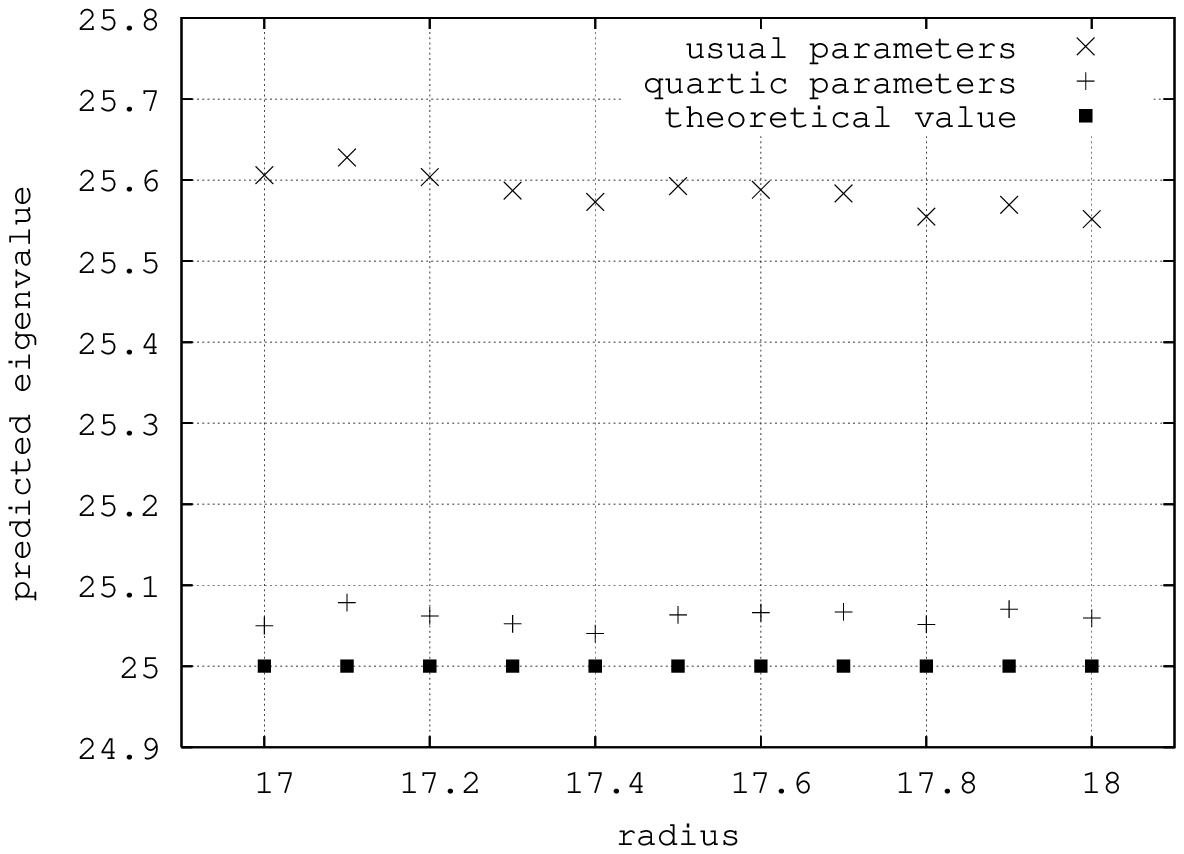}}  
 \caption {  D3Q7 for thermics in a sphere. Eigenmode (in units $\kappa \pi^2 /R^2$) 
for  $n=5$ and   $\ell=0$. Variation of the location of the boundary  
between $R = 17$ and $ R = 18$.    } 
\label{f-d3q7-th-mode-cl}
\end{figure}     

\bigskip \noindent $\bullet$ \quad 
We also made a parameter study of the location of the boundary condition. 
We plot in Figure \ref{f-d3q7-th-mode-cl} 
the ratio  $   \, \smash {   {{\Gamma \, R^2}/({\kappa \, \pi^2 })} } $
with $\, \Gamma \, $ given by  relation (\ref{4.26}). 
We use Bouzidi {\it et al} \cite{bfl} 
boundary procedure with linear interpolation. 
 The fluctuation due to the boundary algorithm is around 0.2 \%. 
The gap between second-order usual computation and new fourth-order computation is of the
order of 2\%. We observe that this gap is 
{\bf one order of magnitude larger}  than the error due to the
choice of the boundary condition estimated from the fluctuations with the imposed radius.


\bigskip \noindent $\bullet$ \quad 
{\bf  The D2Q9 for linearized athermal Navier--Stokes at   fourth-order }

\smallskip \noindent 
We consider now the linear fluid model obtained by a D2Q9 lattice Boltzmann scheme.
The equivalent partial differential equations are given at the order 4 by relations  (\ref{3.29})
to (\ref{3.31}). The dream would be to enforce high order accuracy. However, this is
definitively impossible in the framework considered here due to the never null third-order term
for mass conservation (\ref{3.28}). Recall notation  (\ref{1.9}) for conservative variables:
%
$ \,\, W   \equiv  ( \rho ,\, q_x ,\, q_y )^{\rm  \displaystyle t} \, \,$ 
%
and write the equivalent equations   (\ref{3.29})-(\ref{3.31}) in the synthetic  form:
\begin{equation}      \label{4.30}   \displaystyle        
\partial_t W_k \,+\, \sum_{ j,\, p,\,  q } \,  A_{k p q}^j \, 
  \partial_{x}^p  \partial_{y}^q W_j   \,=\, {\rm O}(\Delta t^4) \,.  
\end{equation} 
We search for a  dissipative mode, {\it id est} a mode for linear incompressible Stokes 
problem under the form
%
$ \,\,   W(t) =  {\rm e}^{   - \Gamma t \,+\, i (k_x \, x 
+ k_y\, y)  } \, \, \widetilde{W} \, .  \,$   
%
Then $\Gamma$ is an eigenvalue of the matrix $A$ defined by  
\begin{equation*}      \label{4.32}   \displaystyle   
A_{k}^j  \,=\,      \sum_{ j,\, p,\,  q } \, 
 A_{k p q}^j \,\, \, (i\, k_x)^p \,(i\, k_y)^q \, . 
\end{equation*} 
We know (see {\it e.g.} Landau and Lifchitz \cite{LL59}) that for  Stokes problem 
(incompressible shear modes),  the relation  
\begin{equation}      \label{4.33}   \displaystyle    
\Gamma \,=\,  \nu \, \big( k_x^2 \,+\, k_y^2 \big)   \,\, 
\end{equation} 
is classical.   Moreover, as a consequence of     (\ref{3.29}) and  (\ref{3.30})
\begin{equation}      \label{4.34}   \displaystyle     
 \nu \,=\, {{ \lambda^2}\over{3}} \, \Delta t \, \sigma_7  \,\,  
\end{equation} 
for a lattice Boltzmann scheme with multiple relaxation times. 

\bigskip \noindent $\bullet$ \quad 
We propose here to 
tune the  parameters $ \, s_{\ell} \, $ in such a way that 
the relation (\ref{4.33})    is enforced for the modes of (\ref{4.30}).
Precisely, we search $ \, s_{\ell} \, $ such that 
\begin{equation}      \label{4.35}   \displaystyle    
 \Delta_m   \, \equiv\,{\rm det} \, \bigg[ A \,-\,  \Big( 
{{ \lambda^2}\over{3}} \, \Delta t \, \sigma_7 \Big) 
 \,  \big( k_x^2 \,+\, k_y^2 \big)  \, {\rm Id} \bigg] \,=\,  
  {\rm O}(\Delta t^7)  \, . 
\end{equation} 
With an elementary formal computation, the third-order term 
$ \,  \Delta_m^3 \, $  of $ \, \Delta_m \,$ relative
to $ \, \Delta t \, $ is equal to  
\begin{equation} \label{4.36}  \left\{ \begin{array}{c} \displaystyle  
 \Delta_m^3  \,=\, -{{\Delta t^3 \,  \lambda^6}\over{108}} \, \sigma_7 \, 
 \big( k_x^2 \,+\, k_y^2 \big)  \,
\Big(   ( -1 \,-\, 4 \,  \sigma_7^2  \,-\, 8 \, \sigma_5  \, \sigma_7 ) 
\,  \big( k_x^4 \,+\, k_y^4 \big)  
\qquad \qquad     \\  \qquad   \qquad   \qquad  \qquad   \qquad       \, 
+ \, 2 \, (1 \,-\,  4 \, \sigma_7^2 \,-\,  4 \, \sigma_5 \, \sigma_7 ) 
\,   k_x^2 \,   k_y^2  \Big) \, . 
\end{array} \right. \end{equation} 
It is then clear that the expression (\ref{4.36}) 
 is identically null for parameters $ \, \sigma_5 \, $ and 
 $ \, \sigma_7 \, $ chosen according to 
\begin{equation}      \label{4.37}   
\sigma_5 \,=\, {{\sqrt{3}}\over{3}} \,, \quad 
\sigma_7 \,=\, {{\sqrt{3}}\over{6}} \, . 
\end{equation} 
With this particular choice of parameters, so-called  quartic in what follows, the
viscosity $ \nu $ in relation  (\ref{4.34})  has the following particular value: 
\begin{equation}      \label{4.38}    
 \nu \,=\,  {{ \lambda^2 \,  \Delta t}\over{\sqrt{108}}} 
\approx \, 0.096225 \,\,  \lambda^2 \,  \Delta t \, . 
\end{equation} 
Then it is very simple to verify that the determinant 
$ \, \Delta_m \, $ is null up to terms of order seven and relation  (\ref{4.35}) is satisfied. 

\begin{figure}[!t]   
\centerline {\includegraphics [width=.73  \textwidth]  {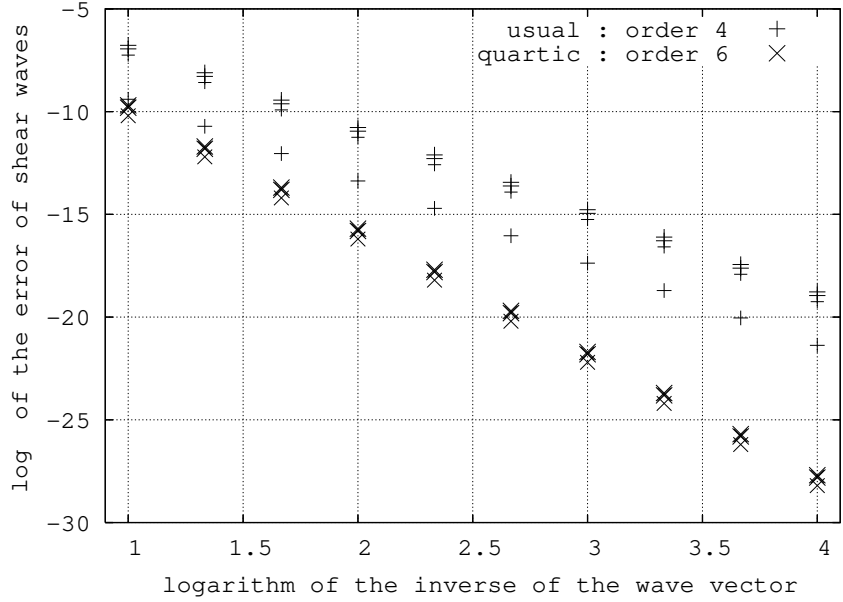}}
 \caption { D2Q9 ``one point'' test case of shear waves for different angles of the wave
 vector.  } 
\label{f-d2q9-unpt}
\end{figure}  
\begin{figure}[!h]       
\centerline {\includegraphics [width=.60 \textwidth,  height=.60 \textwidth] 
{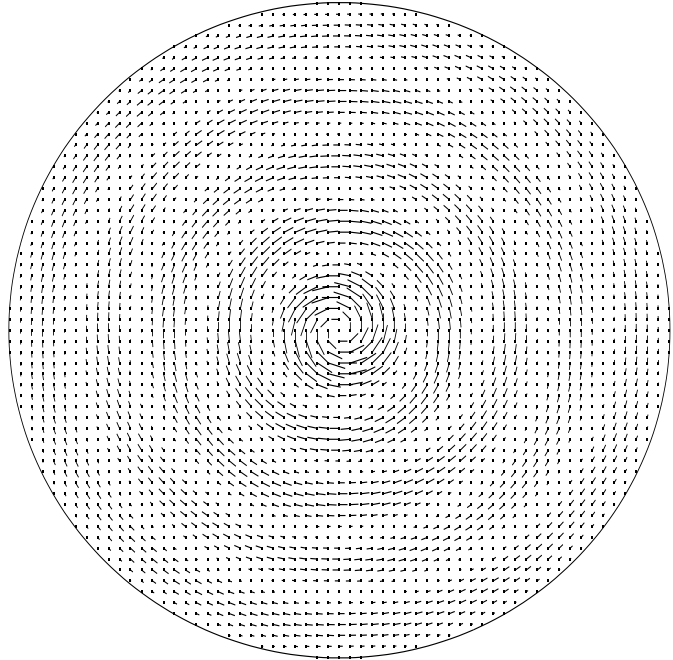}} 
 \caption { D2Q9 scheme for linear Navier--Stokes.  
 Eigenmode $n=5$ $\ell=1$  for the Stokes problem in a circle. } 
\label{f-d2q9-disque}
\end{figure}   
 
\bigskip \noindent $\bullet$ \quad 
As in the particular case of D2Q5 scheme, we have verified with periodic boundary
conditions that the relaxation rate of a transverse wave is determined  with 
error of order six  and relative fourth-order precision, as shown in Figure \ref{f-d2q9-unpt}. 
The  detailed numerical convergence plot 
is very similar to Figure  \ref{f-arno-therm}. 

\begin{figure}[!t]         
\centerline {\includegraphics [width=.76  \textwidth] {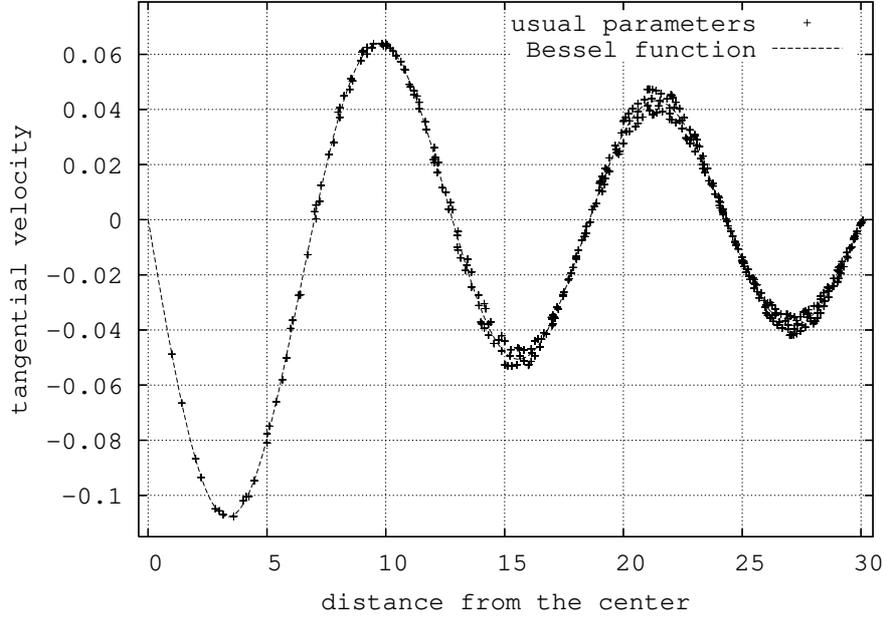}} 
 \caption { D2Q9 scheme for linear Navier--Stokes in a circle.  
 Eigenmode  $ n=5$, $\ell=1$ for the Stokes problem. Usual parameters.  } 
\label{f-d2q9-ns-mode-u}
\end{figure}

\begin{figure}[!h]         
\centerline {\includegraphics [width=.76  \textwidth] {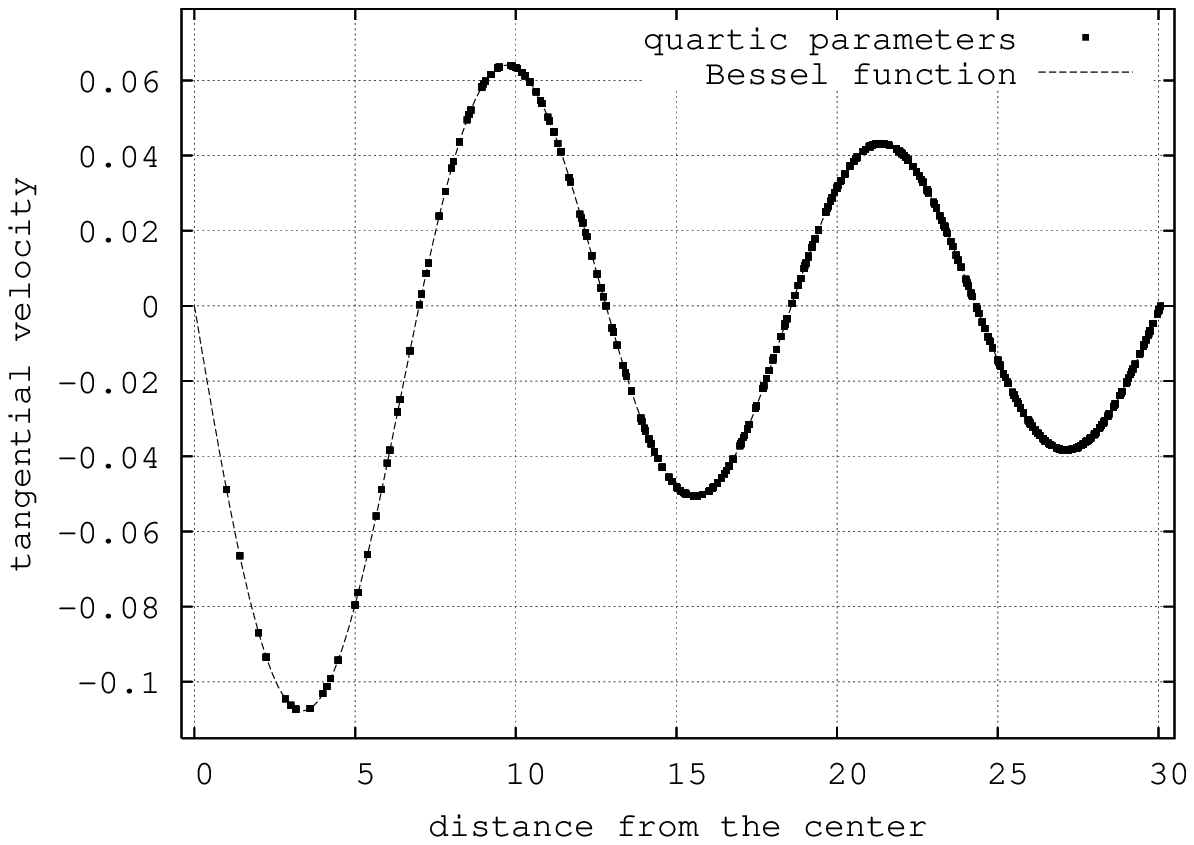}}
 \caption {  D2Q9 scheme for linear Navier--Stokes in a circle. 
 Eigenmode   $ n=5$, $\ell=1$ for the Stokes problem. Quartic parameters.   } 
\label{f-d2q9-ns-mode-q}
\end{figure}

\begin{figure}[!h]          
\centerline {\includegraphics [width=.76 \textwidth] {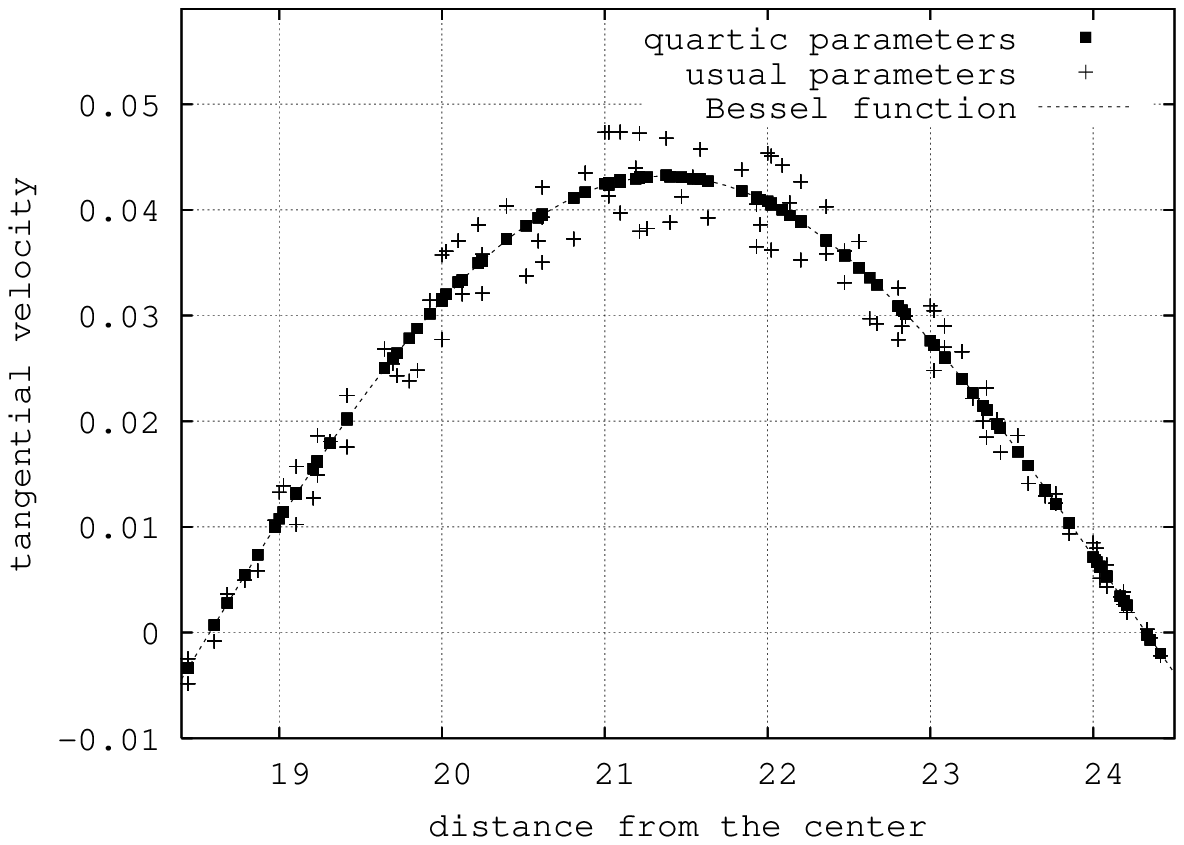}}   
 \caption {   D2Q9 scheme for linear Navier--Stokes in a circle. 
 Eigenmode   $ n=5$, $\ell=1$ for the Stokes problem.  Zoom of the  figures  
\ref{f-d2q9-ns-mode-u} and \ref{f-d2q9-ns-mode-q}.     } 
\label{f-d2q9-ns-mode-zoom}
\end{figure}

\bigskip \noindent $\bullet$ \quad 
We have also validated our results for eigenmodes of the Stokes problem inside  a circle. 
With the notations introduced previously, the eigenvalues $ \, \Gamma \,$ are given 
\cite{LL59}  by 
%
\begin{equation} \label{4.41}
\Gamma \,=\,  \nu \, \Big(  {{ \zeta_\ell^n } \over{ R }}  \Big)^2 \, . 
\end{equation}
The result for $ R = 30.07 $, $\ell=1$ and $ n=5 $ is presented 
in Figure \ref{f-d2q9-disque} 
for the velocity field with a mesh included in a square of size $ \, 61 \times 61 .$ 
 The alternance of directions for the vector field is clearly visible 
  on  Figure \ref{f-d2q9-disque} and we use around seven meshpoints between 
two zeros  of the Bessel function.  
We have compared with the same mesh the results obtained with the lattice Boltzmann scheme 
with the usual parameters that does not satisfy relation  (\ref{4.37}) 
but such that 
%
$ \, \nu  =   {{ \lambda^2 \,  \Delta t}\over{10}}  \, $ 
%
which is very close to   (\ref{4.38})   and quartic  parameters.
The radial profile of the tangential velocity is shown in Figures \ref{f-d2q9-ns-mode-u} 
to  \ref{f-d2q9-ns-mode-zoom}. 
The difference is visually spectacular. 
As for the thermics case, we observe that simple boundary conditions (here we
use those of Bouzidi {\it et al.} \cite{bfl})  prevent fourth-order convergence for the Stokes
problem. Use of more sophisticated boundary conditions  (see 
Ginzburg and d'Humi\`eres  \cite{GD95})   may
help to improve the convergence;  however for models with limited number of
velocities, it is not clear whether the choice of $s$-parameters will be the
same for ``fourth-order volume'' and ``accurate Poiseuille type boundary conditions''.

\bigskip \noindent $\bullet$ \quad 
{\bf The  D3Q19 for linearized athermal Navier--Stokes at   fourth-order }

\noindent
The D3Q19 model 
 is analyzed as was done above for the D2Q9 model. We detail in Appendix~C  
the way to enforce the precision of eigenmodes for the Stokes problem. We obtain a set of
eight equations for the coefficients $ \, \sigma$'s. 
These equations have  only 
one  nontrivial family of solutions  given according to  
\begin{equation}  \label{d3q19-quartic} \left\{   \begin{array}{rcrclcrcl}
\textrm {energy}  && \sigma_4 &=&  \frac{1}{s_{4}} -  \frac{1}{2}   && s_4
&=& \textit {ad libitum}      \\ [1mm]
\textrm {stress tensor} &&  \sigma_5 &=& \displaystyle  1 / \sqrt{12}   && s_5 &=&  3-\sqrt{3}      \\ 
\textrm {energy flux} &&  \sigma_{10} &=& 1 / \sqrt{3}  && s_{10} &=&  4 \sqrt{3} - 6  \\
\textrm {square of energy} &&  \sigma_{13} &=& \frac{1}{s_{13}} -  \frac{1}{2}  && 
s_{13} &=&  \textit {ad libitum}   \\ [1mm]
\textrm {other moments of kinetic energy} &&  \sigma_{14} &=&  \displaystyle  1 / \sqrt{12} 
  && s_{14} &=&  3 - \sqrt{3} \\
\textrm {third-order antisymmetric} &&  \sigma_{16} &=&   
1 / \sqrt{3}    && s_{16} &=&  4 \sqrt{3}-6    \, . 
\end{array} \right.    \end{equation}
Note these results are incompatible with BGK hypothesis (all $\sigma$ equal)
but are compatible with the ``two relaxation times'' hypothesis  which enforces equality of even moments
$ \, \sigma_4 = \sigma_5 = \sigma_{13}  =  \sigma_{14}  $,    
and of odd moments:  $ \,  \sigma_{10}  =  \sigma_{16} .  \,$   
We remark 
that the relaxation rate for energy (linked to the bulk viscosity) is not
constrained. Note   that     the shear viscosity   $\nu$   takes the value $1/\sqrt{108}$
as in  (\ref{4.38}).  
 As for D2Q9 there is no decoupling at order 3 of shear and
acoustic modes, and thus, at least at the present stage,  we make no claim
concerning possible improvements for the acoustic modes. We will study this question in
a forthcoming contribution.

\begin{figure}[!h]           
\centerline {\includegraphics [width=.80 \textwidth ] {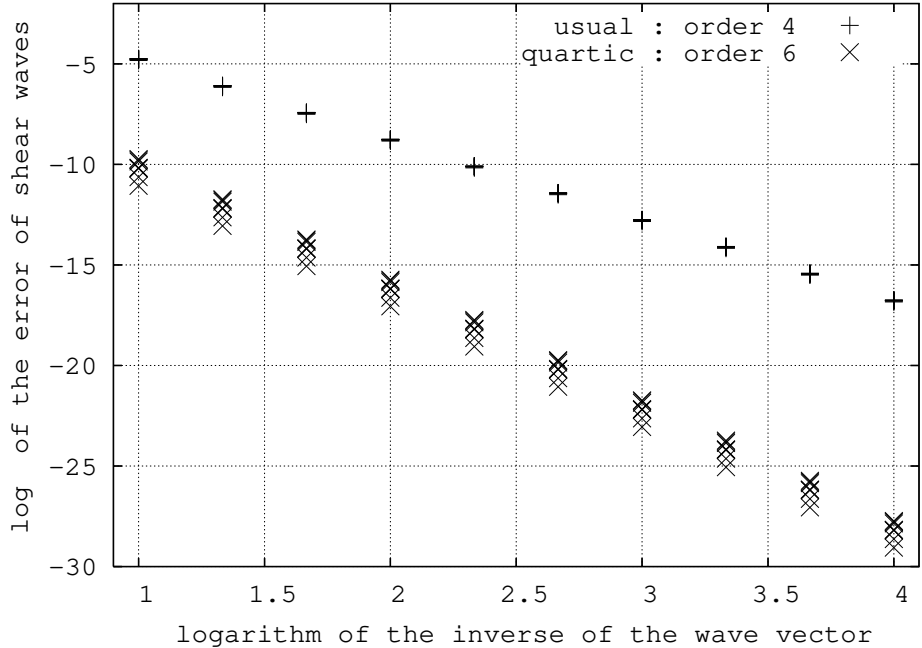}}   
 \caption {  D3Q19 for ``one point'' experiment and various directions of the wave
 vector.  }  
\label{f-d3q19-un-point}
\end{figure}   

\bigskip \noindent $\bullet$ \quad 
We have performed the same kind of numerical analysis as for the 
two-dimensional   D2Q9 case
and find quite similar results. We  illustrate our results first with a ``one point
experiment''. We introduce numerical wave vectors $ \, k \, $ close to zero 
and compute  the eigenmodes numerically. The shear mode is close to 
$\,\,   {{ \lambda^2}\over{3}}  \, \sigma_5   \, \mid k \mid ^2 \,$   
and we plot in Figure~\ref{f-d3q19-un-point} the experimental error. 
With ordinary coefficients, the error is of order 4,  whereas with the so-called 
``quartic coefficients'', the error is of order 6 and the relative error of order 4.

\begin{figure}[!t]            
\centerline {\includegraphics [width=.51 \textwidth, height=.51  \textwidth]  
{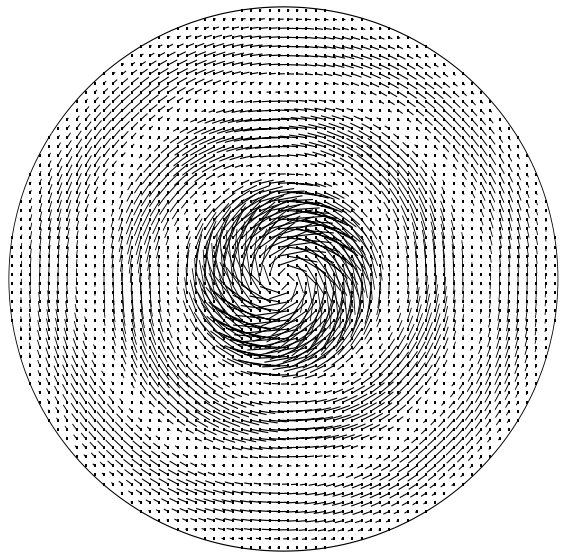}} 
 \caption {  D3Q19 for linear Navier--Stokes in a sphere.
 Eigenmode   $n=3$,   $\ell=1$   for  Stokes problem 
  with  Dirichlet boundary conditions.   Tangential velocity vector field for a plane 
through the center of the sphere. }
\label{f-d3q19-ns-sphere-plane}
\end{figure}    
 
\begin{figure}[!h]           
\centerline {\includegraphics [width=.51 \textwidth, height=.51  \textwidth]  
{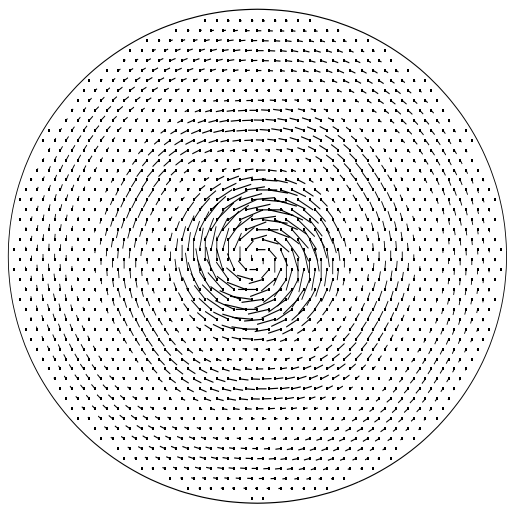}}  
 \caption { D3Q19 for linear Navier--Stokes in a sphere.   
 Eigenmode   $n=3$,   $\ell=1$   for  Stokes problem   
 with  Dirichlet boundary conditions. 
 Tangent vector field for a plane orthogonal to vector $(1,\, 1,\, 1)$.    }  
\label{f-d3q19-ns-sphere-ortho}
\end{figure}        

\begin{figure}[!h]            
\centerline {\includegraphics [width=.76 \textwidth]  {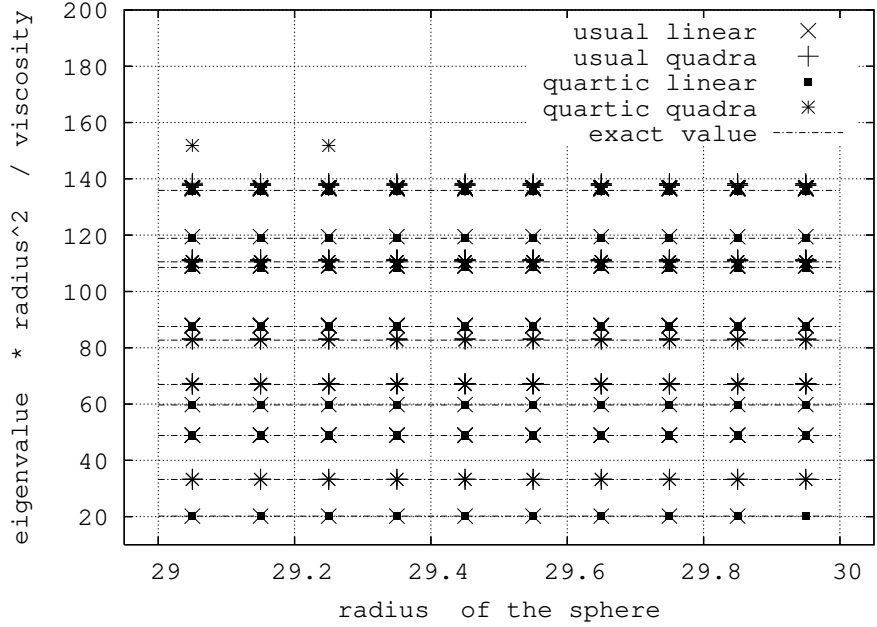}} 
\caption {  D3Q19 for linear Navier--Stokes in a sphere. First eigenmodes 
 for stationary Stokes problem with  Dirichlet boundary conditions.    }  
\label{f-d3q19-ns-sphere-modes}
\end{figure}         
 
\begin{figure}[!h]           
\centerline {\includegraphics [width=.76 \textwidth]  {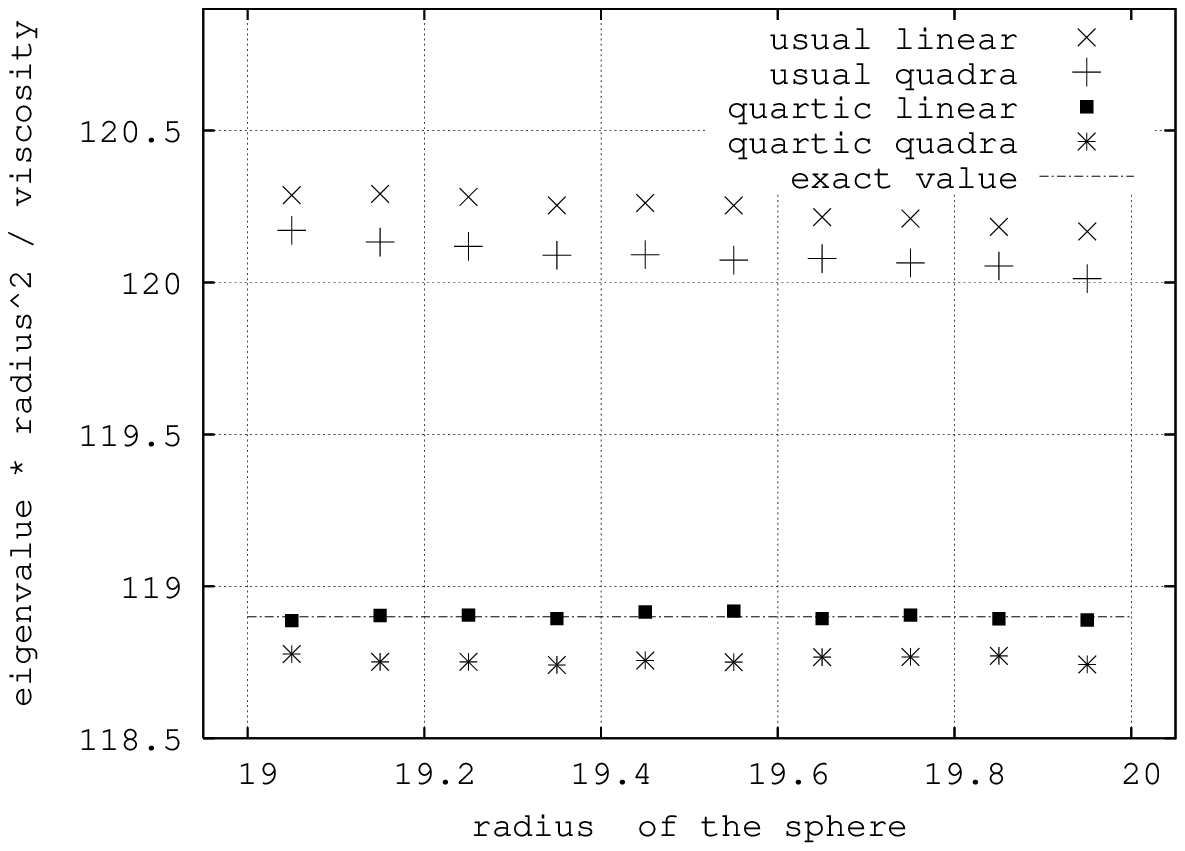}}  
\caption {  D3Q19 for linear Navier--Stokes in a sphere.
  Eigenmode for stationary Stokes problem. 
 Zoom of various schemes for Dirichlet eigenmode close to  118.8998692.     }  
\label{f-d3q19-ns-sphere-119}
\end{figure}     

\bigskip \noindent $\bullet$ \quad 
We also illustrate our results for the problem 
of Stokes modes in a sphere which has an analytical solution in terms of
Bessel functions.
The Stokes problem   searches for  a velocity
field $ {\bf u} (r,\, t)$   with $ {\bf u} =0$ on the surface of a sphere of radius $R$.
An analysis  similar to that for the Stokes problem in a circle, leads to
an eigenvalue problem, with solutions
\begin{equation*} \label{4.43}      
\Gamma \,=\,  \nu  \,   \bigg(  {{  \zeta_{\ell+1/2}^n } \over{ R }}  \bigg)^2 \,, 
\qquad \ell \geq 1 \,,  
\end{equation*}
analogous to (\ref {4.26}), 
with $ \, \smash {  \zeta_{\ell+1/2}^n } \,$ equal to the $n^{\rm th}$ zero of the 
``semi-integer'' Bessel function  $\,   J_{\ell+1/2} \,$ as defined in 
 Abramowitz and Stegun  \cite{AS65}.  
Using the Arnoldi technique, we can determine a few eigenvalues and verify that
they are close to the theoretical formula.    
We present Figures  \ref  {f-d3q19-ns-sphere-plane} and  \ref {f-d3q19-ns-sphere-ortho}  
an example of a typical result   obtained with this framework. 
We find that these eigenvalues have the expected
degeneracy $2 \ell + 1$. Note however  that,  the computations being made for a rather
small radius $R$, there are small splittings of the degenerate eigenvalues due to the
fact that lattice Boltzmann computations have cubic symmetry.

\bigskip \noindent $\bullet$ \quad 
For a more detailed analysis, we take advantage of the symmetry of the Stokes problem
and therefore perform computations on an eighth   of the sphere,  taking proper account of the
symmetry with respect to the planes perpendicular to the coordinates $x,\,\,y,\,z\,,$
through the center of the cube (symmetry or anti-symmetry).
Using four  different combinations of symmetries on the planes,  we can determine all
the eigenvalues, the other combinations leading to the same eigenvalues with
only a permutation in the coordinates for the eigenmodes. Note that
due to the rather high complexity of the Arnoldi prodecure, this
allows a  reduction in computer time of two orders of magnitude. 

\bigskip \noindent $\bullet$ \quad 
We present in Figure \ref {f-d3q19-ns-sphere-modes} 
the effect of boundary conditions for a number of values of the
radius from 29 to 30. We give in Figure  \ref {f-d3q19-ns-sphere-119}  
some details for $R$ between 19 and 20 for the $m=1,\, n=6\,$ mode. 
There are two sets of data, one for usual  $s$-parameters
%
\begin{equation*} \label{4.44}  
s_4 =    1.3 \,,\quad    
s_5  = 1.25   \,,\quad       
s_{10}  = 1.2  \,,\quad    
s_{13}  = 1.4  \,,\quad   
s_{14} =  1.25  \,,\quad  
s_{16} =  1.3     
 \end{equation*}
and one for the quartic  $s$-parameters given precedingly in 
(\ref{d3q19-quartic}) with 
\begin{equation*}  
s_4=1.3   
\,, \quad  s_{13}=1.4 \,.   
 \end{equation*}
Similar work has been done for a cube. The results are published in  Leriche
{\it et al}  in \cite{LLL08}.


\section{Conclusion }

\smallskip  \noindent $\bullet$ \quad 
The expansion of equivalent equations that are satisfied by the mean quantities
determined by the lattice Boltzmann method has been described in this 
contribution  and  
explicit formulae given for a few models up to order four in space derivatives.
Extending either to more complicated models or to higher order derivatives 
is very simple and 
does not imply new conceptual developments, in particular careful treatment
of non commuting terms that appear in the Chapman-Enskog procedure.
The developments imply only simple
algebraic manipulations that can be performed by a ``formal language'' program,
as used here. Note that these developments have a rather high complexity as seen
by the fact that each order takes roughly ten  times as much computer time as the
preceding one.
 
\bigskip \noindent $\bullet$ \quad 
With the Taylor expansion method, we can obtain explicit formulae which, then, unable us
to   tune some parameters of the lattice
Boltzmann scheme initially proposed by d'Humi\`eres in order to capture, up to 
fourth-order  accuracy,  shear waves. 
Of course, this extra-precision obtained with a classical scheme is
possible only if the viscosity is essentially fixed and the expansion done around zero
velocity.  
Even though very few situations were studied here, it can be said that tuning
the accuracy of the ``internal code'' independently from the method  to take
care of boundary conditions allows us to get useful information concerning these
two sources of errors in lattice Boltzmann simulations. Future extension 
of this work will be to try and discriminate between some of the numerous proposed
ways to deal with boundaries to be able to estimate their contributions to
errors in comparison to those due to the ``internal code''.


\section*{Appendix A.   Taylor expansion method}

\smallskip \noindent $\bullet$ \quad 
We start from relation  (\ref{1.21})  
for iteration of the lattice Boltzmann scheme
and take the momentum of order $k$. Then 
\begin{eqnarray*}   
m_k(x,\, t+\Delta t) \,& = & \, \sum_{\ell=0}^J \, M_{k \ell} 
\, f_{\ell}^* (x - v_{\ell} \,  \Delta t, \, t)    \\ 
 \,& = & \, \sum_{\ell=0}^J \,  \sum_{p=0}^J  \,    M_{k \ell} \, M^{-1}_{\ell p} 
\, m_{p}^* (x - v_{\ell} \,  \Delta t, \, t) \qquad {  \textrm {  due to (\ref{1.19}) } } \\ 
 \, & = & \, \sum_{\ell=0}^J \,  \sum_{p=0}^J  \,  \sum_{r=0}^J      M_{k \ell} \, M^{-1}_{\ell p} 
\, \Psi_{p r} \,  m_{r}  (x - v_{\ell} \,  \Delta t, \, t) 
\end{eqnarray*}
due to  (\ref{1.27}). We have 
\begin{equation} \label{6.1} 
m_k(x,\, t+\Delta t) \,=\,  \sum_{\ell = 0}^J  \sum_{p = 0}^J   \sum_{r = 0}^J
M_{k \ell} \, M^{-1}_{\ell p} 
\, \Psi_{p r} \,  m_{r}  (x - v_{\ell} \,  \Delta t, \, t) \,, \qquad 0 \leq k \leq J \,.
\end{equation} 
We expand now momentum $ \,  m_{r}  (x - v_{\ell} \,  \Delta t, \, t) \,$ with 
a Taylor formula of infinite length:
\begin{equation} \label{6.2}   
 m_{r}  (x - v_{\ell} \,  \Delta t, \, t) \,=\, \sum_{q=0}^{+\infty} 
{{(\Delta t)^q}\over{ q \, ! }}  \Big(- \sum_{\alpha=1}^d M_{\alpha \ell} \, \partial_{\alpha}
 \Big)^q \, m_r(x,\,t) \, . 
\end{equation} 
Then due to (\ref{6.1}),  (\ref{6.2}) and  (\ref{2.6}), we have 
%
\begin{equation} \label{6.7}   
m_k(x,\, t+\Delta t) \,=\,   \sum_{    \gamma   }  
 \sum_{\ell = 0}^J  \sum_{p = 0}^J   \sum_{r = 0}^J
\, M_{k \ell} \, M^{-1}_{\ell p} 
\, \Psi_{p r} \,   {{\Delta t ^{\mid \gamma \mid}}\over{\mid \! \gamma \!  \mid ! }} \, 
P_{\ell \gamma} \, \partial_{\gamma} m_r  \,, \qquad 0 \leq k \leq J \, .
\end{equation} 
We can also expand the left hand side of  (\ref{6.7}) and we have finally 
\begin{equation} \label{6.8}    
 \sum_{q=0}^{\infty}  
 {{\Delta t ^{q}}\over{ q ! }} \,  \partial_{t}^q  m_k  \,=\,  
 \sum_{    \gamma   }   \sum_{\ell = 0}^J  \sum_{p = 0}^J   \sum_{r = 0}^J \, 
 M_{k \ell} \, M^{-1}_{\ell p} 
\, \Psi_{p r} \,   {{\Delta t ^{\mid \gamma \mid}}\over{\mid \! \gamma \!  \mid ! }} \, 
P_{\ell \gamma} \, \partial_{\gamma} m_r  \,, \quad 0 \leq k \leq J \, .
\end{equation} 

\bigskip  \noindent $\bullet$ \quad 
We consider relation (\ref{6.8}) at order zero relative to time step for a conserved component
of momentum ({\it id est} $ \, 0 \leq k \equiv i \leq N-1$). The left hand side of (\ref{6.8}) is
equal to $ \, m_i + {\rm O}(\Delta t) \, $ and we have   
\begin{eqnarray*}   
 W_i + {\rm O}(\Delta t) \,& = & \,   \sum_{\ell = 0}^J  \sum_{p = 0}^J   \sum_{r = 0}^J
         \,   M_{i \ell} \,  \, M^{-1}_{\ell p} 
\, \Psi_{p r} \,  m_r  \,  + {\rm O}(\Delta t)  \\ 
 \,& = & \,  \sum_{p = 0}^J   \sum_{r = 0}^J \, 
\delta_{i p} \, \Psi_{p r} \,  m_r  \,  + {\rm O}(\Delta t)  \\  
  \, & = & \,   \sum_{r = 0}^J \, 
 \Psi_{i r} \,  m_r  \,  + {\rm O}(\Delta t)  \qquad   \textrm  {with} 
 \,\, 0 \leq i \leq N     \\  
  \, & = & \,   \sum_{r = 0}^J \,  \delta_{i r} \,  m_r  \,  + {\rm O}(\Delta t)  \qquad   
 {  \textrm {  due to (\ref{1.28}) } }   \\  
  \, & = & \,  m_i  + {\rm O}(\Delta t)
\end{eqnarray*}  
and no information is contained at this first step. Consider now the same development for
$ \, k \geq N$. We pass over some  repeated summations: 
\begin{eqnarray*}   
m_k  + {\rm O}(\Delta t) \, & = & \, M_{k \ell} \,  \, M^{-1}_{\ell p} 
\, \Psi_{p r} \,  m_r  \,  + {\rm O}(\Delta t)  \,  \\ 
 \, & = &   \, \sum_{j=0}^{N-1} \,   M_{k \ell} \,  \, M^{-1}_{\ell p} 
\, \Psi_{p j} \,  m_j  \, + \,   \sum_{r \geq N}  \,M_{k \ell}   \, M^{-1}_{\ell p} 
\, \Psi_{p r} \,  m_r  \,   + {\rm O}(\Delta t)  \,  \\ 
 \,& = &  \, \sum_{j=0}^{N-1} \, \delta_{k p}  
\, \Psi_{p j} \,  W_j  \, + \,   \sum_{r \geq N}  \, M_{k \ell}   \, M^{-1}_{\ell p} 
\, \delta_{p r} \,(1 - s_r) \,   m_r  \,   + {\rm O}(\Delta t)  \\ && 
~ \qquad  \qquad  \qquad  \qquad  \qquad  \qquad  \qquad  \qquad  \qquad  \qquad 
  {  \textrm {  due to (\ref{1.27}) and  (\ref{1.28})} }  \\ 
 \,& = &  \, \sum_{j=0}^{N-1}  \, 
\, \Psi_{k j} \,  W_j  \, + \,      M_{k \ell} \, M^{-1}_{\ell p} 
 \,(1 - s_p) \,   m_p  \,   + {\rm O}(\Delta t)  \\ 
 \,& = &  \,  \delta_{k p} \,  (1 - s_p) \,   m_p  \, + \, 
 \sum_{j=0}^{N-1}  \, \Psi_{k j} \,  W_j  \,   + {\rm O}(\Delta t) \\ 
 \,& = &  \,   (1 - s_k) \,   m_k    \, + \, 
 \sum_{j=0}^{N-1}   \, \Psi_{k j} \,  W_j  
 \,   + {\rm O}(\Delta t)  \,. 
\end{eqnarray*}  
%
We deduce from the previous calculation the relation (\ref{2.11}) with the 
expression  (\ref{2.10}) of the coefficients $ \,B^0_{k j} \, .$   Then   
we can go now one step further. 

\bigskip  \noindent $\bullet$ \quad 
At first order, relation (\ref{6.8}) becomes 
\begin{equation} \label{6.12}   
m_k \,+\, \Delta t \, {{\partial m_k}\over{\partial t}} \,+ \, {\rm O}(\Delta t^2)  \,=\,  
m_k^* \, - \,  \Delta t \, \sum_{\alpha=1}^d  \,
M_{k \ell} \, M^{-1}_{\ell p}  \, \Psi_{p r} \,   
M_{\alpha \ell} \,  \partial_{\alpha}     m_r  
\,+ \, {\rm O}(\Delta t^2)\,. \,
\end{equation}  
%
For conserved variables (\ref{1.9}) ({\it id est} $ \, 0 \leq k\equiv i \leq N-1$), we 
have   after dividing by $\, \Delta t $, 
\begin{eqnarray*}   
{{\partial W_i}\over{\partial t}} \, + \,  {\rm O}(\Delta t)  \, & = & \, 
- \sum_{\alpha=1}^d \,    M_{i \ell} \, M^{-1}_{\ell p}  \, \Psi_{p r} \,  
 M_{\alpha \ell} \,  \partial_{\alpha} m_r   \,+ \, {\rm O}(\Delta t)  \\  
 \, & = & \,  - \sum_{\alpha=1}^d \, \Lambda_{\alpha i}^p \,  \Psi_{p r} \,   
 \partial_{\alpha} m_r   \,+ \, {\rm O}(\Delta t)   \,  
\qquad  {  \textrm {  due to (\ref{1.25}) } }   \\  
 \, & = & \,    \sum_{\alpha=1}^d \, \Lambda_{\alpha i}^p \,  \Big(  
  \sum_{j < N }  \Psi_{p j} \,  
 \partial_{\alpha} W_j   \,+ \,   \sum_{\ell \geq N}  \Psi_{p \ell} \,   
 \partial_{\alpha} m_{\ell} \Big)   \,+ \, {\rm O}(\Delta t)\,  \\  
 \, & = & \,   \sum_{\alpha=1}^d \, \Lambda_{\alpha i}^p \, \sum_{j < N }  \Big(     \Psi_{p j} \,  
 \partial_{\alpha} W_j   \,+ \,   \sum_{\ell \geq N}  \Psi_{p \ell} \,   
 \partial_{\alpha}  \big(  {{1}\over{s_{\ell}}} \,   \Psi_{\ell j} \,  W_j  \big) 
 \Big)   \,+ \, {\rm O}(\Delta t) \,   \\  
 \, & = & \,  \sum_{j = 0 }^{N-1} \,  \,
  \sum_{\alpha=1}^d \,  \sum_{p = 0 }^{J} \, \Lambda_{\alpha i}^p \,  \Big(  \Psi_{p j} \,  + \, 
 \sum_{\ell \geq N}  \Psi_{p \ell} \,    {{1}\over{s_{\ell}}} \,   \Psi_{\ell j} 
 \Big)  \,  \partial_{\alpha} W_j  \,+ \, {\rm O}(\Delta t) \,. 
\end{eqnarray*}    
For an index $ \, \gamma \,$ between $1$ and $d$, we define $ \, A^{\gamma}_{ij} \, $
according to the relation  (\ref{2.13})    
and the previous calculation can be written as a conservation law at first-order
\begin{equation} \label{6.14}   
{{\partial W_i}\over{\partial t}} \, + \, \sum_{\mid \gamma \mid = 1}
A^{\gamma}_{ij} \,  \partial_{\gamma} W_j \,=\,  \, {\rm O}(\Delta t)   \,, 
\quad 0 \leq i \leq N-1  \, . \,   
\end{equation}    

\bigskip  \noindent $\bullet$ \quad 
We start again from relation  (\ref{6.12}) with nonconservative indices $k$ ($k \geq N$): 
\begin{equation*}    
m_k \,=\, - \Delta t   \, {{\partial m_k}\over{\partial t}} \,+ \,  
 (1 - s_k) \,   m_k    \, + \,   \sum_{j=0}^{N-1}   \, \Psi_{k j} \,  W_j   \,
\, - \,  \Delta t  \, \sum_{\alpha=1}^d \, 
M_{k \ell} \, M^{-1}_{\ell p}  \, \Psi_{p r} \,
M_{\alpha \ell} \,  \partial_{\alpha}     m_r   \,+ \, {\rm O}(\Delta t^2)\,.  
\end{equation*}    
Then due to  (\ref{2.11}), 
%
\begin{eqnarray*}   
m_k \, & = & \,   {{1}\over{s_{k}}} \, \bigg(   \Psi_{k j} \,  W_j   
\, - \Delta t  \,  {{1}\over{s_{k}}}  \Psi_{k i} \,   {{\partial W_i}\over{\partial t}}  
\, - \Delta t  \, \Lambda_{\alpha k}^p \,  \Psi_{p r} \,  \partial_{\alpha}  
\Big(  {{1}\over{s_{r}}}    \Psi_{r j} \,  W_j   \Big)   \bigg) 
 \,+ \, {\rm O}(\Delta t^2)  \\ 
 \, & = & \,   {{1}\over{s_{k}}} \, \bigg(   \Psi_{k j} \,  W_j   
\, + {{\Delta t}\over{s_{k}}} \,   \Psi_{k i}  \, \sum_{\mid \gamma \mid = 1}
A^{\gamma}_{ij} \,  \partial_{\gamma} W_j   
\, -   {{\Delta t}\over{s_{r}}}  \, \Lambda_{\gamma k}^p \,  \Psi_{p r} \,  
    \Psi_{r j} \,   \partial_{\gamma}  W_j   \bigg)    
 \,+ \, {\rm O}(\Delta t^2)  \,.\,   
\end{eqnarray*}    
%
We introduce $ \, B^{\gamma}_{kj} \, $ 
for $\, \mid \! \gamma \! \mid = 1 \,$ according to  (\ref{2.15}) 
and due to previous calculation,  relation   (\ref{2.11}) can be extended as 
\begin{equation} \label{6.16}    
m_k  \,=\,   \sum_{0 \leq \mid \gamma \mid \leq 1}  \Delta t ^{\mid \gamma \mid} \, \, 
 B^{\gamma}_{kj} \, \,  \partial_{\gamma} W_j   \,   + \, {\rm O}(\Delta t^2)  \,.  
\end{equation}    

\bigskip  \noindent $\bullet$ \quad
We generalize the relations   (\ref{6.14})  and  (\ref{6.16})  at the order $ \, \sigma \,$ 
 through a recurrence hypothesis   (\ref{2.17})   (\ref{2.18}).  
In order to treat the left-hand side of relation  (\ref{6.8}), we observe that we have 
\begin{eqnarray*}    
\partial_t^2 W_i  \, & = & \, -  \sum_{1 \leq \mid \gamma \mid \leq \sigma}
\,  \Delta t ^{\mid \gamma \mid - 1} \, 
A^{\gamma}_{ij} \,  \partial_{\gamma} \, \Big( \partial_t   W_j \Big)  
\, + \,   {\rm O}(\Delta t^{\sigma})  \\ 
 \, & = &\,  \sum_{1 \leq \mid \delta \mid \leq \sigma}
\,  \Delta t ^{\mid \delta \mid - 1} \, 
A^{\delta}_{i \ell} \,  \partial_{\delta} \, \Big( 
\sum_{1 \leq \mid \varepsilon \mid \leq \sigma}
\,  \Delta t ^{\mid \varepsilon \mid - 1} \, 
A^{\varepsilon}_{\ell j} \,  \partial_{\varepsilon} W_j  \Big)  
\, + \,   {\rm O}(\Delta t^{\sigma}) 
\end{eqnarray*}  
and if we introduce $ \, C^{1,\gamma}_{ij} \,$ according to   (\ref{2.56}) and  
\begin{equation*} \label{6.19}    
C^{2,\gamma}_{ij} \,\,\, \equiv \,  - \sum_{ \mid \delta \mid \geq 1 ,\, 
 \mid  \varepsilon \mid \geq 1 ,\, \delta + \varepsilon = \gamma}
A^{\delta}_{i \ell}  \, \,  A^{\varepsilon}_{\ell j}  \, \,\,, 
\quad 2 \leq \mid \! \gamma  \! \mid  \leq \sigma +1 \,,  \, 
\end{equation*}  
we have for the second time derivative a relation quite analogous to (\ref{2.17}):
\begin{equation*}    
\partial_t^2 W_i  \,\,\, +   \!\!   \sum_{2 \leq \mid \gamma \mid \leq \sigma+1} \, 
 \Delta t ^{\mid \gamma \mid - 2} \, \, 
C^{2, \gamma}_{ij} \, \,  \partial_{\gamma} W_j \,=\,   {\rm O}(\Delta t^{\sigma}) \,, 
\quad 0 \leq i \leq N-1  \,. \,  
\end{equation*}  
This relation can be generalized at an arbitrary order according to 
\begin{equation} \label{6.20} 
\partial_t^q W_i   +   \!\!\!\!   \sum_{q \leq \mid \gamma \mid \leq \sigma+q-1} \, 
 \!\!\!\!  \Delta t ^{\mid \gamma \mid - q} \, \, 
C^{q, \gamma}_{ij} \, \,  \partial_{\gamma} W_j \,=\,   {\rm O}(\Delta t^{\sigma}) \,, 
\quad 0 \leq i \leq N-1  \,. \,  
\end{equation}  
If relation  (\ref{6.20})  is true at order $q$, we have by differentiation with respect  to time, 
\begin{eqnarray*}    
\partial_t^{q+1} W_i  \, & = & \, -   \!\!\!\!   \sum_{q \leq \mid \gamma \mid \leq \sigma+q-1} \, 
 \!\!\!\!  \Delta t ^{\mid \gamma \mid - q} \, \, 
C^{q, \gamma}_{ij} \, \,  \partial_{\gamma}  \, \Big( \partial_t   W_j \Big) 
 \,+\,   {\rm O}(\Delta t^{\sigma})   \\ 
 \, & = & \,      \sum_{q \leq \mid \delta \mid \leq \sigma+q-1} \, 
 \!\!\!\!  \Delta t ^{\mid \delta \mid - q} \, \, 
C^{q, \delta}_{i \ell} \, \,  \partial_{\delta}  \, \Big( 
\sum_{1 \leq \mid \varepsilon \mid \leq \sigma}
\,  \Delta t ^{\mid \varepsilon \mid - 1} \, 
A^{\varepsilon}_{\ell j} \,  \partial_{\varepsilon} W_j  \Big)  
\, + \,   {\rm O}(\Delta t^{\sigma})   \\ 
 \, & \equiv  & \,   \sum_{q+1 \leq \mid \gamma \mid \leq \sigma+q} \, 
 \!\!\!\!  \Delta t ^{\mid \gamma \mid - q-1} \, \, 
C^{q+1, \gamma}_{ij} \, \,  \partial_{\gamma}  \, \Big( \partial_t   W_j \Big) 
 \,+\,   {\rm O}(\Delta t^{\sigma}) 
\end{eqnarray*}  
and relation  (\ref{6.20}) is satisfied at the order  $\, q+1 \,$ 
with $ \, C^{q+1, \gamma}_{ij} \, $ given by the recurrence relation   (\ref{2.21}).    
In an analogous way, we have 
\begin{equation} \label{6.22} 
\partial_t^q  m_k \,=\,   \sum_{q \leq \mid \gamma \mid \leq  \, \sigma+q} 
 \Delta t ^{\mid \gamma \mid - q } \, \, 
 D^{q,\gamma}_{kj} \, \,  \partial_{\gamma} W_j   \, + \, {\rm O}(\Delta t^{\sigma + 1})
 \,, \qquad  k \geq N \, ,  
\end{equation} 
with  $ \, D^{0,\gamma}_{kj} \, $ defined according to (\ref{2.57}). 
If the relation  (\ref{6.22})  is satisfied at order $q$, we have by differentiation relative to time, 
\begin{eqnarray*}     
\partial_t^{q+1} m_k \, & = & \,   \sum_{q \leq \mid \gamma \mid \leq  \, \sigma+q} 
 \Delta t ^{\mid \gamma \mid - q } \, \, 
 D^{q,\gamma}_{kj} \, \,  \partial_{\gamma}  \Big( \partial_t   W_j \Big) 
  \, + \, {\rm O}(\Delta t^{\sigma + 1})  \\ 
\, & = & \,  -  \sum_{q \leq \mid \delta \mid \leq  \, \sigma+q} 
 \Delta t ^{\mid \delta \mid - q } \, \, 
 D^{q,\delta}_{k \ell} \, \,  \partial_{\delta}  \Big( 
\sum_{1 \leq \mid \varepsilon \mid \leq \sigma}
\,  \Delta t ^{\mid \varepsilon \mid - 1} \, 
A^{\varepsilon}_{\ell j} \,  \partial_{\varepsilon} W_j  \Big)  
\, + \, {\rm O}(\Delta t^{\sigma + 1})   \\ 
\, & \equiv & \,  \sum_{q+1 \leq \mid \gamma \mid \leq  \, \sigma+q+1} 
 \Delta t ^{\mid \gamma \mid -( q + 1) } \, \, 
 D^{q+1,\gamma}_{kj} \, \,  \partial_{\gamma} W_j  
 \, + \, {\rm O}(\Delta t^{\sigma + 1})
\end{eqnarray*}  
with coefficients $\,  D^{q+1,\gamma}_{kj} \,$ determined according to 
the relation  (\ref{2.55}).  We observe that for the particular value 
$ \, \mid \! \gamma \! \mid = \sigma +1 \,$  the coefficient 
 $\,  D^{q+1,\gamma}_{kj} \,$ 
is well defined for $ \, 0 \leq q   \leq \sigma \, $. In other words,  
the coefficient $\,  D^{q,\gamma}_{kj} \,$ is well  defined 
for $ \, 1 \leq q   \leq  \, \mid \! \gamma \! \mid \, $.

\bigskip  \noindent $\bullet$ \quad
We verify now by induction that the recurrence relations  (\ref{2.17}) 
and  (\ref{2.18}) are satisfied.
It is the case at  order 1 as we have shown in   (\ref{6.14})  and  (\ref{6.16}). 
We first consider a label $i$ such that $ 0 \leq i \leq N-1$. Then according to  (\ref{6.8}), we have
at the order $ \, \sigma + 2 \, $: 
\begin{equation*}  \left\{   \begin{array} {c} \displaystyle 
W_i  \,+\, \Delta t \,  {{\partial W_i}\over{\partial t}} \, + \, 
\sum_{q=2}^{\sigma+1} \,  {{\Delta t ^{q}}\over{ q ! }} \,  \partial_{t}^q W_i 
\,+\,  {\rm O}(\Delta t^{\sigma + 2}) \, \,   =  
\qquad  \qquad  \qquad  \\   [1mm]  \displaystyle \qquad
   W_i \,+\,    M_{i \ell} \, M^{-1}_{\ell p}  \, \Psi_{p r} 
  \sum_{1 \leq \mid \delta \mid \leq \sigma+1} \, 
 \frac{\Delta t ^{\mid \delta \mid}}{\mid  \delta   \mid ! } \, 
P_{\ell \delta} \, \, \partial_{\delta} \bigg( 
 \sum_{0 \leq  \mid \varepsilon \mid \leq \sigma} \, \Delta t^{ \mid \varepsilon  \mid} \, 
B^{ \varepsilon}_{rj} \,  \partial_{ \varepsilon} W_j  \bigg)
+  {\rm O}(\Delta t^{\sigma + 2}) \, .  
\end{array} \right.  \end{equation*}
We use relation  (\ref{6.20})  for  the left hand side of previous relation. We get after dividing
by $ \, \Delta t \,$  
\begin{equation*}  \left\{   \begin{array} {c} \displaystyle 
\displaystyle  {{\partial W_i}\over{\partial t}} \, - \, 
 \sum_{q=2}^{\sigma+1} \,  {{\Delta t ^{q-1}}\over{ q \, ! }} \,  
 \!\!\!\!   \sum_{q \leq \mid \gamma \mid \leq \sigma+q-1} \, 
 \!\!\!\!  \Delta t ^{\mid \gamma \mid - q} \, \, 
C^{q, \gamma}_{ij} \, \,  \partial_{\gamma} W_j   
\,+\,  {\rm O}(\Delta t^{\sigma + 1}) \, \, =   
\qquad  \qquad  \qquad  \\   [1mm]  \displaystyle \qquad   
 \sum_{1 \leq \mid \delta \mid \leq \sigma+1, \,  0 \leq \mid \varepsilon \mid \leq \sigma } 
  \,  M_{i \ell} \, M^{-1}_{\ell p} \, \Psi_{p r} \, P_{\ell \delta} \, 
  {{\Delta t ^{\mid \delta \mid +  \mid \varepsilon  \mid -1}}
\over{ \mid \! \delta \!  \mid \, !  }}  \, B^{ \varepsilon}_{rj} \, 
\partial_{\delta +  \varepsilon}  W_j  
\,  +\,  {\rm O}(\Delta t^{\sigma + 1}) \, .  
\end{array} \right.  \end{equation*} 
and the relation  (\ref{2.17}) is extended one step further with a coefficient 
$ \, A^{\gamma}_{ij} \, $ defined for $ \, \mid \! \gamma  \!  \mid \, = \, q+1 \,$ 
 by the recurrence relation   (\ref{2.24}).      
For the nonconserved moments ($k \geq N$), the relation  (\ref{6.8}) can be written at the order 
$\, \sigma + 2$ as 
\begin{equation*}  \left\{   \begin{array} {c} \displaystyle 
m_k \,+\,  \sum_{q=1}^{\sigma+1} \,  {{\Delta t ^{q}}\over{ q \, ! }} \,  
\partial_t^q m_k \,+\,   {\rm O}(\Delta t^{\sigma + 2}) \,\,=  
\qquad  \qquad  \qquad  \qquad  \qquad \qquad  \qquad  \qquad   \\   [1mm]  \displaystyle \,  
  (1-s_k)\, m_k  \,+ \!\!   \sum_{1 \leq \mid \delta \mid \leq \sigma + 1 } \, 
 M_{k \ell} \, M^{-1}_{\ell p} 
\, \Psi_{p r} \,   {{\Delta t ^{\mid \delta \mid}}\over{\mid \! \delta \!  \mid ! }} \, 
    \displaystyle 
P_{\ell \delta} \, \partial_{\delta}   \bigg( 
 \sum_{0 \leq \mid \varepsilon \mid \leq  \, \sigma}    
\Delta t ^{\mid \varepsilon \mid} \, \, 
 B^{\varepsilon}_{rj} \, \,  \partial_{\varepsilon} W_j \bigg)  
\, + \,  {\rm O}(\Delta t^{\sigma + 2})  \,.    
\end{array}  \right. \end{equation*}
We use the relation   (\ref{6.22})  and we deduce: 
\begin{equation*}  \left\{   \begin{array} {c} \displaystyle 
s_k \, m_k  \,\, =  \displaystyle  
\,\, -  \sum_{q=1}^{\sigma+1} \,  {{\Delta t ^{q}}\over{ q \, ! }} \,  
 \sum_{q \leq \mid \gamma \mid \leq  \, \sigma+q} 
 \Delta t ^{\mid \gamma \mid - q } \, \, 
 D^{q,\gamma}_{kj} \, \,  \partial_{\gamma} W_j   
\qquad  \qquad  \qquad  \qquad  \qquad  \qquad  \qquad  
  \\   [1mm]      \displaystyle  \qquad  \qquad  \qquad  
   \, + \sum_{1 \leq \mid \delta \mid \leq \sigma + 1,\,
0 \leq \mid \varepsilon \mid \leq  \, \sigma  } \, 
  {{\Delta t ^{\mid \delta \mid +\mid \varepsilon \mid }}
\over{\mid \! \delta \!  \mid ! }} \,   
 M_{k \ell} \, M^{-1}_{\ell p}  \, \Psi_{p r} \, P_{\ell \delta} \,  B^{\varepsilon}_{rj}
 \, \partial_{\delta + \varepsilon} W_j \,+\, 
 {\rm O}(\Delta t^{\sigma + 2}) \, .    
\end{array}  \right. \end{equation*} 
We set, with $ \, \,  \mid \! \gamma  \!  \mid \, = \, \sigma +1 \,, \,\,  
k \geq N  \,, \,\,   0 \leq  j \leq N-1 , \, $ 
\begin{equation*} \label{6.25} 
B^{\gamma}_{kj} \,= \displaystyle 
\, {{1}\over{s_k}} \,  \bigg( -  \sum_{1 \leq q \leq  \sigma + 1 } \, 
 {{1}\over{ q \, ! }} \,  D^{q,\gamma}_{kj} 
 \,+\, \sum_{1 \leq \mid \delta \mid \leq \sigma + 1,\,
0 \leq \mid \varepsilon \mid \leq  \, \sigma,\,  \delta +  \varepsilon = \gamma}
 \,    {{1} \over{\mid \! \delta \!  \mid ! }} \,  
 M_{k \ell} \, M^{-1}_{\ell p}  \, \Psi_{p r} \, P_{\ell \delta} \,  B^{\varepsilon}_{rj}
  \bigg) \, 
\end{equation*} 
and the relation   (\ref{2.25})  is established by induction. 
$\hfill \square$


\section*{Appendix B.\\  Notations for classical lattice Boltzmann   schemes  }

\noindent 
In order  to define precisely our results,  the numbering of degrees of freedom
must be defined and we make this point  precise  in this Appendix with the help of 
usual graphics.   The  choice of moments, {\it id
  est} the $M$ matrix (relation (\ref{1.8})) is also  made explicit.

\begin{figure}[!h]  
\centerline {\includegraphics[width=0.35 \textwidth ]{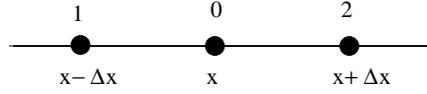}}
\smallskip
\caption {Stencil for the D1Q3 lattice Boltzmann scheme  } 
\label{fd1q3}
\end{figure} 

\bigskip \noindent $\bullet$ \quad 
{\bf   D1Q3 for advective thermics  }

\smallskip \noindent   
Recall first that the D1Q3  lattice Boltzmann scheme ($J=2$ in relation  (\ref{1.8})) 
uses three neighbours for a given node $x$:  the vertex 
$x$ itself and the first neighbours located at $\, \pm \Delta x \,$ from $x$ (see
Figure \ref{fd1q3}).  We introduce $ \, \lambda \,$ as in (\ref{1.23}) and adopt a labelling  for
matrix $M$ of relation  (\ref{1.8}) as in  Figure~ \ref{fd1q3}: 
\begin{equation} \label{7.1}  
M \,=\,  \begin{pmatrix}     1 & 1 & 1 \cr -\lambda & 0 &  \lambda \cr 
 \lambda^2 / 2 & 0 & \lambda^2 / 2  \end{pmatrix}   \, . \,  
\end{equation} 
%

\begin{figure}[!t]  
 \centerline 
{ {\includegraphics[width=.35  \textwidth] {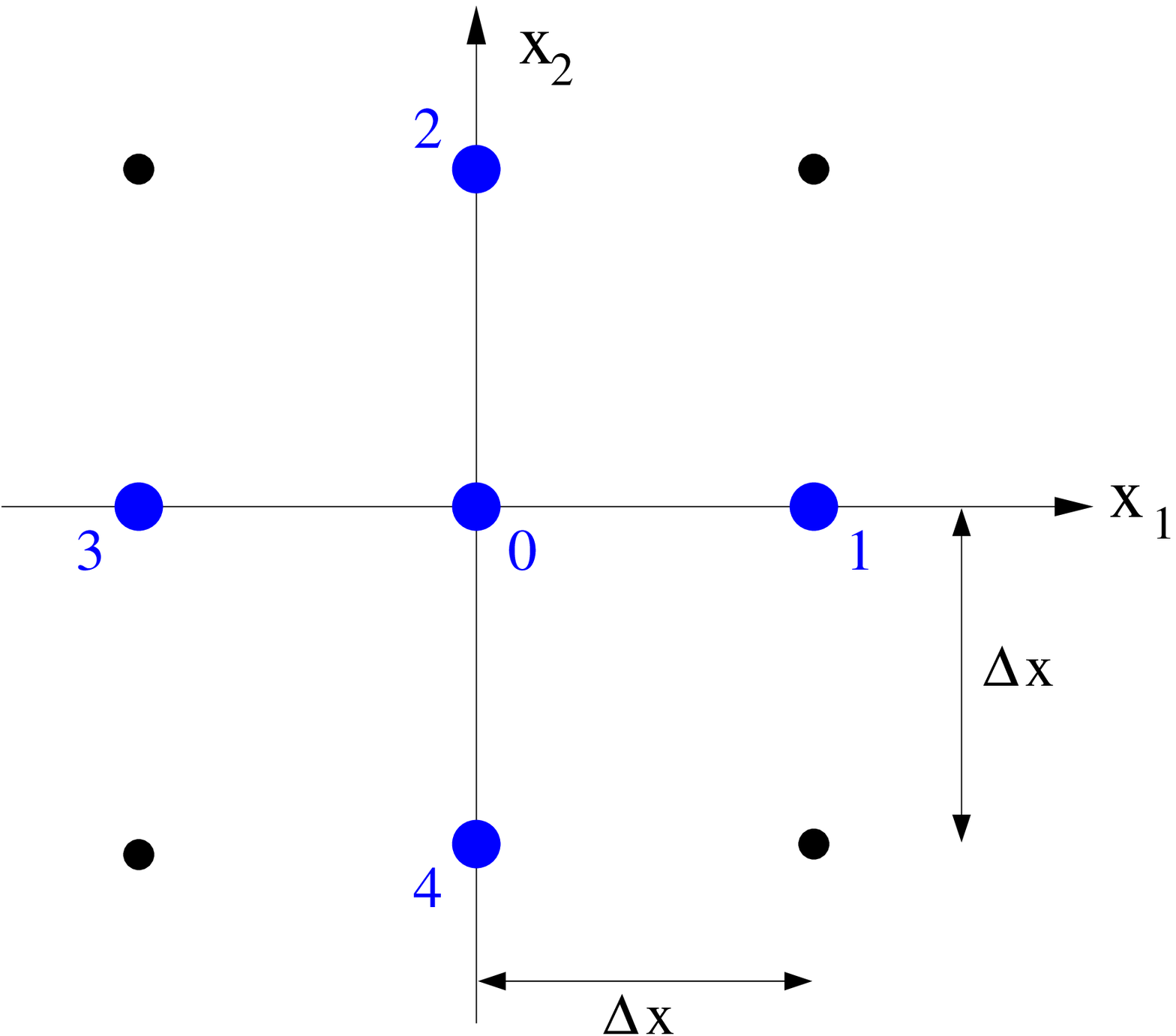}}  \qquad 
 {\includegraphics[width=.35   \textwidth]{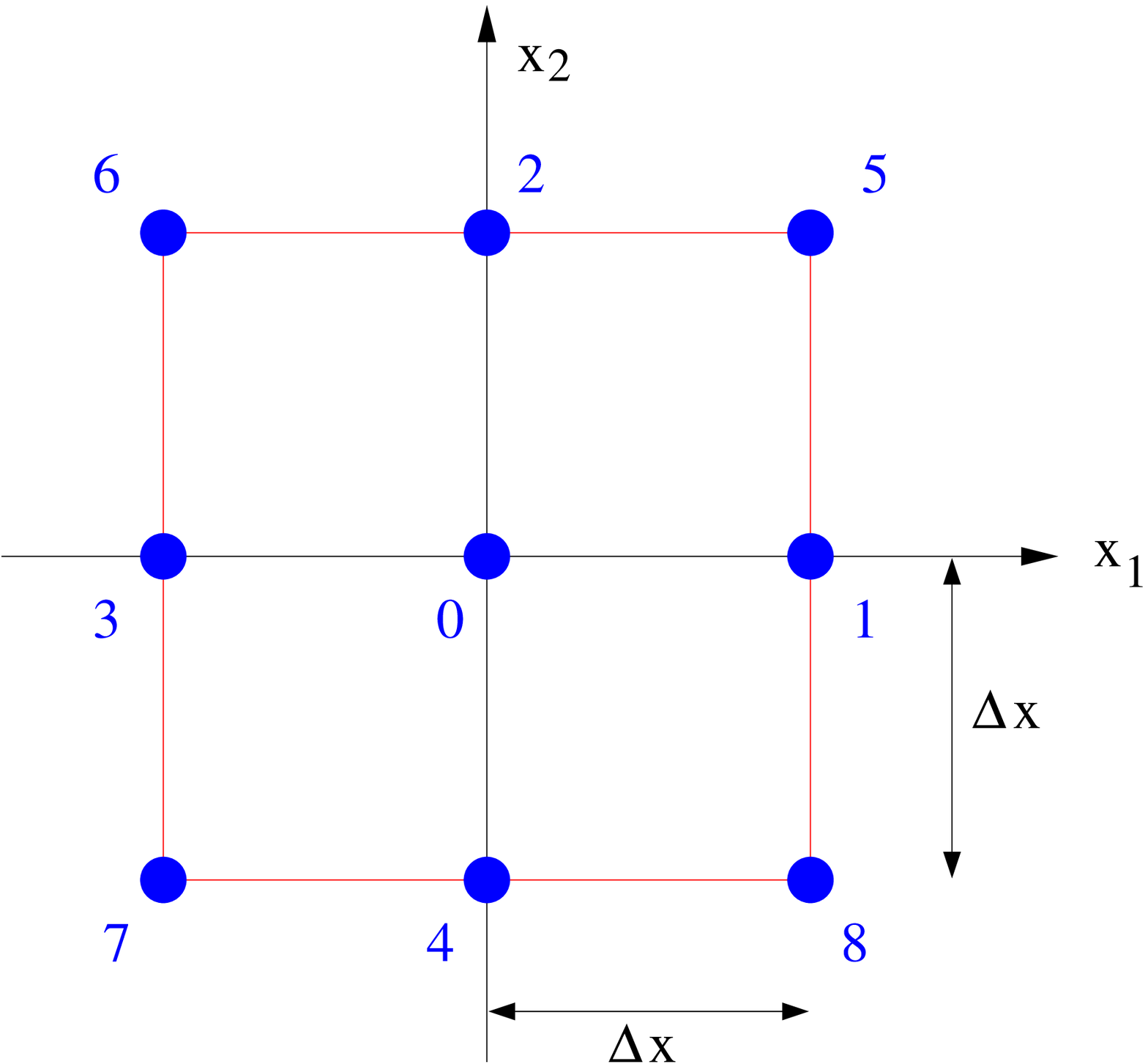}}}  
\smallskip
\caption {Stencils for  D2Q5 and D2Q9 lattice Boltzmann schemes  } 
\label{fd2q5q9}
\end{figure}    

\bigskip \noindent $\bullet$ \quad 
{\bf   D2Q5 for classical thermics  }
 
\smallskip \noindent
We have now four ($J=4$) nontrivial possible directions for propagation of particles 
(Figure \ref{fd2q5q9}, left). 
We adopt for the $M$ matrix of relation  (\ref{1.8}) the following choice:
\begin{equation}  \label{7.14}  
M   \,=\, \begin{pmatrix}   1 &  1 & 1 & 1 & 1  \cr 
                      0 &  \lambda & 0 & -\lambda & 0   \cr 
                      0 &  0 & \lambda & 0 & -\lambda  \cr 
                      -4 &  1 & 1 & 1 & 1  \cr  
                       0 &   1 & -1 & 1 & -1  \end{pmatrix}   \, . 
\end{equation}   

\bigskip \noindent $\bullet$ \quad 
{\bf  D2Q9 for  classical thermics  }

\smallskip \noindent  
The lattice Boltzmann model  D2Q9 is obtained from the 
  D2Q5 model by adding four velocities 
along the diagonals (Figure \ref{fd2q5q9}, right). 
The evaluation of matrix $M$ is entirely nontrivial. We
refer the reader to \cite{LL00},   and the reader can also consult our introduction
\cite{Du07}.  We have:   
\begin{equation}  \label{7.19} 
M   \,=\, \begin{pmatrix} 
1 &  1 & 1 & 1 & 1 & 1 & 1 & 1 & 1 \cr 
 0 &  \lambda & 0 & -\lambda & 0 & \lambda & -\lambda & -\lambda  & \lambda \cr 
 0 &  0 & \lambda & 0 & -\lambda & \lambda & \lambda & -\lambda & -\lambda \cr 
 -4 &  -1 & -1 & -1 & -1 & 2 & 2 & 2 & 2 \cr 
 4 &  -2 & -2 & -2 & -2 & 1 & 1 & 1 & 1 \cr 
 0 &  -2 & 0 & 2 & 0 & 1 & -1 & -1 & 1 \cr 
 0 &  0 & -2 & 0 & 2 & 1 & 1 & -1 & -1 \cr 
 0 &  1 & -1 & 1 & -1 & 0 & 0 & 0 & 0 \cr 
 0 &  0 & 0 & 0 & 0 & 1 & -1 & 1 & -1  \end{pmatrix}  \, . 
\end{equation} 
%

\begin{figure}[!t]  
 \centerline 
{ {\includegraphics[width=.45  \textwidth] {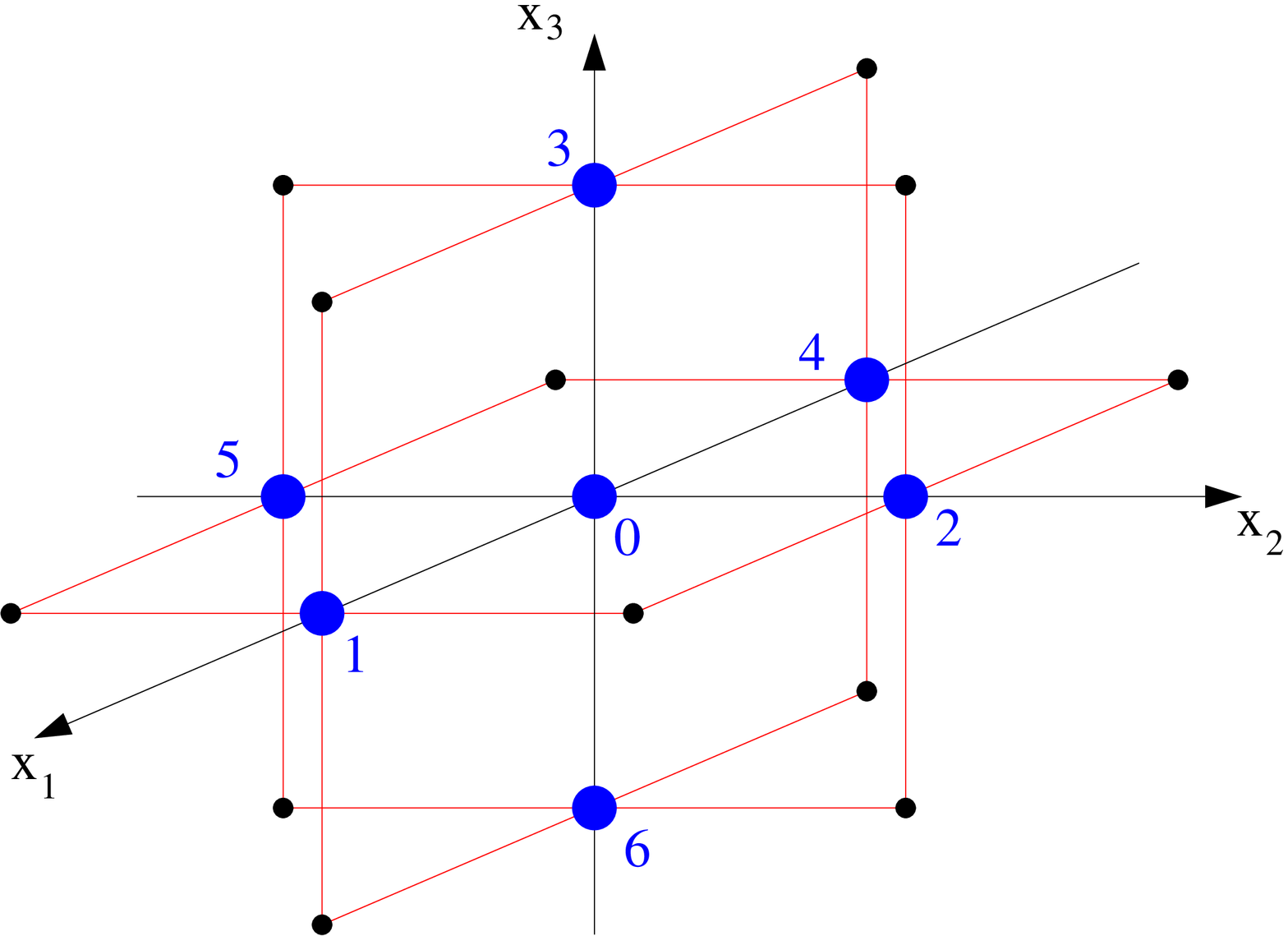}}  \qquad 
 {\includegraphics[width=.45   \textwidth]{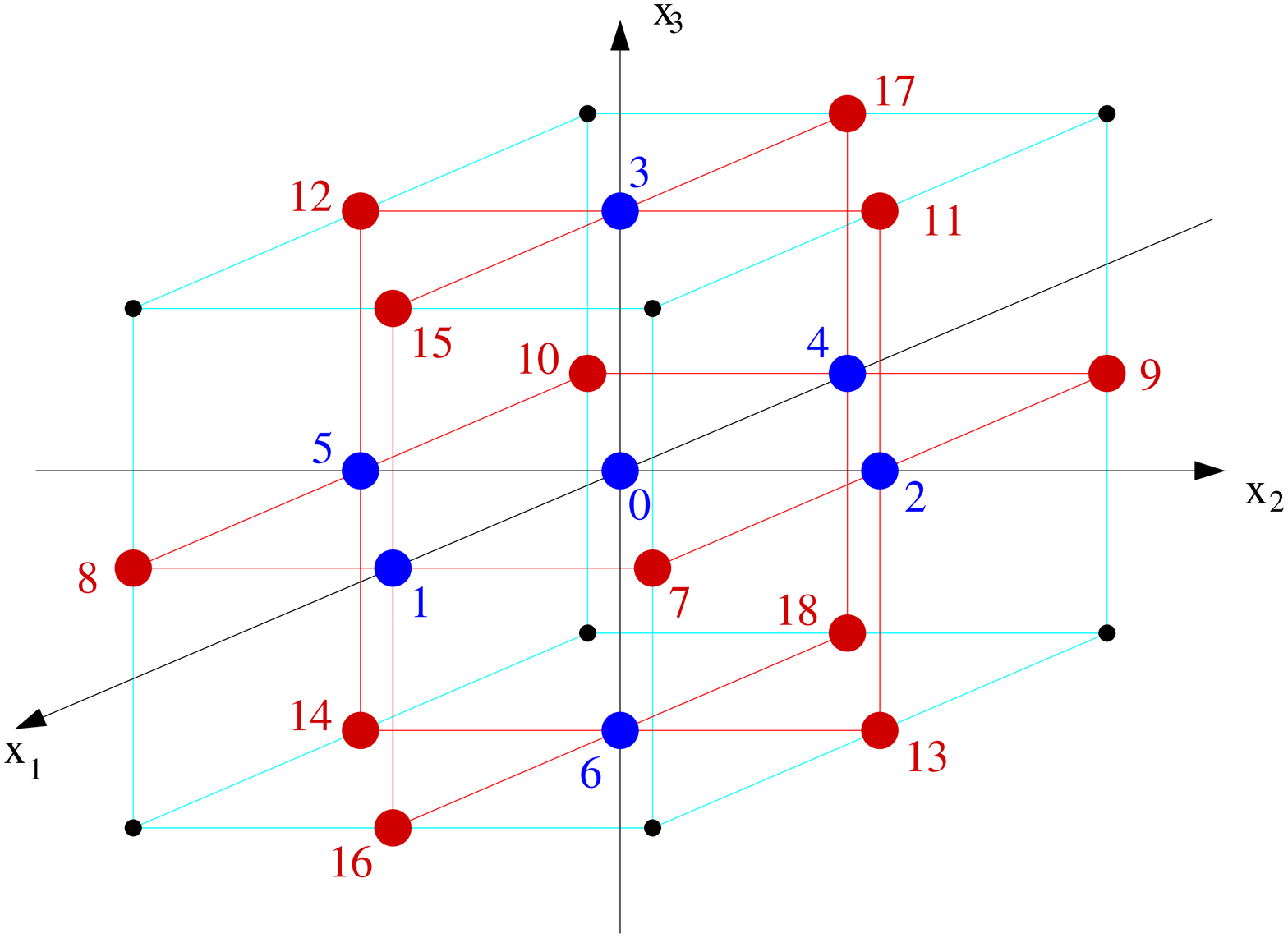}}}  
\smallskip
\caption {Stencils for  D3Q7  and D3Q19  lattice Boltzmann schemes  } 
\label{fd3q7q19}
\end{figure}    

\bigskip \noindent $\bullet$ \quad 
{\bf   D3Q7 for pure  thermics }

\noindent 
For three-dimensional thermics, one only needs  a seven point scheme and can use the so-called
D3Q7 lattice Boltzmann scheme whose stencil is described in 
the left part of  Figure \ref{fd3q7q19}. 
The matrix is not very difficult to construct. We follow \cite{LL03}: 
\begin{equation}  \label{7.36}  
 M   \,=\,  \begin {pmatrix}
1 &  1 & 1 & 1 & 1 & 1 & 1  \cr 
 0 &  \lambda & 0 & 0  & -\lambda &   0  & 0  \cr 
 0 &  0 & \lambda & 0 & 0 & -\lambda &   0  \cr 
 0 &  0 &  0 &   \lambda  &  0 & 0  & -\lambda   \cr 
 0  &  -1    &    -1  &  2 & -1  & -1            & 2   \cr      
 0  &   1    &    -1  &  0 & 1   &  -1       & 0        \cr     
 -6 &   1    &    1   &  1   &  1       & 1      &  1   \end {pmatrix}  \, . 
\end{equation} 
%
 
  \bigskip  \noindent $\bullet$ \quad 
{\bf  D3Q19 for linearized Navier--Stokes }

\noindent 
The D3Q19 Lattice Boltzmann scheme is described with details 
{\it e.g.} in J.~T\"olke {\it et al}  \cite{TKSR02} and the stencil 
is presented in  Figure \ref{fd3q7q19} (right).
The matrix $M$ that parameterizes the transformation (\ref{1.8}) looks
 like this: 
\begin{equation*}  \label{M-d3q19}  
 M   \,=\,   \left[ \begin {array} {ccccccccccccccccccc} 
 \scriptscriptstyle  1 \pet   1 \pet   1 \pet 1 \pet 1 \pet 1 \pet 1 \pet 1 
\pet  1 \pet 1 \pet 1 \pet 1 \pet  1 \pet 1 \pet  1 \pet  1 \pet 1 \pet 1 \pet 1   \\
 \scriptscriptstyle  0 \pet \lambda \pet 0 \pet 0 \pet -\lambda \pet 0 \pet 0 \pet  \lambda \pet  \lambda  \pet  -\lambda \pet -\lambda  
\pet 0 \pet 0 \pet 0 \pet 0 \pet    \lambda \pet  \lambda \pet   -\lambda \pet  -\lambda    \\  
 \scriptscriptstyle 0 \pet 0 \pet   \lambda \pet 0 \pet 0 \pet  -\lambda \pet 0  \pet \lambda \pet  -\lambda \pet  \lambda \pet  -\lambda \pet  \lambda \pet  -\lambda \pet 
 \lambda \pet  -\lambda \pet   0 \pet 0 \pet 0 \pet 0  \\  
 \scriptscriptstyle 0 \pet 0 \pet 0 \pet  \lambda \pet 0 \pet 0 \pet   -\lambda \pet 0 \pet 0 \pet 0 \pet 0 \pet  \lambda \pet  \lambda  \pet
-\lambda \pet -\lambda \pet    \lambda \pet  -\lambda \pet  \lambda \pet  -\lambda  \\  
 \scriptscriptstyle    -30 \lambda^2 \pet      -11 \lambda^2  \pet    -11  \lambda^2 \pet 
   -11 \lambda^2  \pet     -11 \lambda^2  \pet      -11 \lambda^2  \pet     -11 \lambda^2 
\pet     8 \lambda^2 \pet   8  \lambda^2 \pet    8  \lambda^2 \pet    8  \lambda^2
 \pet    8  \lambda^2 \pet    8  \lambda^2 \pet    8  \lambda^2 
\pet      8  \lambda^2  \pet     8  \lambda^2 \pet    8  \lambda^2 \pet     8     \lambda^2 
\pet 8  \lambda^2  \\ 
 \scriptscriptstyle 0 \pet 2 \lambda^2 \pet - \lambda^2 \pet - \lambda^2
\pet 2 \lambda^2 \pet - \lambda^2 \pet - \lambda^2 
 \pet  \lambda^2 \pet  \lambda^2 \pet  \lambda^2 \pet  \lambda^2 
 \pet  -2 \lambda^2 \pet  -2 \lambda^2 \pet  -2 \lambda^2 \pet  -2 \lambda^2
 \pet  \lambda^2 \pet  \lambda^2 \pet  \lambda^2 \pet  \lambda^2  \\ 
 \scriptscriptstyle 0 \pet  0 \pet  \lambda^2 \pet  - \lambda^2 \pet  
0 \pet  \lambda^2 \pet  - \lambda^2 
 \pet  \lambda^2 \pet  \lambda^2 \pet  \lambda^2 \pet  \lambda^2  
 \pet 0  \pet 0  \pet 0   \pet 0 
 \pet  \lambda^2 \pet  \lambda^2 \pet  \lambda^2 \pet  \lambda^2 \\ 
 \scriptscriptstyle 0 \pet  0 \pet  0 \pet  0 \pet  0 \pet  0 \pet  0 \pet 
 \lambda^2 \pet - \lambda^2 \pet - \lambda^2 \pet  \lambda^2 \pet  
 0 \pet  0 \pet  0 \pet  0 \pet  0 \pet  0 \pet  0 \pet 0 \\ 
 \scriptscriptstyle 0 \pet  0 \pet  0 \pet  0 \pet  0 \pet  0 \pet  0 \pet  
0 \pet  0 \pet  0 \pet  0 \pet  
 \lambda^2 \pet - \lambda^2 \pet - \lambda^2 \pet  \lambda^2 \pet  
 0 \pet  0 \pet  0 \pet  0  \\  
 \scriptscriptstyle 0 \pet  0 \pet  0 \pet  0 \pet  0 \pet  0 \pet  0 \pet  
0 \pet  0 \pet  0 \pet  0 \pet 0 \pet  0 \pet  0 \pet  0 \pet  
 \lambda^2 \pet - \lambda^2 \pet - \lambda^2 \pet  \lambda^2  \\ 
 \scriptscriptstyle 0 \pet  -4 \lambda^3  \pet  0 \pet  0 \pet  4 \lambda^3  \pet  0 \pet  0 \pet
\lambda^3  \pet \lambda^3  \pet -\lambda^3  \pet -\lambda^3  \pet  
 0 \pet  0 \pet  0 \pet  0 \pet
\lambda^3  \pet \lambda^3  \pet -\lambda^3  \pet -\lambda^3   \\ 
 \scriptscriptstyle 0 \pet  0 \pet   -4 \lambda^3  \pet 0 \pet 0 \pet 4 \lambda^3 \pet 0  \pet 
\lambda^3  \pet - \lambda^3  \pet \lambda^3  \pet - \lambda^3  \pet 
\lambda^3  \pet - \lambda^3  \pet \lambda^3  \pet - \lambda^3  \pet   
 0 \pet  0 \pet  0 \pet  0     \\ 
\scriptscriptstyle 0 \pet  0 \pet  0 \pet  -4 \lambda^3  \pet 0 \pet 0 \pet 4 \lambda^3  \pet 
 0 \pet  0 \pet  0 \pet  0 \pet 
\lambda^3  \pet \lambda^3  \pet -\lambda^3  \pet -\lambda^3  \pet 
\lambda^3  \pet -\lambda^3  \pet \lambda^3  \pet -\lambda^3  \\ 
\scriptscriptstyle 12  \lambda^4 \pet -4  \lambda^4 \pet -4  \lambda^4 \pet -4  \lambda^4
 \pet -4  \lambda^4 \pet -4  \lambda^4 \pet -4  \lambda^4 
\pet \lambda^4 \pet \lambda^4 \pet \lambda^4 \pet \lambda^4 
\pet \lambda^4 \pet \lambda^4 \pet \lambda^4 \pet \lambda^4 
\pet \lambda^4 \pet \lambda^4 \pet \lambda^4 \pet \lambda^4 \\ 
\scriptscriptstyle 0  \pet -4  \lambda^4   \pet 2  \lambda^4  \pet 2  \lambda^4 
\pet -4  \lambda^4   \pet 2  \lambda^4  \pet 2  \lambda^4 
\pet \lambda^4 \pet \lambda^4 \pet \lambda^4 \pet \lambda^4 
\pet -2 \lambda^4 \pet -2 \lambda^4 \pet -2 \lambda^4 \pet -2 \lambda^4 
\pet \lambda^4 \pet \lambda^4 \pet \lambda^4 \pet \lambda^4  \\ 
\scriptscriptstyle 0 \pet 0  \pet -2 \lambda^4  \pet 2 \lambda^4 
\pet 0  \pet -2 \lambda^4  \pet 2 \lambda^4
\pet \lambda^4 \pet \lambda^4 \pet \lambda^4 \pet \lambda^4 
\pet 0 \pet 0 \pet 0 \pet 0 
\pet -\lambda^4 \pet -\lambda^4 \pet -\lambda^4 \pet -\lambda^4  \\ 
\scriptscriptstyle 0 \pet 0 \pet 0 \pet 0 \pet 0 \pet 0 \pet  0 
\pet  \lambda^3  \pet  \lambda^3 \pet - \lambda^3  \pet  -\lambda^3  \pet 
 0 \pet 0 \pet 0 \pet 0 
\pet  -\lambda^3  \pet  -\lambda^3 \pet  \lambda^3  \pet  \lambda^3 \\ 
\scriptscriptstyle 0 \pet 0 \pet 0 \pet 0 \pet 0 \pet 0 \pet  0 
\pet  -\lambda^3  \pet  \lambda^3 \pet - \lambda^3  \pet  \lambda^3  \pet  
\lambda^3 \pet - \lambda^3  \pet  \lambda^3 \pet - \lambda^3  \pet 0 \pet 0 \pet 0 \pet 0  \\ 
\scriptscriptstyle 0 \pet 0 \pet 0 \pet 0 \pet 0 \pet 0 \pet  0 \pet 0 \pet 0 \pet 0 \pet 0
\pet  -\lambda^3  \pet - \lambda^3 \pet \lambda^3  \pet \lambda^3 
\pet  \lambda^3  \pet - \lambda^3 \pet \lambda^3  \pet -\lambda^3 
 \end {array}  \right] \, .  \end{equation*}  
Due to the important number of moments, we detail in this sub-section the way the previous
matrix is obtained. First, velocities $\, v_j^{\alpha} \, $ for 
$\, 0 \leq j \leq J\equiv 18 $ and  $\, 1 \leq \alpha \leq 3 \, $ 
are naturally associated with Figure \ref{fd3q7q19}.  
 The four first moments $\, \rho \, $ and $ \, q^{\alpha}\,$ are determined according to 
(\ref{1.10}) and (\ref{1.11}) and the associated elements for matrix $\, M \,$ are given
 in (\ref{1.12}) and (\ref{1.13}).  
The construction of other moments uses the  tensorial nature of the 
variety of moments that can be constructed, as analyzed by Rubinstein and Luo 
\cite{RL08}: scalar fields are naturally coupled with one another  and similarly  for 
vector fields,  and so on. So components of kinetic energy are introduced: 
\begin{equation} \label{3-kinetic}  
\widetilde {M} _{4 j} \,=\,  19 \,   \sum_{\alpha}  \mid  v_j^{\alpha}  \mid ^2  
\, , \qquad 0 \leq j \leq J \, . 
\end{equation} 
The entire set of second-order tensors is completed according to 
\begin{equation} \label{3-ordre2}  \left\{ \begin{array}{rcl}  
\widetilde {M} _{5 j} &  \,=\, &  2 \,  (v_j^{1})^2 \,-\,  (v_j^{2})^2\,-\,  (v_j^{3})^2
\\ \widetilde {M} _{6 j} &  \,=\, &  (v_j^{2})^2\,-\,  (v_j^{3})^2    
\\  \widetilde {M} _{7 j} &  \,=\,&   v_j^{1}  \,  v_j^{2}   \,, \quad  
  \widetilde {M} _{8 j}   \,=\,   v_j^{2}  \,  v_j^{3}   \,, \quad    
     \widetilde {M} _{9 j}    \,=\,    v_j^{3}  \,  v_j^{1}  
\, , \qquad  \qquad  0 \leq j \leq J \, .  
\end{array} \right.    \end{equation} 
The three components of heat flux are  defined by  
 
\begin{equation} \label{3-flux-chaleur}  \left\{ \begin{array}{c}  \displaystyle 
\widetilde {M} _{10 \, j}  =
 \displaystyle  5 \,   v_j^{1} \, \sum_{\alpha}  \mid  v_j^{\alpha}  \mid ^2  \,,
\quad   \widetilde {M} _{11 \,   j}   =
\displaystyle   5 \,    v_j^{2} \, \sum_{\alpha}  \mid  v_j^{\alpha}  \mid ^2   \,,
\quad    \widetilde {M} _{12 \,   j}  =
\displaystyle  5 \,    v_j^{3} \, \sum_{\alpha}  \mid  v_j^{\alpha}  \mid ^2  \,, 
\\ \hfill  
 \qquad  \qquad  \qquad  \qquad \qquad  \qquad  0 \leq j \leq J \, .  
\end{array} \right.    \end{equation}  
We finally obtain the moments of higher degree: the square of the kinetic energy
\begin{equation} \label{3-energy-square} 
  \widetilde {M} _{13 \,   j}   \,=\,  
 \frac{21}{2} \,  \Big(    \sum_{\alpha}  \mid  v_j^{\alpha}  \mid ^2 \Big)^2    
\, , \qquad 0 \leq j \leq J \, , 
\end{equation}  
second-order moments ``weighted'' by  kinetic energy: 
\begin{equation} \label{3-xxe-yye}  \left\{ \begin{array}{rcl}  \displaystyle 
\widetilde {M} _{14 \, j} &  \,=\, &  \displaystyle   
  3  \, \Big(  2 \,  (v_j^{1})^2 \,-\,  (v_j^{2})^2\,-\,  (v_j^{3})^2 \Big) \, \,  
  \sum_{\alpha}  \mid  v_j^{\alpha}  \mid ^2 
\\   \widetilde {M} _{15 \,   j} &  \,=\, &   \displaystyle   
  3   \, \big(  (v_j^{2})^2\,-\,  (v_j^{3})^2  \big) \, 
  \sum_{\alpha}  \mid  v_j^{\alpha}  \mid ^2  
\, , \qquad  \qquad  0 \leq j \leq J \, ,  
\end{array} \right.    \end{equation}  
and third-order anti-symmetric moments:  
\begin{equation} \label{3-third}  \left\{ \begin{array}{rcl}  \displaystyle 
\widetilde {M} _{16 \, j} &  \,=\, &  \displaystyle   
 v_j^{1} \, \big(  (v_j^{2})^2  \,-\,  (v_j^{3})^2 \big) 
\\   \widetilde {M} _{17 \,   j} &  \,=\, &   \displaystyle  
 v_j^{2} \, \big(  (v_j^{3})^2  \,-\,  (v_j^{1})^2 \big)  
\\   \widetilde {M} _{18 \,   j} &  \,=\, &   \displaystyle  
 v_j^{3} \, \big(  (v_j^{1})^2  \,-\,  (v_j^{2})^2 \big)  
\, , \qquad   0 \leq j \leq J \, .  
\end{array} \right.    \end{equation}  
Then matrix $\, M \,$ is orthogonalized from relations 
 (\ref{1.12}),  (\ref{1.13}),  (\ref{3-kinetic}), (\ref{3-ordre2}),
 (\ref{3-flux-chaleur}), (\ref{3-energy-square}), (\ref{3-xxe-yye}) and  
 (\ref{3-third}) with a Gram-Schmidt classical algorithm: 
  \begin{equation*} \label{gram-schmidt}  
 M _{i j}  \,=\,   \widetilde {M} _{i j} - \sum_{\ell < i}  g_{i \ell} 
\, M_{\ell j}  \,,  \quad  i \geq 4    \, .  
 \end{equation*} 
The coefficients $\,  g_{i \ell} \,$ 
 are computed recursively in order to force orthogonality: 
\begin{equation*} \label{orthogo}  
\sum_{j = 0}^J \, M _{i j} \, M _{k j} \, = \, 0 \, \qquad   \textrm {for } i \not = k  \, .
\end{equation*} 
%


\section*{Appendix C. Quartic parameters in three dimensions }

\smallskip \noindent $\bullet$ \quad 
We use the equivalent equations of lattice Boltzmann scheme D3Q19 obtained previously in
the following way. We consider the vector of conserved variables (\ref{1.9}): 
%
$ \,  W  \equiv   ( \rho ,\, q_x ,\, q_y,\, q_z )^{\rm  \displaystyle t}  . $  
%
We write the equivalent partial differential equations    under the combined   form:
 \begin{equation}   \label{3d-eq-equiv}
\partial_t W_k \,+\, \sum_{ j,\, p,\,  q ,\,  r} \,  A_{k p q r}^j \, 
  \partial_{x}^p  \partial_{y}^q \partial_{y}^r  W_j   \,=\, {\rm O}(\Delta t^4) \,. 
\end{equation} 
We search dissipative mode  solution of  (\ref{3d-eq-equiv})  under the form
%
$ \,   W(t)=   {\rm e}^{   - \Gamma t \,+\, i (k_x \, x 
+ k_y\, y + k_z\, z)  } \,   \widetilde{W} \, . $ 
%
Then $\Gamma$ is an eigenvalue of the matrix $A$ defined by  
 \begin{equation*}   \label{d3q19-det-matrice}
A_{k}^j  \,=\,      \sum_{ p,\,  q,\, r  } \, 
 A_{k p q}^j \,\, \, (i\, k_x)^p \,(i\, k_y)^q  \,(i\, k_z)^r  \, . 
\end{equation*}

\bigskip \noindent $\bullet$ \quad 
We wish to solve this dispersion equation with a high order of accuracy, {\it id est} in
our present case: 
 \begin{equation}     \label{d3q19-det}  
 \Delta   \, \equiv\,{\rm det} \, \left[ A \,-\,  \Gamma   \, {\rm Id} \right] \,=\,  
  {\rm O}(\Delta t^7)  \, . 
\end{equation} 
We impose also that this eigenvalue is {\bf double} as classical for shear waves in three
dimensions \cite{LL59}: 
%
$\, \frac {\rm d}{{\rm d}\Gamma}  \left( {\rm det} \, 
\left[ A \,-\,  \Gamma   \, {\rm Id} \right] \right) \,\approx\,   0  \, .  \, $ 
%
The first nontrivial term in powers of $ \, \Delta t \,$ for this derivative of the
determinant is the term of order 3. Then we force  
 \begin{equation}     \label{d3q19-derdet}     
\frac {\rm d}{{\rm d}\Gamma}  \left( {\rm det} \, 
\left[ A \,-\,  \Gamma   \, {\rm Id}   \right] \right)  \,=\,     {\rm O}(\Delta t^4)  \, . 
\end{equation} 
For  Stokes problem  (incompressible shear modes) and D3Q19 lattice Boltzmann 
d'Humi\`eres scheme,  we have \cite{QDL92}:  
 \begin{equation}      \label{d3q19-Gamma}     
\Gamma \,\, \equiv \,\,  \nu  \, \mid k \mid ^2 \,\,=\,\, 
  {{ \lambda^2}\over{3}} \, \Delta t \, \sigma_5   
\,  \big( k_x^2 \,+\, k_y^2 \,+\, k_z^2 \big)  \, . 
\end{equation}

\bigskip \noindent $\bullet$ \quad 
We solve the set (\ref{d3q19-det})   (\ref{d3q19-derdet})   (\ref{d3q19-Gamma})   of
equations for all values of  the time step $ \, \Delta t \, $.  
We obtain in this way a set of eight  algebraic equations: 
 \begin{equation*}   \left\{ \begin{array}{rcl}
 2 \, \sigma_5 \, \sigma_{10}-4 \, \sigma_5^2+6 \, \sigma_5 \, \sigma_{16} &\,=\, & 1 \\ 
80 \, \sigma_5^4-32 \, \sigma_5^3 \, \sigma_{10}+24 \, \sigma_5^2 \, \sigma_{10} \, \sigma_{16}
+12 \,  \sigma_{14} \, \sigma_{16} \, \sigma_5^2
-8 \, \sigma_5^2 -4 \, \sigma_5^2 \, \sigma_{10}^2 &&\\
+12 \, \sigma_5^2 \, \sigma_{16}^2  
-12 \, \sigma_5^2 \, \sigma_{14} \, \sigma_{10}  
-12 \, \sigma_5 \, \sigma_{16} \, \sigma_{14} \, \sigma_{10} 
+6 \, \sigma_5 \, \sigma_{14} \, \sigma_{10}^2 &&\\
-8 \, \sigma_5 \, \sigma_{16}+6 \, \sigma_5 \, \sigma_{16}^2
 \, \sigma_{14} -\sigma_{14} \, \sigma_{16}+\sigma_{14} \, \sigma_{10}+1 &\,=\,& 0  \\ 
-48 \, \sigma_5^5 \, \sigma_{10}+44 \, \sigma_5^4 \, \sigma_{10}^2
+2000 \, \sigma_5^5 \, \sigma_{16}
+95 \, \sigma_5^2-16 \, \sigma_5^4 \, \sigma_{14} \, \sigma_{10} &&\\
+292 \, \sigma_{14} \, \sigma_{16} \, \sigma_5^2
+68 \, \sigma_5^2 \, \sigma_{14} \, \sigma_{10}
-272 \, \sigma_5^4 \, \sigma_{16} \, \sigma_{14}
-1032 \, \sigma_5^3 \,  \sigma_{16}^2 \, \sigma_{14}  &&\\
+56 \, \sigma_5^3 \, \sigma_{14} \, \sigma_{10}^2
-320 \, \sigma_5^6-1048 \, \sigma_5^4 \, \sigma_{10} \, \sigma_{16}
+\sigma_{14}^2+60 \, \sigma_5^2 \, \sigma_{16}^2 \, \sigma_{14}^2
-16 \, \sigma_5 \, \sigma_{16} \, \sigma_{14}^2       &&\\
+72 \, \sigma_5^2 \, \sigma_{14}^2 \, \sigma_{10} \, \sigma_{16}
-8 \, \sigma_5 \, \sigma_{14}^2 \, \sigma_{10}+24 \, \sigma_5^3 \, \sigma_{14}
+12 \, \sigma_5^2 \, \sigma_{14}^2 \, \sigma_{10}^2
-248 \, \sigma_5^4  &&\\
-464 \, \sigma_5^3 \, \sigma_{16} \, \sigma_{14} \, \sigma_{10}+148 \, \sigma_5^3 \, \sigma_{10}
-1284 \, \sigma_{16} \, \sigma_5^3+4284 \, \sigma_{16}^2 \, \sigma_5^4-20 \, \sigma_5 \, 
\sigma_{14}  &\,=\,& 0  
  \end{array}  \right. \end{equation*} 
  \begin{equation*}   \left\{ \begin{array}{rcl}
\big( -1 +  2 \, \sigma_5 \, \sigma_{10}-4 \, \sigma_5^2+6 \, \sigma_5 \, \sigma_{16} 
\big) \, 
\big( 2 \, \sigma_5 \, \sigma_{10} + 2 \,   \sigma_{14} \, \sigma_{10}    &&\\
- 2 \,  \sigma_5^2  - 10 \,  \sigma_5 \, \sigma_{16} 
-2 \,  \sigma_{14} \, \sigma_{16}  + 3 \big) \,    &\,=\,& 0  \\ 
96 \, \sigma_5^5 \, \sigma_{10}+24 \, \sigma_5^4 \, \sigma_{10}^2
-1920 \, \sigma_5^5 \, \sigma_{16}
+ 98 \, \sigma_5^2+24 \, \sigma_5^4 \, \sigma_{14} \, \sigma_{10}
+350 \, \sigma_{14} \, \sigma_{16} \, \sigma_5^2 &&\\
+34 \, \sigma_5^2 \, \sigma_{14} \, \sigma_{10}
+264 \, \sigma_5^4 \, \sigma_{16} \, \sigma_{14}
-1524 \, \sigma_5^3 \, \sigma_{16}^2 \, \sigma_{14}
+12 \, \sigma_5^3 \, \sigma_{14} \, \sigma_{10}^2 &&\\
+240 \, \sigma_5^6-576 \, \sigma_5 ^4 \, \sigma_{10} \, \sigma_{16}
+\sigma_{14}^2+102 \, \sigma_5^2 \, \sigma_{16}^2 \, \sigma_{14}^2
-20 \, \sigma_5 \, \sigma_{16} \, \sigma_{14}^2
+36 \, \sigma_5^2 \, \sigma_{14}^2 \, \sigma_{10} \, \sigma_{16} &&\\  
-4 \, \sigma_5 \, \sigma_{14}^2 \, \sigma_{10}
-24 \, \sigma_5^3 \, \sigma_{14}+6 \, \sigma_5^2 \, \sigma_{14}^2 \, \sigma_{10}^2
+240 \, \sigma_5^4-216 \, \sigma_5^3 \, \sigma_{16} \, \sigma_{14} \, \sigma_{10}
+72 \, \sigma_5^3 \, \sigma_{10} &&\\
- 1488 \, \sigma_{16} \, \sigma_5^3
+5688 \, \sigma_{16}^2 \, \sigma_5^4-20 \, \sigma_5 \, \sigma_{14}  &\,=\,& 0  \\ 
 -\sigma_5+6 \, \sigma_{16} \, \sigma_5^2+2 \, \sigma_5^2 \, \sigma_{10}-4 \, \sigma_5^3 &\,=\,& 0  \\ 
2 \, \sigma_5^2 \, \sigma_{10}-2 \, \sigma_{16} \, \sigma_5^2+\sigma_5-\sigma_5 \, \sigma_{16} \, 
\sigma_{14}+\sigma_5 \, \sigma_{14} \, \sigma_{10}-12 \, \sigma_5^3   &\,=\,& 0  \\ 
10 \, \sigma_5 \, \sigma_{16} \, \sigma_{14}
+2 \, \sigma_5 \, \sigma_{14} \, \sigma_{10}
+11 \, \sigma_5-\sigma_{14}+8 \, \sigma_5^3
-82 \, \sigma_{16} \, \sigma_5^2+6 \, \sigma_5^2 \, \sigma_{10} &\,=\,& 0 \, .  \\ 
 \end{array}  \right. \end{equation*}
These equations have  only 
one  nontrivial family of solutions  given by   (\ref{d3q19-quartic}).  
%


\medskip

\end{document}